\documentclass[leqno]{amsart}

\usepackage[foot]{amsaddr}

\usepackage{color}
\usepackage{geometry}
\usepackage{subfiles}
\usepackage{graphicx}
\usepackage[section]{placeins}
\usepackage{hyperref}
\usepackage{stmaryrd}
\usepackage{amsmath,amssymb,amsfonts,amsthm,textcomp}
\usepackage{mathtools}
\usepackage{tensor}
\usepackage{multicol}
\usepackage{subfig}
\usepackage{csquotes}
\usepackage{stackrel}
\usepackage{tikz}
\usepackage{siunitx}
\usepackage{soul}
\usepackage{pdfpages}
\usepackage{dsfont}
\usepackage{booktabs}
\graphicspath{
              {Images/}
             }

\newfont{\tenbfsl}{cmbxti9 scaled 1200}
\newfont{\tenbbb}{msbm10}
\newfont{\svnbbb}{msbm8}

\DeclareMathOperator*{\argmin}{arg\,min}

\newcommand{\bs}[1]{\boldsymbol{#1}}
\newcommand{\cl}[1]{\mathcal{#1}}

\newcommand{\bb}[1]{\mathbb{#1}}

\newcommand{\di}[1]{\,\mathrm{d}{#1}}

\newcommand{\trans}{\scriptscriptstyle\mskip-1mu\top\mskip-2mu}

\newtheorem{cor}{Corollary}
\newtheorem{lem}{Lemma}
\newtheorem{prop}{Proposition}
\newtheorem{rmk}{Remark}
\newtheorem{set}{Definition}

\newtheorem{asu}{Assumption}

\begin{document}

\title[Multilevel randomized quasi-Monte Carlo for nested integration]{Multilevel randomized quasi-Monte Carlo estimator for nested integration}
\author{Arved Bartuska$^{\#,1,2}$, Andr\'{e} Gustavo Carlon$^1$, Luis Espath$^3$, Sebastian Krumscheid$^5$, \& Ra\'{u}l Tempone$^{1,2,4}$}
\address{$^1$Department of Mathematics, RWTH Aachen University, Geb\"{a}ude-1953 1.OG, Pontdriesch 14-16, 161, 52062 Aachen, Germany}
\address{$^2$King Abdullah University of Science \& Technology (KAUST), Computer, Electrical and Mathematical Sciences \& Engineering Division (CEMSE), Thuwal 23955-6900, Saudi Arabia}
\address{$^3$School of Mathematical Sciences, University of Nottingham, Nottingham, NG7 2RD, United Kingdom}
\address{$^4$Alexander von Humboldt Professor in Mathematics for Uncertainty Quantification, RWTH Aachen University, Germany}
\address{$^5$Scientific Computing Center, and Institute for Applied and Numerical Mathematics, Karlsruhe Institute of Technology, Germany}
\email{$^\#$arved.bartuska@gmail.com}

\date{\today}

\begin{abstract}
\noindent
Nested integration problems arise in various scientific and engineering applications, including Bayesian experimental design, financial risk assessment, and uncertainty quantification. These nested integrals take the form $\int f\left(\int g(\bs{y},\bs{x})\di{}\bs{x}\right)\di{}\bs{y}$, for nonlinear $f$, making them computationally challenging, particularly in high-dimensional settings. Although widely used for single integrals, traditional Monte Carlo (MC) methods can be inefficient when encountering complexities of nested integration. This work introduces a novel multilevel estimator, combining deterministic and randomized quasi-MC (rQMC) methods to handle nested integration problems efficiently. In this context, the inner number of samples and the discretization accuracy of the inner integrand evaluation constitute the level. We provide a comprehensive theoretical analysis of the estimator, deriving error bounds demonstrating significant reductions in bias and variance compared with standard methods. The proposed estimator is particularly effective in scenarios where the integrand is evaluated approximately, as it adapts to different levels of resolution without compromising precision.
We verify the performance of our method via numerical experiments, focusing on estimating the expected information gain of experiments. When applied to Gaussian noise in the experiment, a truncation scheme ensures finite error bounds. The results reveal that the proposed multilevel rQMC estimator outperforms existing MC and rQMC approaches, offering a substantial reduction in computational costs and offering a powerful tool for practitioners dealing with complex, nested integration problems across various domains.
$\cdot$
65C05 
$\cdot$
65D30 
$\cdot$
65D32 
\end{abstract}

\maketitle
\section*{Acknowledgments}

This publication is based upon work supported by the King Abdullah University of Science and Technology (KAUST) Office of Sponsored Research (OSR) under Award No.~OSR-2019-CRG8-4033, the Alexander von Humboldt Foundation, the Deutsche Forschungsgemeinschaft (DFG, German Research Foundation) -- 333849990/GRK2379 (IRTG Hierarchical and Hybrid Approaches in Modern Inverse Problems), and was partially supported by the Flexible Interdisciplinary Research Collaboration Fund at the University of Nottingham Project ID 7466664. The authors thank Prof.~Abdul-Lateef Haji-Ali for his helpful and constructive comments.

\tableofcontents                        


\section{Introduction}
In this work, we address the computational challenges of nested integration problems of the form \\
$\int f\left(\int g(\bs{y},\bs{x})\di{}\bs{x}\right)\di{}\bs{y}$ where $f$ is a nonlinear function. Such problems frequently arise in engineering, mathematical finance, and medical decision-making. Traditional quadrature methods are constrained by the curse of dimensionality, making them impractical as the dimension increases. In contrast, Monte Carlo (MC) methods offer a more scalable solution but may become less efficient when managing complex nested structures. Our goal is to improve the efficiency of these estimators by developing particularly tailored, hierarchical quasi-MC (QMC) techniques.

The MC estimators are widely used for standard (i.e., nonnested) integration problems due to their scalability. However, achieving a desired error tolerance $TOL > 0$ with standard MC comes with a computational cost of $\mathcal{O}(TOL^{-2})$. The QMC methods, utilizing low-discrepancy sequences, can provide improved convergence for smooth integrands, reducing computational effort for instance to almost $\mathcal{O}(TOL^{-1})$, where the multiplicative constant usually depends exponentially on the dimension of the integration domain. Dimension-independent constant terms are possible for certain integrands when weighted function spaces are employed \cite{Dic13}. Moreover, higher-order QMC methods permit even faster convergence for certain integrands \cite{Dic13}.
Despite these improvements, challenges remain, especially when dealing with singularities. Indeed, standard error bounds for QMC methods via the Koksma--Hlawka inequality \cite{Hla61} become obsolete for singular integrands that exhibit infinite Hardy--Krause variation \cite{Caf98}. A useful error bound for such integrands was introduced by Owen \cite{Owe06} and later generalized by He et al.~\cite{He23}.

However, randomized QMC (rQMC) methods combine the efficiency of QMC with important properties of MC methods, including its unbiasedness and practical, variance-based error estimates. Moreover, the computational cost of such methods can be as low as almost $\mathcal{O}(TOL^{-2/3})$ for highly regular integrands. Nevertheless, for nested integration, the nonlinear function $f$ that evaluates the inner integrand introduces a bias in the outer integral estimation, even if the inner integral is estimated using an unbiased method. Thus, the number of inner and outer samples of numerical estimators must be carefully controlled to reduce the overall error. 

Ryan~\cite{Rya03} introduced the double-loop MC (DLMC) estimator for the expected information gain (EIG) of an experiment, which is a comprehensive tool to evaluate the quality of an experimental setup, or design \cite{Cha95, Lin56}. Estimating the EIG efficiently and accurately enables design optimization in a Bayesian setting, and hence efficient data collection in real-world applications. We consider integrands whose approximate evaluation requires a numerical solution to a partial differential equation (PDE), for instance, via the finite element method (FEM). The DLMC estimator was further developed in~\cite{Bec18} with regard to additional bias error stemming from the FEM approximation of the experiment model. The computational cost of this estimator to achieve an error tolerance $TOL$ is $\cl{O}(TOL^{-3})$ with no finite element approximation, and is $\cl{O}(TOL^{-3-\gamma/\eta_{\rm{w}}})$ with such an approximation, where the term $\gamma/\eta_{\rm{w}}$ relates to computational work and accuracy of the FEM and will be specified later. 

A multilevel MC (MLMC) estimator for the EIG, based on the error analysis by Giles \cite{Gil08}, was introduced in \cite{Bec20}. This estimator addressed the inner sampling and FEM discretization error. An MLMC estimator combined with antithetic sampling for EIG estimation was also proposed in~\cite{God20}, where the inexact evaluation of the inner integrand was not considered. Attempts to combine multilevel estimators with rQMC were made in \cite{Kuo15} where dimension truncation via the Karhunen--Lo\'{e}ve expansion was also considered; however, there was no nested integration structure. The works \cite{God18} and \cite{Xu20} introduced estimators for nested integrals based on the MLMC method, where the number of samples in the inner QMC estimator defines the levels. On the other hand, \cite{Fan22} proposed a multilevel rQMC estimator where the levels are defined by the number of samples in an inner MC estimator. Neither of these studies addressed FEM discretizations. Estimators for nested integrals with more than one nesting were analyzed in \cite{Rai18} and \cite{Sye23} for general applications. The latter of these studies developed a randomized (i.e., unbiased) MLMC method, achieving a complexity independent of the number of nestings. Randomized MLMC estimators for financial applications were proposed in \cite{Spe23, Gil24} and for Bayesian experimental design applications in \cite{God22}. The recent work \cite{Du24} proposed a multi-index MC (MIMC) estimator for EIG with discretized inner integrand that can achieve a computational cost of $\cl{O}(TOL^{-2})$. MIMC methods are an extension to MLMC methods, where an optimization among multi-indices of accuracies is performed \cite{haji2016multi}. This scheme allows for exploitation of anisotropies in the error convergence among different kinds of discretizations present in the integrand. The randomized double-loop QMC (rDLQMC) estimator was introduced in \cite{Bar23}, where rQMC methods were applied to estimate both integrals. This estimator achieves an optimal computational cost of almost $\cl{O}(TOL^{-10/9})$ for a high-dimensional example where both the inner and outer integrands are polynomials and a computational cost of almost $\cl{O}(TOL^{-4/3-\gamma/\eta_{\rm{w}}})$ for EIG applications. 

The two estimators based on rQMC methods with either a full tensor product structure or a sparse tensor product structure proposed in \cite{Kaa24} were geared toward inner integrands that could be evaluated exactly, and achieved computational work of almost $\cl{O}(TOL^{-2})$ and $\cl{O}(TOL^{-1})$, respectively. Moreover, the FEM-based approximation was treated as the exact integrand in \cite{Kaa24} whenever needed. These estimators use a truncation of the Gaussian observation noise to address singularities arising in the EIG application. As the error tolerance $TOL$ approaches zero, this truncation extends towards the boundaries of the integration domain, and as a result, a factor is introduced in their root-mean-square-error bound that grows exponentially with the noise dimension (see \cite[Theorem 6.2]{Kaa24}). Such behavior is not observed for the truncation used in the present work for error bounds based on the Hardy--Krause variation or Owen's boundary growth condition, see Corollary~\ref{cor:EIG.total.error} and the discussion in~\cite{Bar23}.

As our main contributions, we propose a novel estimator including a multilevel structure, combined with QMC and rQMC methods, analyze the resulting bias and statistical errors, and derive the optimal number of inner and outer samples and the discretization parameter. Moreover, we derive conditions for applying this estimator to the EIG of an experiment based on a truncation scheme introduced in a previous work \cite{Bar23}. We then demonstrate that our estimator achieves computational cost of nearly $\cl{O}(TOL^{-1-\gamma/\eta_{\rm{w}}})$ for EIG estimation in the presence of FEM discretization, and almost $\cl{O}(TOL^{-1})$ without such discretization. Moreover, a computational cost of nearly $\cl{O}(TOL^{-8/9})$ is achieved for an example where the inner and outer integrands are polynomials. Whenever the evaluation of the integrand that is approximated via MC (or rQMC) methods requires additional discretization methods, such as the inner rQMC method in this setting, multilevel techniques employ a telescoping structure of increasing accuracy of the integrand approximation to reduce the variance of the combined estimator. The integral is approximated for a very coarse inner integrand, and the result is then corrected via increasing levels of difference estimators.
The primary application of the estimator presented in this study is to approximate the EIG, where the nonlinear function separating the inner and outer integrals is the logarithm. The regularity assumptions stated in this work should be understood in this context; however, the analysis is conducted in general terms and can be adapted to more specific settings.

This work is structured as follows. In Section~\ref{sec:overview}, we provide a brief overview of numerical approximation methods of integrals and introduce useful notation. As the primary contribution, we present the multilevel double-loop QMC (MLDLQMC) estimator with inexact sampling in Section~\ref{sec:MLrDLQMC}, along with asymptotic error bounds of the bias and statistical error. This section concludes with a bound on the total work of this estimator for a given error tolerance based on a nearly optimal allocation of the inner and outer samples. In Section~\ref{sec:EIG.estimation}, we demonstrate how to apply the MLDLQMC estimator to estimate the EIG of an experiment. This analysis involves a truncation of the Gaussian observation noise present in the experiment model. Section~\ref{sec:numerical.results} provides numerical results to highlight the practical applicability and competitiveness of the proposed estimators for a Bayesian experimental design setting.

\section{Overview of numerical integration methods}\label{sec:overview}
For the reader's convenience, we provide an overview of existing numerical methods to estimate integrals and nested integrals that are relevant to this work. Moreover, we provide error bounds and analyze the computational cost associated with these methods.

\subsection{Monte Carlo method}\label{sec:MC}
Given an integral
\begin{equation}\label{eq:integral}
  I=\int_{[0,1]^d}\varphi(\bs{z})\di{}\bs{z},
\end{equation}
where $\varphi:[0,1]^d\to\bb{R}$ is square-integrable and $d$ is a positive integer, the MC estimator
\begin{equation}\label{eq:monte.carlo}
I_{\rm{MC}}\coloneqq\frac{1}{M}\sum_{m=1}^{M}\varphi(\bs{z}^{(m)})
\end{equation}
uses independent and identically distributed (iid) samples $\bs{z}^{(m)}$, where $1\leq m\leq M$, from the uniform distribution $\cl{U}\left([0,1]^d\right)$, to approximate $I$ in \eqref{eq:integral}. The central limit theorem (CLT) \cite{Dur19} reveals that the error of the MC estimator is proportional to the estimator variance $\bb{V}[I_{\rm{MC}}]$ in a probabilistic sense:
\begin{equation}\label{eq:CLT.prob}
    \bb{P}\left(|I_{\rm{MC}}-I| \le C_{\alpha} \sqrt{\bb{V}[I_{\rm{MC}}]}\right) \ge 1 - \alpha,
\end{equation}
as $M\to\infty$, where $\alpha\in(0,1)$ is a parameter, $C_{\alpha} = \Phi^{-1}(1-\alpha/2)$, $\Phi^{-1}$ denotes the inverse cumulative distribution function (CDF) of the standard normal distribution, and $C_{\alpha}\to\infty$ as $\alpha\to 0$. The bound~\eqref{eq:CLT.prob} motivates the use of variance reduction techniques to control the error. Independence of the samples in the MC estimator ensures that
\begin{align}\label{eq:MC.Var}
    \bb{V}[I_{\rm{MC}}]=\frac{\bb{V}[\varphi]}{M}.
\end{align}
Increasing the number of samples $M$ is a straightforward way to reduce the variance \eqref{eq:MC.Var}. The cost of the MC estimator is proportional to the number of times the function $\varphi$ is evaluated, i.e., $M$. Thus, for a linear increase in cost, we obtain a linear decrease in variance, and a sublinear decrease in the error \eqref{eq:CLT.prob}, namely, at rate $1/2$. A total number of samples proportional to $TOL^{-2}$ is required to achieve an error less than a prescribed tolerance $TOL>0$, i.e.,
\begin{equation}
    C_{\alpha} \sqrt{\frac{\bb{V}[\varphi]}{M}}<  TOL.
\end{equation}
Thus, the total computational work of the MC estimator is
\begin{equation}
    W_{MC}^{\ast}=\cl{O}\left(TOL^{-2}\right),
\end{equation}
as $TOL\to 0$.

\subsection{Quasi-Monte Carlo method}\label{sec:QMC}
The convergence rate of the MC error is independent of the dimension $d$. However, the conventional rate of $1/2$ is not large enough for many practical applications. One primary concern with iid samples is that they do not fill the integration domain well. The QMC method,
\begin{equation}\label{eq:quasi.monte.carlo}
 I_{\rm{Q}} \coloneqq \frac{1}{M}\sum_{m=1}^M\varphi(\bs{z}^{(m)}),
\end{equation}
where
\begin{equation}
    \bs{z}^{(m)}\coloneqq\bs{t}^{(m)}, \quad 1\leq m\leq M,
\end{equation}
instead uses deterministic points $\bs{t}^{(m)}$, where $1\leq m\leq M$. These are known as low-discrepancy points with the following property:
\begin{align}\label{eq:star.discrepancy}
D^{\ast}(\bs{t}^{(1)},\ldots,\bs{t}^{(M)}){}&\coloneqq \sup_{(s_1,\ldots,s_d)\in[0,1]^d}\left|\frac{1}{M}\sum_{m=1}^M\mathds{1}_{\{\bs{t}^{(m)}\in\prod_{i=1}^d [0,s_i)\}}-\prod_{i=1}^ds_i\right|\;,\nonumber\\
{}&\leq \frac{C_{\epsilon,d}}{M^{1-\epsilon}},
\end{align}
for any $\epsilon>0$, where $C_{\epsilon,d}>0$ is independent of $M$, and $C_{\epsilon,d}\to\infty$ as $\epsilon\to 0$, where $D^{\ast}$ is the star-discrepancy~\cite[Def.~(5.2)]{Caf98}. An important class of low-discrepancy sequences are digital sequences \cite{Nie92, Owe95, Owe03, Dic10}. We apply a digital sequence called the Sobol' sequence \cite{Sob67} throughout this work. However, other low-discrepancy sequences could also be used as long as the bound in~\eqref{eq:star.discrepancy} is satisfied. The approximation error of the QMC estimator is bounded by the Koksma--Hlawka inequality \cite{Hla61, Nie92, Caf98},
\begin{equation}\label{eq:QMC.bound}
    |I_{\rm{Q}}-I|\leq D^{\ast}(\bs{t}^{(1)},\ldots,\bs{t}^{(M)})V_{\rm{HK}}(\varphi),
\end{equation}
where $V_{\rm{HK}}(\varphi)$
is the recursively defined Hardy--Krause variation~\cite[Def.~(5.8)]{Caf98} for continuously differentiable $\varphi$:
\begin{equation}\label{eq:VHK}
    V_{\rm{HK}}(\varphi)\coloneqq \int_{[0,1]^d}\left|\left(\prod_{i=1}^d\frac{\partial}{\partial z_i}\right)\varphi(\bs{z})\right|\di{}\bs{z}+\sum_{i=1}^dV_{\rm{HK}}(\varphi(\bs{z})|_{z_i=1}),
\end{equation}
with $\varphi(\bs{z})|_{z_i=1}:[0,1]^{d-1}\to\bb{R}$ representing the restriction of $\varphi(\bs{z})$ to $\varphi(z_1,\ldots,z_i=1,\ldots,z_d)$. For $\varphi:[0,1]\to\bb{R}$, the Hardy--Krause variation is defined as follows:
\begin{equation}
    V_{\rm{HK}}(\varphi)\coloneqq \int_{[0,1]}\left|\frac{\di{}}{\di{}z}\varphi(z)\right|\di{}z.
\end{equation}
To achieve an error less than $TOL$, the QMC estimator consequently has computational work
\begin{equation}
    W_{Q}^{\ast}=\cl{O}\left(TOL^{-\frac{1}{1-\epsilon}}\right)
\end{equation}
as $TOL\to 0$ for integrands with bounded Hardy--Krause variation \eqref{eq:VHK}, based on the bounds in~\eqref{eq:star.discrepancy} and~\eqref{eq:QMC.bound}, where the multiplicative constant may depend on the integrand dimension $d$.
The main concerns with the QMC method are, in no particular order, that:
\begin{itemize}
    \item There may be a multiplicative constant in the error depending on the dimension $d$;
    \item The first point of many popular low-discrepancy sequences is the origin, i.e., $\bs{t}^{(1)}=\bs{0}$, including the Sobol' sequence;
    \item The Hardy--Krause variation in the error bound in~\eqref{eq:QMC.bound} is usually more restrictive and more challenging to compute or estimate than the estimator variance arising in \eqref{eq:CLT.prob}, which plays a crucial role in the error of the rQMC method;
    \item Many popular integrands exhibit infinite Hardy--Krause variation.
\end{itemize}
The first concern can be remedied for certain integrands and specialized QMC methods~\cite{Dic13}. The second concern can be elegantly addressed by introducing randomization to the deterministic low-discrepancy points, as will be discussed in the following section. The third concern can also be partially addressed by randomization. The sample variance taken over multiple independent randomizations, as will be elaborated in the next section, serves as a practical error estimate through Chebyshev's inequality. Regarding the fourth concern, there are certain types of popular integrands, e.g., integrands weighted by a normal density rather than the uniform one, which must be mapped back to the interval $[0,1]^d$. This mapping results in infinite Hardy--Krause variation of the resulting integrand, even though the corresponding QMC error is finite.

\subsection{Randomized Quasi-Monte Carlo method}\label{sec:rQMC}
The rQMC estimator for the integral \eqref{eq:integral} with $R$ randomizations is defined as follows:
\begin{equation}
 I_{\rm{rQ}}^{(R)} \coloneqq \frac{1}{R}\sum_{r=1}^R\frac{1}{M}\sum_{m=1}^M\varphi(\bs{z}^{(r,m)}),
\end{equation}
where
\begin{equation}
    \bs{z}^{(r,m)}\coloneqq\bs{\rho}^{(r)}(\bs{t}^{(m)}), \quad 1\leq r\leq R, \, 1\leq m\leq M.
\end{equation}
This estimator uses independent randomizations $\bs{\rho}^{(r)}\in\Pi$, where $1\leq r\leq R$, of the deterministic low-discrepancy sequence $\bs{t}^{(m)}$, where $1\leq m\leq M$, and $\Pi$ denotes the space of randomizations. Owen's scrambling for digital sequences \cite{Nie92, Owe95, Owe03} is employed throughout this work, where a random permutation is applied to the base-two digits of $\bs{t}^{(m)}$, where $1\leq m\leq M$. This randomization scheme preserves the low-discrepancy property with probability one (see \cite{Owe95, Dic10} for details on this method). It holds that
\begin{align}\label{eq:rQMC.Var}
    \bb{V}[I_{\rm{rQ}}^{(R)}]{}&=\bb{V}\left[\frac{1}{R}\sum_{r=1}^{R}\frac{1}{M}\sum_{m=1}^{M}\varphi(\bs{z}^{(r,m)})\right],\nonumber\\
    {}&=\frac{1}{R}\bb{V}\left[\frac{1}{M}\sum_{m=1}^{M}\varphi(\bs{\rho}(\bs{t}^{(m)}))\right],\nonumber\\
    {}&=\frac{1}{R}\bb{E}\left[\left|\frac{1}{M}\sum_{m=1}^{M}\varphi(\bs{\rho}(\bs{t}^{(m)}))-\bb{E}\left[\frac{1}{M}\sum_{m=1}^{M}\varphi(\bs{\rho}(\bs{t}^{(m)}))\right]\right|^2\right],\nonumber\\
    {}&=\frac{1}{R}\bb{E}\left[\left|\frac{1}{M}\sum_{m=1}^{M}\varphi(\bs{\rho}(\bs{t}^{(m)}))-I\right|^2\right],\nonumber\\
    {}&\leq \frac{1}{R}\bb{E}\left[\left|D^{\ast}\left(\bs{\rho}(\bs{t}^{(1)}),\ldots,\bs{\rho}(\bs{t}^{(M)})\right)V_{\rm{HK}}(\varphi)\right|^2\right],
\end{align}
where the last line follows from the Koksma--Hlawka inequality \eqref{eq:QMC.bound} and the variance and expectation are taken over the randomization $\bs{\rho}$. The variance of the rQMC estimator can thus be bounded as
\begin{equation}\label{eq:rQMC.Var.HK}
    \bb{V}[I_{\rm{rQ}}^{(R)}]\leq \frac{C_{\epsilon,d}^2(V_{\rm{HK}}(\varphi))^2}{RM^{2-2\epsilon}},
\end{equation}
for any $\epsilon>0$ and $C_{\epsilon,d}\to\infty$ as $\epsilon\to 0$ via~\eqref{eq:star.discrepancy} for all sequences that maintain the low-discrepancy property in $L^2(\Pi,\bb{R})$ under randomization, such as the Sobol' sequence with Owen's scrambling. It follows that increasing the number of QMC points is often more efficient than increasing the number of randomizations. There are two important concerns with the Hardy--Krause variation. First, several important classes of functions have infinite Hardy--Krause variation and finite variance. Thus, for such cases, the Koksma--Hlawka inequality does not provide a useful error bound. Second, it is challenging to obtain practical estimates of the Hardy--Krause variation, whereas the variance can be estimated using the sample variance. Regarding the latter point, the rQMC variance $\bb{V}[I_{\rm{rQ}}^{(R)}]$ can be estimated using $R$ independent randomizations whereas the bound in~\eqref{eq:rQMC.Var} is only for theoretical purposes. This situation presents a trade-off: increasing the number of randomizations improves the reliability of the variance estimate but only improves the estimate of the integral, i.e., expectation, at the slower MC rate. A naive approximation based on the CLT for finite $R$ can be problematic for certain integrands \cite{Tuf04, Lec10}. Chebyshev's inequality serves as an alternative, at the expense of a larger confidence constant $C_{\alpha}=1/\sqrt{\alpha}$. The number of randomizations $R$ plays a lesser role in theoretical error bounds, and we thus introduce the rQMC estimator with one randomization as
\begin{equation}\label{eq:randomized.quasi.monte.carlo}
    I_{\rm{rQ}}\coloneqq I_{\rm{rQ}}^{(1)}.
\end{equation}
As for the first point mentioned above, the following result in \eqref{eq:RMSE} provides an error bound for a specific, yet common, type of integrands that have infinite Hardy--Krause variation. We begin by introducing Assumption~\ref{eq:boundary.growth} (from \cite{Owe06}):
\begin{asu}[Boundary growth condition]\label{eq:boundary.growth}
Let $\varphi:[0,1]^d\to\bb{R}$ and assume that there exists $0<b<\infty$ and $A_i>0$ for $1\leq i\leq d$, where
\begin{equation}\label{eq:A.max}
A_{\rm{max}}\coloneqq\max_{1\leq i\leq d}A_i<\frac{1}{2},
\end{equation}
such that
\begin{equation}
    \left|\left(\prod_{j\in u}\frac{\partial}{\partial z_j}\right)\varphi(\bs{z})\right|\leq b\prod_{i=1}^d\min (z_i,1-z_i)^{-A_i-\mathds{1}_{\{i\in u\}}}
\end{equation}
for all $\bs{z}=(z_1,\ldots,z_d)\in(0,1)^d$ and all $u\subseteq\{1,\ldots,d\}$ with the convention that $(\prod_{j\in u}\partial/\partial z_j)\varphi\equiv \varphi$ for $u=\emptyset$.
\end{asu}
Assumption~\ref{eq:boundary.growth} implies, among other things, that $\varphi\in L^2([0,1]^{d})$, i.e., that $\varphi(\bs{Z})$, where $\bs{Z}$ is a uniform random variable, has finite variance. This allows for the application of standard error bounds for the MC method. Next, the following proposition is stated (from \cite{He23}):
\begin{prop}[Thm.~3.1 in \cite{He23}]\label{cor:he}
Given Assumption~\ref{eq:boundary.growth}, the rQMC estimator \eqref{eq:randomized.quasi.monte.carlo} satisfies the following $L^2([0,1]^{d})$ error bound:
\begin{equation}\label{eq:RMSE}
        \bb{E}[|I_{\rm{rQ}}-I|^2]^{\frac{1}{2}}\leq b B_A C_{\epsilon,d} M^{-1+\epsilon+A_{\rm{max}}},
    \end{equation}
for all $\epsilon>0$, where $B_A\to\infty$ as $\min_{1\leq i\leq d}A_i\to 0$ and $C_{\epsilon,d}\to\infty$ as $\epsilon\to 0$.
    \end{prop}
This proposition is an extension of an earlier result from \cite{Owe06}, where a similar bound was provided for the $L^1([0,1]^{d})$ error. 
The rQMC method offers another important benefit over the deterministic QMC method. If the integrand satisfies an even stricter regularity requirement than the bounded Hardy--Krause variation stated above, convergence rates of the estimator variance of almost three are possible. In particular, for continuously differentiable integrands $\varphi$, let the generalized Vitali variation of order one in dimension $d$ (see \cite[Chapter 13]{Dic10} and \cite[Appendix F]{Bar23}) be defined as follows:
    \begin{equation}\label{def:GVV}
    V^{(d)}(\varphi)\coloneqq \left(\int_{[0,1]^d}\left|\frac{\partial^d\varphi}{\partial z_1\ldots \partial z_d}(\bs{z})\right|^2\di{}\bs{z}\right)^{\frac{1}{2}}.
    \end{equation}
    Next, let $\cl{I}_{d}\coloneqq\{1,\ldots,d\}$, and $\mathfrak{u}$ be the set of all subsets of $\cl{I}_d$ and define
    \begin{equation}
        \varphi_{\mathfrak{u}}(\bs{z}_{\mathfrak{u}})\coloneqq \int_{[0,1]^{d-|\mathfrak{u}|}}\varphi(\bs{z})\di{}\bs{z}_{\cl{I}_d\setminus\mathfrak{u}},
    \end{equation}
    where $|\mathfrak{u}|$ denotes the cardinality of $\mathfrak{u}$ and $\bs{z}_{\mathfrak{u}}$ contains those entries of $\bs{z}$ that are in $\mathfrak{u}$ whereas $\bs{z}_{\cl{I}_d\setminus\mathfrak{u}}$ contains those entries of $\bs{z}$ that are not in $\mathfrak{u}$. The generalized Hardy--Krause variation of order one of $\varphi$ on $[0,1]^d$ is subsequently defined as follows:
    \begin{equation}
        V_{\rm{GHK}}(\varphi)\coloneqq \left(\left|\int_{[0,1]^d}\varphi(\bs{z})\di{}\bs{z}\right|^2+\sum_{\emptyset\neq \mathfrak{u}\subseteq\cl{I}_{d}}\left(V^{(|\mathfrak{u}|)}(\varphi_{\mathfrak{u}})\right)^2\right)^{\frac{1}{2}}.
    \end{equation}
The above definition essentially differs from the one in~\ref{eq:VHK} in that the mixed derivatives must be square-integrable. 
\begin{prop}[Theorem 6.25 in~\cite{Dic13}]\label{prop:GHK}
For integrands $\varphi$ with bounded generalized Hardy--Krause variation of order one, the variance of an rQMC estimator with $R=1$ randomizations and $M$ Sobol' points randomized via Owen's scrambling is bounded as follows:
\begin{equation}\label{prop.BGHK.var}
    \bb{V}[I_{\rm{rQ}}]\leq \left(V_{\rm{GHK}}(\varphi)\right)^2C_{\epsilon, d}^2M^{-3+2\epsilon},
\end{equation}
where $C_{\epsilon, d}\to\infty$ as $\epsilon\to0$.
\end{prop}
From the above results, it follows that the computational work of the rQMC estimator to achieve a root mean square error less than $TOL$ is
\begin{equation}
    W_{rQ}^{\ast}=\cl{O}\left(TOL^{-\frac{1}{1-\epsilon}}\right)
\end{equation}
as $TOL\to 0$ for integrands with bounded Hardy--Krause variation, and
\begin{equation}
    W_{rQ}^{\ast}=\cl{O}\left(TOL^{-\frac{1}{1-\epsilon-A_{\rm{max}}}}\right)
\end{equation}
as $TOL\to 0$ for integrands with singularities at the boundary satisfying Assumption~\ref{eq:boundary.growth}, where $A_{\rm{max}}>0$ as in~\eqref{eq:A.max} indicates the severity of the singularity. The condition that $A_{\rm{max}}<1/2$ also implies that the rQMC error converges faster than the MC error. Finally, the computational work of the rQMC estimator to achieve a root mean square error less than $TOL$ is
\begin{equation}
    W_{rQ}^{\ast}=\cl{O}\left(TOL^{-\frac{2}{3-2\epsilon}}\right)
\end{equation}
as $TOL\to 0$ for integrands with bounded generalized Hardy--Krause variation of order one.

\subsection{Double-loop Monte Carlo method for nested integrals}\label{sec:DLMC}
Following \cite{Bar23}, a nested integral is defined as follows:
\begin{set}[Nested integral]
Let $f:\bb{R}\to\bb{R}$ be a nonlinear function and let $g:[0,1]^{d_1}\times[0,1]^{d_2}\to\bb{R}$ define a nonlinear relation between $\bs{x}\in[0,1]^{d_2}$ and $\bs{y}\in[0,1]^{d_1}$, where $d_1, d_2$ are positive integers. Then, a nested integral is defined as follows:
\begin{align}\label{eq:double.loop}
  I {}&= \int_{[0,1]^{d_1}}f\left(\int_{[0,1]^{d_2}}g(\bs{y},\bs{x})\di{}\bs{x}\right)\di{}\bs{y}.
\end{align}
\end{set}
Depending on the context, it may be beneficial to write the nested integral as a nested expectation instead, that is, 
\begin{equation}
I=\bb{E}[f(\bb{E}[g(\bs{y},\bs{x})|\bs{y}])],
\end{equation}
where $\bs{y}\sim\cl{U}\left([0,1]^{d_1}\right)$ and $\bs{x}\sim\cl{U}\left([0,1]^{d_2}\right)$. Moreover, the notation $\bar{g}:[0,1]^{d_1}\mapsto \bb{R}$ is defined as follows:
\begin{align}\label{eq:bar.g}
\bar{g}(\bs{y}){}&\coloneqq\mathbb{E}[g(\bs{y},\bs{x})|\bs{y}],\nonumber\\
{}&=\int_{[0,1]^{d_2}} g(\bs{y},\bs{x})\di{}\bs{x}.
\end{align}
Ryan~\cite{Rya03} introduced the DLMC estimator for nested integrals:
\begin{equation}\label{eq:dlmc}
I_{\rm{DLMC}} \coloneqq \frac{1}{N}\sum_{n=1}^Nf\left(\frac{1}{M}\sum_{m=1}^Mg(\bs{y}^{(n)},\bs{x}^{(n,m)})\right),
\end{equation}
where $\bs{y}^{(n)}$, for $1\leq n\leq N$, are iid samples from $\cl{U}\left([0,1]^{d_1}\right)$. Moreover, $\bs{x}^{(n,m)}$, for $1\leq n\leq N$ and $1\leq m\leq M$, are iid samples from $\cl{U}\left([0,1]^{d_2}\right)$ for the specific application of EIG estimation where $f$ is the logarithm and $g$ is the likelihood of the data observation. For large values of $N$ and $M$, this estimator has a total error based on the CLT \cite{Rya03,Bec18} of
\begin{equation}
    |I_{\rm{DLMC}}-I|\leq \frac{C_{\rm{DL,B}}}{M}+
    C_{\alpha}
    \sqrt{\frac{C_{\rm{DL,V}}^{(1)}}{N}+\frac{C_{\rm{DL,V}}^{(2)}}{NM}}, \quad \text{with probability at least $1-\alpha$},
\end{equation}
for some $\alpha\in(0,1)$, where the first term is a bias stemming from the variance of the inner MC method, and the second term is derived from the total variance of the DLMC estimator, and $C_{\alpha}\to\infty$ as $\alpha\to 0$. Moreover, $C_{\rm{DL,B}}>0$ and $C_{\rm{DL,V}}^{(1)}, C_{\rm{DL,V}}^{(2)}>0$ are constants independent of $N$ and $M$. A total of $N\times M$ samples is required for this estimator, and thus the total work to achieve an error less than $TOL$ is:
\begin{equation}
    W_{DL}^{\ast}=\cl{O}\left(TOL^{-3}\right)
\end{equation}
as $TOL\to 0$.

\subsection{Double-loop randomized quasi-Monte Carlo method for nested integrals}\label{sec:rDLQMC}
The rDLQMC estimator as introduced in \cite{Bar23} is:
\begin{equation}\label{eq:dlqmc.S.R}
I_{\rm{rDLQ}}^{(S,R)} \coloneqq \frac{1}{S}\sum_{s=1}^S\frac{1}{N}\sum_{n=1}^Nf\left(\frac{1}{R}\sum_{r=1}^R\frac{1}{M}\sum_{m=1}^Mg(\bs{y}^{(s,n)},\bs{x}^{(s,n,r,m)})\right),
\end{equation}
where the samples take the following form
\begin{equation}
    \bs{y}^{(s,n)}\coloneqq\bs{\tau}^{(s)}(\bs{u}^{(n)}), \quad 1\leq s\leq S, \, 1\leq n\leq N.
\end{equation}
Here, $\bs{\tau}^{(s)}$, with $1\leq s\leq S$, are independent randomizations of the Sobol' sequence $\bs{u}^{(n)}$ in $d_1$, where $1\leq n\leq N$, using Owen's scrambling, and
\begin{equation}\label{eq:x.nrm}
    \bs{x}^{(s,n,r,m)}\coloneqq\bs{\rho}^{(s,n,r)}(\bs{t}^{(m)}), \quad 1\leq s\leq S, \, 1\leq n\leq N, \, 1\leq r\leq R, \, 1\leq m\leq M,
\end{equation}
where $\bs{\rho}^{(s,n,r)}$, with $1\leq s\leq S$, $1\leq n\leq N$, and $1\leq r\leq R$ are independent randomizations of the Sobol' sequence $\bs{t}^{(m)}$ in $d_2$, where $1\leq m\leq M$, using Owen's scrambling. Assuming that the nested integrand $f(\bar{g})$ satisfies Assumption~\ref{eq:boundary.growth} with $\bar{g}$ as in~\eqref{eq:bar.g}, for $0<A_{\rm{max}}<1/2$, whereas the inner integrand $g$ has bounded generalized Hardy--Krause variation of order one, this estimator has a total error based on Chebyshev's inequality \cite{Bar23} of
\begin{equation}
    |I_{\rm{rDLQ}}^{(1,1)}-I|\leq \frac{C_{\rm{DLQ,B}}}{M^{3-2\epsilon}}+
    C_{\alpha}
    \sqrt{\frac{C_{\rm{DLQ,V}}^{(1)}}{N^{2-2\epsilon-2A_{\rm{max}}}}+\frac{C_{\rm{DLQ,V}}^{(2)}}{NM^{3-2\epsilon}}}, \quad \text{with probability at least $1-\alpha$},
\end{equation}
for some $\alpha\in(0,1)$, $C_{\alpha}\to\infty$ as $\alpha\to 0$, for  $S=R=1$ randomizations, where $C_{\rm{DLQ,B}}>0$ and $C_{\rm{DLQ,V}}^{(1)}, C_{\rm{DLQ,V}}^{(2)}>0$ are constants independent of $N$ and $M$. The total work to achieve an error less than $TOL$ based on this error bound is:
\begin{equation}
    W_{rDLQ}^{\ast}=\cl{O}\left(TOL^{-\frac{1}{1-\epsilon-A_{\rm{max}}}-\frac{1}{3-2\epsilon}}\right)
\end{equation}
as $TOL\to 0$. This setting is also analyzed in particular for the EIG estimation in \cite{Bar23}.

\subsection{Multilevel double-loop Monte Carlo method for nested integrals}\label{sec:MLDL}
The multilevel method for nested integrals \cite{God20} constructs a sequence of increasingly accurate MC estimators for the inner integral to harness variance reduction achieved by taking differences of adjacent levels evaluated for the same random variables. The general form of the MLMC estimator for a nonnested integral is given as follows:
\begin{equation}\label{eq:ML}
    I_{\rm{ML}}\coloneqq \frac{1}{N_0}\sum_{n=1}^{N_{0}}\varphi_{0}\left(\bs{z}^{(0,n)}\right)+\sum_{\ell=1}^{L}\frac{1}{N_{\ell}}\sum_{n=1}^{N_{\ell}}\left[\varphi_{\ell}\left(\bs{z}^{(\ell,n)}\right)-\varphi_{\ell-1}\left(\bs{z}^{(\ell,n)}\right)\right],
\end{equation}
where $\bs{z}^{\ell,n}$ are iid samples from $\cl{U}([0,1]^d)$ for $1\leq n\leq N_{\ell}$, and $1\leq\ell\leq L$, $L\in\bb{N}$. The multilevel double-loop MC (MLDLMC) estimator introduced in \cite{God20} for nested integrals defines the level $\varphi_{\ell}$ as follows:
\begin{equation}
    \varphi_{\ell}(\bs{z}^{(\ell,n)})\coloneqq f\left(\frac{1}{M_{\ell}}\sum_{m=1}^{M_{\ell}}g\left(\bs{y}^{(\ell,n)},\bs{x}^{(\ell,n,m)}\right)\right),
\end{equation}
where $\bs{z}^{(\ell,n)}\coloneqq(\bs{y}^{(\ell,n)},\bs{x}^{(\ell,n,m)})$ and $M_{\ell}=M_02^{\ell}$, for $M_0>0$. An antithetic construction first introduced in \cite{Gil14} to increase the coupling, and thus reduce the variance further, was used to obtain computational work:
\begin{equation}
    W_{\rm{ML}}^{\ast}=\cl{O}\left(TOL^{-2}\right)
\end{equation}
to achieve an error tolerance $TOL>0$ as $TOL\to 0$ for estimating the EIG. Note that the asymptotic work presented in \cite{God20} was derived under the assumption that exact evaluations of $g$ were available. Depending on the considered application, this assumption may be inadequate, because the considered experiment generating the data may involve numerical approximations.

\subsection{Inexact sampling}\label{sec:Inexact.sampling}
Inner integrands $g$ encountered in many practical applications cannot be evaluated exactly at a given point $\bs{z}$. Instead, an approximation $g_{h}$ is introduced, where $g_{h}\to g$ as $h\to 0$ in a suitable sense, for a discretization parameter $h>0$. This setting is referred to as inexact sampling and examples include finite difference methods, FEM, and time stepping methods. This discretization introduces bias to the estimators discussed in the previous sections, and thus affects the computational work required to achieve an error tolerance $TOL$. To account for this additional bias, the following assumptions are considered in this work:
\begin{asu}[Weak convergence rate]\label{asu:weak.rate}
    Let $g_{h}$ be an approximation of $g:[0,1]^{d}\mapsto \bb{R}$, and $h>0$ such that the weak error of $g_{h}$ is given by
\begin{equation}
\left|\int_{[0,1]^{d_1}}f\left(\int_{[0,1]^{d_2}}g_h(\bs{y},\bs{x})\di{}\bs{x}\right)\di{}\bs{y}-\int_{[0,1]^{d_1}}f\left(\int_{[0,1]^{d_2}}g(\bs{y},\bs{x})\di{}\bs{x}\right)\di{}\bs{y}\right| \leq C_{\rm{w}}h^{\eta_{\rm{w}}},
\end{equation}
 where $\eta_{\rm{w}}>0$ is the weak $h$-convergence rate, and $C_{\rm{w}}>0$ is a constant independent of $h$. The work of evaluating $g_{h}$ is assumed to be of order $\cl{O}(h^{-\gamma})$, for $\gamma > 0$ uniformly in $[0,1]^d$.
\end{asu}
In the inexact sampling case, an additional term ${-\gamma/\eta_{\rm{w}}}$ is added to the exponent of $TOL$ in the total computational work of the MC, QMC, rQMC, DLMC, and rDLQMC estimators \cite{Bec18,Bar23}, which is motivated by the EIG application. Note that only the inner integrand in the nested setting is assumed to require discretization throughout this work. The work \cite{Bec20} proposed an MLMC estimator for nested integration applied to the EIG, where the data generation involves approximating a PDE via the FEM. The accuracy level in the error bound proposed in this work was determined by the inner number of samples and the discretization parameter of the FEM. However, the total computational work is then bounded for a fixed inner number of samples as $\cl{O}(TOL^{-2})$ for large and intermediate error tolerances $TOL>0$, and $\cl{O}(TOL^{-3})$ asymptotically as $TOL\to 0$. This change in complexity regime is due to the effects of Laplace-based importance sampling used in this estimator, which renders just one inner sample sufficient for $TOL$ large enough.
\begin{rmk}[Laplace-based importance sampling]
   Laplace-based importance sampling is an efficient variance reduction technique for EIG estimation \cite{Lon13, Bec18, Bec20}. However, combining this method with multilevel strategies and especially with rQMC methods is nontrivial \cite{Ouy23}. Thus, we leave the study of this method for future work.
\end{rmk}

\section{Multilevel double-loop randomized quasi-Monte Carlo method for nested integrals with inexact sampling}\label{sec:MLrDLQMC}
As the primary focus of this work we analyzed the multilevel double-loop randomized quasi-Monte Carlo (MLDLQMC) estimator with inexact sampling. 
This estimator combines the multilevel technique with both deterministic QMC and rQMC methods. 
At all levels $\ell$, where $0\leq\ell\leq L$, the outer integral is estimated using the rQMC method, allowing for a practical error estimate via the sample variance. This estimate is theoretically justified by Proposition~\ref{cor:he}. The inner estimation at levels $\ell$, where $0\leq\ell\leq L-1$, is based on QMC methods. Due to the deterministic nature of these methods, the variances of the estimators for those levels depend only on the outer samples. The multilevel structure ensures that only the final level $L$ contributes to the bias of the overall estimator.
The rQMC method is applied for the inner integral approximation only at level $L$, allowing for practical error estimates of the bias and variance of the estimator.
\subsection{Estimator introduction}
\begin{set}[Multilevel double-loop randomized quasi-Monte Carlo estimator with inexact sampling]
The MLDLQMC estimator with inexact sampling is defined as follows:
\begin{align}
    I_{\rm{MLDLQ}}^{(S,R)}\coloneqq{}& \frac{1}{S}\sum_{s=1}^{S}\frac{1}{N_0}\sum_{n=1}^{N_0}f\left(\frac{1}{M_0}\sum_{m=1}^{M_0}g_{h_0}\left(\bs{y}^{(0,s,n)},\bs{x}^{(m)}\right)\right)\nonumber\\
    {}&+
    \sum_{\ell=1}^{L-1} \frac{1}{S}\sum_{s=1}^{S}
    \frac{1}{N_{\ell}}\sum_{n=1}^{N_{\ell}}\left[f\left(\frac{1}{M_{\ell}}\sum_{m=1}^{M_{\ell}}g_{h_{\ell}}\left(\bs{y}^{(\ell,s,n)},\bs{x}^{(m)}\right)\right) -f\left(\frac{1}{M_{\ell-1}}\sum_{m=1}^{M_{\ell-1}}g_{h_{\ell-1}}\left(\bs{y}^{(\ell,s,n)},\bs{x}^{(m)}\right)\right)\right]\nonumber\\
    {}&+\frac{1}{S}\sum_{s=1}^{S}\frac{1}{N_{L}}\sum_{n=1}^{N_{L}}\left[f\left(\frac{1}{R}\sum_{r=1}^{R}\frac{1}{M_{L}}\sum_{m=1}^{M_{L}}g_{h_{L}}\left(\bs{y}^{(L,s,n)},\bs{x}^{(s,r,m)}\right)\right)\right. \nonumber\\
    {}&\left.\quad\quad\quad\quad\quad\quad\quad\quad-f\left(\frac{1}{M_{L-1}}\sum_{m=1}^{M_{L-1}}g_{h_{L-1}}\left(\bs{y}^{(L,s,n)},\bs{x}^{(m)}\right)\right)\right].
\end{align}
Here, $g_{h_{\ell}}$ is an approximation of $g$ such that Assumption~\ref{asu:weak.rate}
holds for $0\leq \ell\leq L$, and
\begin{equation}\label{eq:y.ell.h}
    \bs{y}^{(\ell,s,n)}\coloneqq \bs{\tau}^{(\ell,s)}(\bs{u}^{(n)}), \quad 0\leq \ell \leq L, \,1\leq s\leq S, \, 1\leq n\leq N_{\ell},
\end{equation}
where $\bs{\tau}^{(\ell,s)}$ are iid Owen's scramblings of points in the Sobol' sequence $\bs{u}^{(n)}$ in dimension $d_1$, for $1\leq s\leq S$, $1\leq n\leq N_{\ell}$, and $0\leq \ell\leq L$, such that
\begin{equation}
\bb{E}[D^{\ast}(\bs{y}^{(\ell,s,1)},\ldots,\bs{y}^{(\ell,s,N_{\ell})})]\leq \frac{C_{\epsilon,d_1}}{N_{\ell}^{1-\epsilon}}, \quad 0\leq \ell \leq L, \,1\leq s\leq S,
\end{equation} 
for all $\epsilon>0$, where $C_{\epsilon,d_1}>0$ is independent of $N_{\ell}$, where $0\leq \ell \leq L$, and $C_{\epsilon,d_1}\to\infty$ as $\epsilon\to 0$. Moreover,
\begin{equation}\label{eq:x.ell.h.det}
    \bs{x}^{(m)}\coloneqq \bs{t}^{(m)}, \quad 1\leq m\leq M_{\ell},
\end{equation}
are points in the Sobol' sequence $\bs{t}^{(m)}$ in dimension $d_2$, for $1\leq m\leq M_{\ell}$, $0\leq \ell\leq L-1$, such that
\begin{equation}
D^{\ast}(\bs{x}^{(1)},\ldots,\bs{x}^{(M_{\ell})})\leq \frac{C_{\epsilon,d_2}}{M_{\ell}^{1-\epsilon}}, \quad 0\leq \ell \leq L,
\end{equation} 
for all $\epsilon>0$, where $C_{\epsilon,d_2}>0$ is independent of $M_{\ell}$, where $0\leq \ell \leq L$, and $C_{\epsilon,d_2}\to\infty$ as $\epsilon\to 0$. Finally, the following construction is applied for the randomized inner points at level $L$:
\begin{equation}\label{eq:x.ell.h.det.L}
    \bs{x}^{(s,r,m)}\coloneqq \bs{\rho}^{(s,r)}(\bs{t}^{(m)}), \quad 1\leq s\leq S, \, 1\leq r\leq R, \, 1\leq m\leq M_{L},
\end{equation}
where $\bs{\rho}^{(s,r)}$ are iid Owen's scramblings of points in the Sobol' sequence $\bs{t}^{(m)}$ in dimension $d_2$, for $1\leq s\leq S$, $1\leq r\leq R$, and $1\leq m\leq M_{L}$, such that
\begin{equation}
\bb{E}[D^{\ast}(\bs{x}^{(s,r,1)},\ldots,\bs{x}^{(s,r,M_L)})]\leq \frac{C_{\epsilon,d_2}}{M_{L}^{1-\epsilon}}, \quad 1\leq s\leq S, \, 1\leq r\leq R,
\end{equation} 
for all $\epsilon>0$, where $C_{\epsilon,d_2}>0$ is independent of $M_{L}$ and $C_{\epsilon,d_2}\to\infty$ as $\epsilon\to 0$.
\end{set}
For now, we set:
\begin{equation}\label{eq:M.ell.shape}
    M_{\ell}\coloneqq M_0 2^{\ell}, \quad 1\leq \ell\leq L,
\end{equation}
where $M_0\geq 1$, and
\begin{equation}\label{eq:h.ell.shape}
    h_{\ell}\coloneqq h_02^{-\ell}, \quad 1\leq \ell\leq L,
\end{equation}
where $h_0>0$.
The choices~\eqref{eq:M.ell.shape} and~\eqref{eq:h.ell.shape} are not specific to the MLDLQMC estimator and are revisited later (see Remark~\ref{rmk:optimal.M.h}). It is crucial for the quality of the estimator that all point sets used have the low-discrepancy property~\eqref{eq:star.discrepancy}; however, it would be possible to use different points $\bs{x}^{(m)}$, where $1\leq m\leq M_{\ell}$, at different levels $1\leq \ell\leq L-1$, as long as the telescoping property is preserved. The following notations are defined for convenience:
\begin{equation}
I_{\rm{MLDLQ}}^{(S,R)}= I_{\rm{MLDLQ}_0^0}^{(S,R)}+\sum_{\ell=1}^{L}\left[I_{\rm{MLDLQ}_{\ell}^{\ell}}^{(S,R)}-I_{\rm{MLDLQ}_{\ell}^{\ell-1}}^{(S,R)}\right],
\end{equation}
where
\begin{equation}
    I_{{\rm{MLDLQ}}_{\ell}^{k}}^{(S,R)}\coloneqq \frac{1}{S}\sum_{s=1}^{S}
    \frac{1}{N_{\ell}}\sum_{n=1}^{N_{\ell}}f\left(\frac{1}{M_{k}}\sum_{m=1}^{M_{k}}g_{h_{k}}\left(\bs{y}^{(\ell,s,n)},\bs{x}^{(m)}\right)\right), \quad 1\leq\ell\leq L,\, 1\leq k\leq L-1,
\end{equation}
and
\begin{equation}
    I_{{\rm{MLDLQ}}_{L}^{L}}^{(S,R)}\coloneqq \frac{1}{S}\sum_{s=1}^{S}
    \frac{1}{N_{L}}\sum_{n=1}^{N_{L}}f\left(\frac{1}{R}\sum_{r=1}^{R}\frac{1}{M_{L}}\sum_{m=1}^{M_{L}}g_{h_{L}}\left(\bs{y}^{(L,s,n)},\bs{x}^{(s,r,m)}\right)\right).
\end{equation}
Here, two additional indices were introduced to distinguish between levels of the estimator and the two terms within each level difference. The MLDLQMC estimator with $S=R=1$ randomizations is defined as follows:
\begin{equation}\label{eq:MLrQ.h.det}
    I_{\rm{MLDLQ}}\coloneqq I_{\rm{MLDLQ}}^{(1,1)},
\end{equation}
with $I_{{\rm{MLDLQ}}_{\ell}^{k}}$ defined analogously for $1\leq\ell,k\leq L$,
and $\bar{g}_h$ is defined analogously to~\eqref{eq:bar.g} for $\bs{y}\in[0,1]^{d_1}$ and $\bs{x}\in[0,1]^{d_2}$. We proceed by bounding the total error of the MLDLQMC estimator based on certain assumptions that are verified for the EIG application, where $f$ is the logarithm and $g$ is the likelihood combined with appropriate inverse CDFs. The total error is split into bias and statistical errors as follows:
\begin{equation}\label{eq:MLQMC.error}
|I_{\rm{MLDLQ}}-I|\leq\underbrace{|\mathbb{E}[I_{\rm{MLDLQ}}]-I|}_{\text{bias error}}+\underbrace{|I_{\rm{MLDLQ}}-\mathbb{E}[I_{\rm{MLDLQ}}]|}_{\text{statistical error}}.
\end{equation}
\subsection{Bias error bounds}
This study follows the work by \cite{Bar23} and restates the necessary assumptions for completeness to bound the bias of the MLDLQMC estimator.
\begin{asu}[Integrability assumption for the composition of the outer and inner integrand]\label{asu:f.3.diff}
There exists $0<\delta_{K} <1$ and $1\leq p\leq\infty$ such that $\tilde{f}(\cdot)\equiv f(K\bar{g}_{h}(\cdot)):[0,1]^{d_1}\to\bb{R}$ is an element of the Sobolev space $W^{3,p}([0,1]^{d_1})$ for all $K\in(1-\delta_{K},1+\delta_{K})$, where $\bar{g}_{h}$ is as in \eqref{eq:bar.g}.
\end{asu}
\begin{asu}[Inverse inequality for the inner integrand]\label{asu:g.L4}
    Let $g_{h}(\cdot,\cdot):[0,1]^{d_1}\times[0,1]^{d_2}\to\bb{R}$ and assume that there exists $0\leq k<\infty$ such that
    \begin{equation}\label{eq:reg:g}
    \sup_{\bs{y}\in[0,1]^{d_1}} \left\lvert \frac{V_{\rm{HK}}(g_{h}(\bs{y},\cdot))}{ \bar{g}_{h}(\bs{y})} \right\rvert \leq k\;,
\end{equation}
with $V_{\rm{HK}}$ as in \eqref{eq:VHK}, $g_h$ and $\bar{g}_h$ as in Assumption~\ref{asu:weak.rate} and \eqref{eq:bar.g} as $h\to 0$. Furthermore, if Assumption~\ref{asu:f.3.diff} holds, assume that there exists $1\leq q\leq\infty$ with $1/p+1/q=1$ for $p$ as in Assumption~\ref{asu:f.3.diff} such that $\bar{g}_{h}\in L^{3q}([0,1]^{d_1})$.
\end{asu}
Together, these assumptions ensure that the bias from the inner sampling in the MLDLQMC estimator remains bounded. The assumptions depend on the particular integrands $f$ and $g$ and must be verified for each application. For the EIG, it is demonstrated that these assumptions hold in \cite{Bar23} and also in Section~\ref{sec:EIG.estimation} and Appendix~\ref{app:vcond} of this work, contingent on certain auxiliary assumptions on the experiment under consideration.
The expectation of the MLDLQMC estimator is equal to the expectation of the corresponding single-level estimator at the finest level $L$, that is,
\begin{equation}
    \mathbb{E}[I_{\rm{MLDLQ}}]=\mathbb{E}[I_{{\rm{MLDLQ}}_L^L}].
\end{equation}
The bound on the bias of the MLDLQMC estimator follows immediately from \cite[Proposition~3]{Bar23} with $M= M_L$ and $h= h_L$ as stated below.
\begin{prop}[Bias of the MLDLQMC estimator with inexact sampling]\label{prop:B.MLrQ-Q.h}
Given Assumption~\ref{asu:weak.rate} for $g_{h_L}(\cdot,\cdot):[0,1]^{d_1}\times[0,1]^{d_2}\to\bb{R}$ and Assumptions~\ref{asu:f.3.diff} and~\ref{asu:g.L4}, if $|f'''|$ is monotonic, the bias of the MLDLQMC estimator \eqref{eq:MLrQ.h.det} with inexact sampling has the following upper bound:
\begin{equation}\label{eq:Bias.constraint.h.det}
  |\mathbb{E}[I_{\rm{MLDLQ}}]-I|\leq C_{\mathrm{w}}h_L^{\eta_{\rm{w}}}+ \frac{\bb{E}[|\bar{g}_{h}|^2\left|f''(\bar{g}_{h})\right|]k^2C_{\epsilon,d_2}^2}{2M_L^{2-2\epsilon}}+\frac{\bb{E}[\left|\bar{g}_{h}\right|^3\left|f'''(K\bar{g}_{h})\right|]k^3C_{\epsilon,d_2}^3}{6M_L^{3-3\epsilon}},
\end{equation}
for any $\epsilon>0$, where $C_{\epsilon,d_2}\to\infty$ as $\epsilon\to 0$.
\end{prop}
\subsection{Statistical error bounds}
The statistical error of the MLDLQMC estimator is bounded in probability using Chebyshev's inequality as follows:
\begin{equation}\label{eq:statistical.error}
    \bb{P}\left(|I_{\rm{MLDLQ}}-\bb{E}[I_{\rm{MLDLQ}}]|\leq C_{\alpha}\sqrt{\bb{V}[I_{\rm{MLDLQ}}]}\right)\geq 1-\alpha
\end{equation}
for $\alpha\in(0,1)$, where $C_{\alpha}=1/\sqrt{\alpha}$ and $\bb{V}[I_{\rm{MLDLQ}}]$ is the estimator variance. The randomizations of the outer samples are independent for each level; thus, the total variance of the MLDLQMC estimator is decomposed into the variances of the individual levels, that is,
\begin{align}\label{eq:variance.levels}
    \bb{V}\left[I_{\rm{MLDLQ}}\right]={}& \bb{V}\left[I_{\rm{MLDLQ}_0^0}\right]+\sum_{\ell=1}^{L-1}\bb{V}\left[I_{\rm{MLDLQ}_{\ell}^{\ell}}-I_{\rm{MLDLQ}_{\ell}^{\ell-1}}\right]+\bb{V}\left[I_{{\rm{MLDLQ}}_{L}^{L}}-I_{{\rm{MLDLQ}}_{L}^{L-1}}\right].
\end{align}
The terms appearing in \eqref{eq:variance.levels} are bounded individually in the following Propositions~\ref{prop:V.MLrQ-Q.0}-\ref{prop:V.ell.MLQ.h.L}, leading to an upper bound of the total variance of the MLDLQMC estimator.
\begin{asu}\label{asu:boundary.0}
    Let $\varphi_{h_{0},M_{0}}^{(0)}$ be the auxiliary function
    \begin{equation}\label{eq:varphi.h.M.0}
        \varphi_{h_{0},M_{0}}^{(0)}(\cdot)\coloneqq f\left(\frac{1}{M_{0}}\sum_{m=1}^{M_{0}}g_{h_{0}}(\cdot,\bs{x}^{(m)})\right).
    \end{equation}
    Moreover, let $h_{0}>0$, $M_{0}\geq 1$ and $\bs{x}^{(m)}$ be as in~\eqref{eq:x.ell.h.det}, where $1\leq m\leq M_{0}$. Then, assume that $\varphi_{h_{0},M_{0}}^{(0)}$ satisfies Assumption~\ref{eq:boundary.growth}, that is, there exists $0<b^{(0)}<\infty$ and $A_i^{(0)}>0$ for $1\leq i\leq d_1$, where
\begin{equation}\label{eq:A.max.0}
A_{\rm{max}}^{(0)}\coloneqq\max_{1\leq i\leq d_1}A_i^{(0)}<\frac{1}{2},
\end{equation}
such that 
    \begin{equation}
    \left|\left(\prod_{j\in u}\frac{\partial^{u}}{\partial y_j}\right)\varphi_{h_{0},M_{0}}^{(0)}(\bs{y})\right|\leq b^{(0)}\prod_{i=1}^{d_1}\min (y_i,1-y_i)^{-A^{(0)}_i-\mathds{1}_{\{i\in u\}}}
\end{equation}
for all $\bs{y}=(y_1,\ldots,y_{d_1})\in(0,1)^{d_1}$ and all $u\subseteq\{1,\ldots,{d_1}\}$.
\end{asu}
The zeroth level of the MLDLQMC estimator may be regarded as an rQMC estimator~\eqref{eq:randomized.quasi.monte.carlo}, with $N=N_0$, for a particular choice of integrand. In view of Assumption~\ref{asu:boundary.0}, the bound on the variance contribution of the zeroth level of the MLDLQMC estimator follows immediately from Proposition~\ref{cor:he} as stated next.
\begin{prop}[Variance of level 0 of the MLDLQMC estimator with inexact sampling]\label{prop:V.MLrQ-Q.0}
Let $\varphi_{h_0,M_0}^{(0)}$ as in \eqref{eq:varphi.h.M.0} satisfy Assumption~\ref{asu:boundary.0}. The variance contribution of the zeroth level of the MLDLQMC estimator \eqref{eq:MLrQ.h.det} with inexact sampling then has the following upper bound:
\begin{equation}\label{eq:variance.det.0}
  \bb{V}\left[I_{\rm{MLDLQ}_0^0}\right]\leq \frac{(b^{(0)})^2B_{A^{(0)}}^2C_{\epsilon, d_1}^2}{N_0^{2-2\epsilon-2A^{(0)}_{\rm{max}}}},
\end{equation}
for any $\epsilon>0$, where $C_{\epsilon,d_1}\to\infty$ as $\epsilon\to 0$, and $B_{A^{(0)}}\to\infty$ as $\min_{1\leq i\leq d_1}A^{(0)}_i\to 0$.
\end{prop}
The variance contributions of the levels $\ell$, where $1\leq \ell\leq L$, display a different behavior, as the multilevel structure introduces correlations. Nonetheless, these levels are similar in structure to the rQMC estimator~\eqref{eq:randomized.quasi.monte.carlo} with $N=N_{\ell}$, for $1\leq\ell\leq L-1$. The following assumptions are made to analyze these variances:
\begin{asu}\label{asu:boundary.ell.l.o.t.}
    Let $\varphi_{h_{\ell},M_{\ell}}^{(I)}$, for $0\leq \ell\leq L$, be the auxiliary function
    \begin{equation}\label{eq:varphi.h.M.l.o.t.}
        \varphi_{h_{\ell},M_{\ell}}^{(I)}(\cdot)\coloneqq f'\left(\bar{g}(\cdot)\right)\left(\frac{1}{M_{\ell}}\sum_{m=1}^{M_{\ell}}g_{h_{\ell}}(\cdot,\bs{x}^{(m)})-\bar{g}(\cdot)\right).
    \end{equation}
    Moreover, let $h_{\ell}>0$, $M_{\ell}\geq 1$, $\bar{g}$ as in~\eqref{eq:bar.g}, and $\bs{x}^{(m)}$ be as in~\eqref{eq:x.ell.h.det}, where $1\leq m\leq M_{\ell}$, and $0\leq \ell\leq L$. Then, assume that $\varphi_{h_{\ell},M_{\ell}}^{(I)}$, for $0\leq \ell\leq L$, satisfies Assumption~\ref{eq:boundary.growth}, that is,
there exists $0<b^{(I)}<\infty$ and $A_i^{(I)}>0$ for $1\leq i\leq d_1$, where
\begin{equation}\label{eq:A.max.I}
A_{\rm{max}}^{(I)}\coloneqq\max_{1\leq i\leq d_1}A_i^{(I)}<\frac{1}{2},
\end{equation}
such that
    \begin{equation}
    \left|\left(\prod_{j\in u}\frac{\partial^{u}}{\partial y_j}\right)\varphi_{h_{\ell},M_{\ell}}^{(I)}(\bs{y})\right|\leq \left(C_{\epsilon,d_2}M_{\ell}^{-1+\epsilon}+C_{\eta_{{\rm{s}},d_1}}h_{\ell}^{\eta_{{\rm{s}},d_1}}\right)b^{(I)}\prod_{i=1}^{d_1}\min (y_i,1-y_i)^{-A^{(I)}_i-\mathds{1}_{\{i\in u\}}}
\end{equation}
for all $\bs{y}=(y_1,\ldots,y_{d_1})\in(0,1)^{d_1}$ and all $u\subseteq\{1,\ldots,{d_1}\}$, where $0<\eta_{{\rm{s}},d_1}\leq \eta_{\rm{w}}$ is the $h-$convergence rate in a strong sense and $C_{\eta_{{\rm{s}},d_1}}>0$ is independent of $h$.
\end{asu}
\begin{asu}\label{asu:boundary.ell}
    Let $\varphi_{h_{\ell},M_{\ell}}^{(II)}$, for $0\leq \ell\leq L$, be the auxiliary function
    \begin{equation}\label{eq:varphi.h.M}
        \varphi_{h_{\ell},M_{\ell}}^{(II)}(\cdot)\coloneqq\left(\frac{1}{M_{\ell}}\sum_{m=1}^{M_{\ell}}g_{h_{\ell}}(\cdot,\bs{x}^{(m)})-\bar{g}(\cdot)\right)^2\int_0^1f''\left(\bar{g}(\cdot)+s \left(\frac{1}{M_{\ell}}\sum_{m=1}^{M_{\ell}}g_{h_{\ell}}(\cdot,\bs{x}^{(m)})-\bar{g}(\cdot)\right)\right)(1-s)\di{}s.
    \end{equation}
    Moreover, let $h_{\ell}>0$, $M_{\ell}\geq 1$, $\bar{g}$ as in~\eqref{eq:bar.g}, and $\bs{x}^{(m)}$ be as in~\eqref{eq:x.ell.h.det}, where $1\leq m\leq M_{\ell}$, and $0\leq \ell\leq L$. Then, assume that $\varphi_{h_{\ell},M_{\ell}}^{(II)}$, for $0\leq \ell\leq L$, satisfies Assumption~\ref{eq:boundary.growth}, that is, there exists $0<b^{(II)}<\infty$ and $A_i^{(II)}>0$ for $1\leq i\leq d_1$, where
\begin{equation}\label{eq:A.max.II}
A_{\rm{max}}^{(II)}\coloneqq\max_{1\leq i\leq d_1}A_i^{(II)}<\frac{1}{2},
\end{equation}
such that
    \begin{equation}
    \left|\left(\prod_{j\in u}\frac{\partial^{u}}{\partial y_j}\right)\varphi_{h_{\ell},M_{\ell}}^{(II)}(\bs{y})\right|\leq \left(C_{\epsilon,d_2}^2M_{\ell}^{-2+2\epsilon}+C_{\eta_{{\rm{s}},d_1}}^2h_{\ell}^{2\eta_{{\rm{s}},d_1}}\right)b^{(II)}\prod_{i=1}^{d_1}\min (y_i,1-y_i)^{-A^{(II)}_i-\mathds{1}_{\{i\in u\}}}
\end{equation}
for all $\bs{y}=(y_1,\ldots,y_{d_1})\in(0,1)^{d_1}$ and all $u\subseteq\{1,\ldots,{d_1}\}$.
\end{asu}
\begin{prop}[Variance of level $\ell$ for $1\leq\ell\leq L-1$ of the MLDLQMC estimator with inexact sampling]\label{prop:V.ell.MLQ.h.det}
Let $\varphi_{h_{\ell},M_{\ell}}^{(I)}$ as in \eqref{eq:varphi.h.M.l.o.t.}, where $0\leq \ell\leq L-1$, satisfy Assumption~\ref{asu:boundary.ell.l.o.t.} and  $\varphi_{h_{\ell},M_{\ell}}^{(II)}$ as in \eqref{eq:varphi.h.M}, where $0\leq \ell\leq L-1$, satisfy Assumption~\ref{asu:boundary.ell}. The variance contribution of the $\ell$th level of the MLDLQMC estimator \eqref{eq:MLrQ.h.det} with inexact sampling, where $1\leq \ell\leq L-1$, then has the following upper bound:
\begin{align}\label{eq:variance.ell.h.det}
        \bb{V}\left[I_{\rm{MLDLQ}_{\ell}^{\ell}}-I_{\rm{MLDLQ}_{\ell}^{\ell-1}}\right]\leq{}&\frac{\left(\frac{C_{\epsilon,d_2}^2}{M_{\ell}^{2-2\epsilon}}+C_{\eta_{{\rm{s}},d_1}}^2h_{\ell}^{2\eta_{{\rm{s}},d_1}}\right)(b^{(I)})^2B_{A^{(I)}}^2C_{\epsilon,d_1}^2}{N_{\ell}^{2-2\epsilon-2A^{(I)}_{\rm{max}}}}\nonumber\\
        {}&+\frac{\left(\frac{C_{\epsilon,d_2}^4}{M_{\ell}^{4-4\epsilon}}+C_{\eta_{{\rm{s}},d_1}}^4h_{\ell}^{4\eta_{{\rm{s}},d_1}}\right)(b^{(II)})^2B_{A^{(II)}}^2C_{\epsilon,d_1}^2}{N_{\ell}^{2-2\epsilon-2A^{(II)}_{\rm{max}}}},
    \end{align}
    for any $\epsilon>0$, where  $B_{A^{(I)}}\to\infty$ as $\min_{1\leq i\leq d_1}A^{(I)}_i\to 0$, $B_{A^{(II)}}\to\infty$ as $\min_{1\leq i\leq d_1}A^{(II)}_i\to 0$, and $C_{\epsilon,d_1}, C_{\epsilon,d_2}\to\infty$ as $\epsilon\to 0$.
\end{prop}

\begin{proof}
The following derivations were performed for the choices of $M_{\ell}$ and $h_{\ell}$ detailed in~\eqref{eq:M.ell.shape} and~\eqref{eq:h.ell.shape}, where $0\leq\ell\leq L-1$. We omitted the index $\ell\in\{0,\ldots,L-1\}$ for simplicity. The following Taylor expansion is considered for $1\leq n\leq N$:
\begin{align}\label{eq:Taylor.IM.2}
  f\left(\frac{1}{M}\sum_{m=1}^Mg_h(\bs{y}^{(n)},\bs{x}^{(m)})\right)={}&f\left(\bar{g}(\bs{y}^{(n)})\right)+ f'\left(\bar{g}(\bs{y}^{(n)})\right)\left(\frac{1}{M}\sum_{m=1}^Mg_h(\bs{y}^{(n)},\bs{x}^{(m)})-\bar{g}(\bs{y}^{(n)})\right)\nonumber\\
  {}&+ \left(\frac{1}{M}\sum_{m=1}^Mg_h(\bs{y}^{(n)},\bs{x}^{(m)})-\bar{g}(\bs{y}^{(n)})\right)^2\nonumber\\
  {}&\quad\times\int_0^1f''\left(\bar{g}(\bs{y}^{(n)})+s \left(\frac{1}{M}\sum_{m=1}^Mg_h(\bs{y}^{(n)},\bs{x}^{(m)})-\bar{g}(\bs{y}^{(n)})\right)\right)(1-s)\di{}s,
\end{align}
and the expansion:
\begin{align}\label{eq:Taylor.IM.2.I}
  f\left(\frac{1}{\frac{M}{2}}\sum_{m=1}^{\frac{M}{2}}g_{2h}(\bs{y}^{(n)},\bs{x}^{(m)})\right)={}&f\left(\bar{g}(\bs{y}^{(n)})\right) + f'\left(\bar{g}(\bs{y}^{(n)})\right)\left(\frac{1}{\frac{M}{2}}\sum_{m=1}^{\frac{M}{2}}g_{2h}(\bs{y}^{(n)},\bs{x}^{(m)})-\bar{g}(\bs{y}^{(n)})\right)\nonumber\\
  {}&+ \left(\frac{1}{\frac{M}{2}}\sum_{m=1}^{\frac{M}{2}}g_{2h}(\bs{y}^{(n)},\bs{x}^{(m)})-\bar{g}(\bs{y}^{(n)})\right)^2\nonumber\\
  {}&\quad\times\int_0^1f''\left(\bar{g}(\bs{y}^{(n)})+s \left(\frac{1}{\frac{M}{2}}\sum_{m=1}^{\frac{M}{2}}g_{2h}(\bs{y}^{(n)},\bs{x}^{(m)})-\bar{g}(\bs{y}^{(n)})\right)\right)(1-s)\di{}s.
\end{align}
Substituting the Taylor expansions \eqref{eq:Taylor.IM.2} and \eqref{eq:Taylor.IM.2.I} into the variance results in the following:
\begin{align}\label{eq:Taylor.variance.ell}
    {}&\bb{V}\left[\frac{1}{N}\sum_{n=1}^{N}f\left(\frac{1}{M}\sum_{m=1}^{M}g_{h}\left(\bs{y}^{(n)},\bs{x}^{(m)}\right)\right)-f\left(\frac{1}{\frac{M}{2}}\sum_{m=1}^{\frac{M}{2}}g_{2h}\left(\bs{y}^{(n)},\bs{x}^{(m)}\right)\right)\right]\nonumber\\
    {}&=\bb{V}\left[\frac{1}{N}\sum_{n=1}^{N}f'\left(\bar{g}(\bs{y}^{(n)})\right)\left(\frac{1}{M}\sum_{m=1}^{M}g_{h}(\bs{y}^{(n)},\bs{x}^{(m)})-\bar{g}(\bs{y}^{(n)})\right)\right.\nonumber\\
    {}&\quad\quad-\frac{1}{N}\sum_{n=1}^{N}f'\left(\bar{g}(\bs{y}^{(n)})\right)\left(\frac{1}{\frac{M}{2}}\sum_{m=1}^{\frac{M}{2}}g_{2h}(\bs{y}^{(n)},\bs{x}^{(m)})-\bar{g}(\bs{y}^{(n)})\right)\nonumber\\
    {}&\quad\quad+\frac{1}{N}\sum_{n=1}^{N}\left(\frac{1}{M}\sum_{m=1}^Mg_h(\bs{y}^{(n)},\bs{x}^{(m)})-\bar{g}(\bs{y}^{(n)})\right)^2\nonumber\\
    {}&\quad\quad\quad\times\int_0^1f''\left(\bar{g}(\bs{y}^{(n)})+s \left(\frac{1}{M}\sum_{m=1}^Mg_h(\bs{y}^{(n)},\bs{x}^{(m)})-\bar{g}(\bs{y}^{(n)})\right)\right)(1-s)\di{}s\nonumber\\
    {}&\quad\quad-\frac{1}{N}\sum_{n=1}^{N}\left(\frac{1}{\frac{M}{2}}\sum_{m=1}^{\frac{M}{2}}g_{2h}(\bs{y}^{(n)},\bs{x}^{(m)})-\bar{g}(\bs{y}^{(n)})\right)^2\nonumber\\
    {}&\left.\quad\quad\quad\times\int_0^1f''\left(\bar{g}(\bs{y}^{(n)})+s \left(\frac{1}{\frac{M}{2}}\sum_{m=1}^{\frac{M}{2}}g_{2h}(\bs{y}^{(n)},\bs{x}^{(m)})-\bar{g}(\bs{y}^{(n)})\right)\right)(1-s)\di{}s\right].
    \end{align}
    Applying the Cauchy--Schwarz inequality (see \cite[Appendix~A]{Bar23}) to the above variance yields the following:
    \begin{align}\label{eq:var.ell.split}
    {}&\bb{V}\left[\frac{1}{N}\sum_{n=1}^{N}f\left(\frac{1}{M}\sum_{m=1}^{M}g_{h}\left(\bs{y}^{(n)},\bs{x}^{(m)}\right)\right)-f\left(\frac{1}{\frac{M}{2}}\sum_{m=1}^{\frac{M}{2}}g_{2h}\left(\bs{y}^{(n)},\bs{x}^{(m)}\right)\right)\right]\nonumber\\
    {}&\leq4\bb{V}\left[\frac{1}{N}\sum_{n=1}^{N}f'\left(\bar{g}(\bs{y}^{(n)})\right)\left(\frac{1}{M}\sum_{m=1}^{M}g_{h}(\bs{y}^{(n)},\bs{x}^{(m)})-\bar{g}(\bs{y}^{(n)})\right)\right]\nonumber\\
    {}&+4\bb{V}\left[\frac{1}{N}\sum_{n=1}^{N}f'\left(\bar{g}(\bs{y}^{(n)})\right)\left(\frac{1}{\frac{M}{2}}\sum_{m=1}^{\frac{M}{2}}g_{2h}(\bs{y}^{(n)},\bs{x}^{(m)})-\bar{g}(\bs{y}^{(n)})\right)\right]\nonumber\\
    {}&+4\bb{V}\left[\frac{1}{N}\sum_{n=1}^{N}\left(\frac{1}{M}\sum_{m=1}^Mg_h(\bs{y}^{(n)},\bs{x}^{(m)})-\bar{g}(\bs{y}^{(n)})\right)^2\right.\nonumber\\
    {}&\left.\quad\quad\quad\times\int_0^1f''\left(\bar{g}(\bs{y}^{(n)})+s \left(\frac{1}{M}\sum_{m=1}^Mg_h(\bs{y}^{(n)},\bs{x}^{(m)})-\bar{g}(\bs{y}^{(n)})\right)\right)(1-s)\di{}s\right]\nonumber\\
    {}&+4\bb{V}\left[\frac{1}{N}\sum_{n=1}^{N}\left(\frac{1}{\frac{M}{2}}\sum_{m=1}^{\frac{M}{2}}g_{2h}(\bs{y}^{(n)},\bs{x}^{(m)})-\bar{g}(\bs{y}^{(n)})\right)^2\right.\nonumber\\
    {}&\left.\quad\quad\quad\times\int_0^1f''\left(\bar{g}(\bs{y}^{(n)})+s \left(\frac{1}{\frac{M}{2}}\sum_{m=1}^{\frac{M}{2}}g_{2h}(\bs{y}^{(n)},\bs{x}^{(m)})-\bar{g}(\bs{y}^{(n)})\right)\right)(1-s)\di{}s\right].
\end{align}
Taking
\begin{equation}
    \varphi_{h,M}^{(I)}(\cdot)\coloneqq f'\left(\bar{g}(\cdot)\right)\left(\frac{1}{M}\sum_{m=1}^{M}g_{h}(\cdot,\bs{x}^{(m)})-\bar{g}(\cdot)\right),
\end{equation}
it follows from Assumption~\ref{asu:boundary.ell.l.o.t.} that
\begin{align}
    {}&\left|\left(\prod_{j\in u}\frac{\partial^{u}}{\partial y_j}\right)\varphi_{h,M}^{(I)}(\bs{y})\right|\leq \left(C_{\epsilon,d_2}M^{-1+\epsilon}+C_{\eta_{{\rm{s}},d_1}}h^{\eta_{{\rm{s}},d_1}}\right)\tilde{b}^{(I)}\prod_{i=1}^{d_1}\min (y_i,1-y_i)^{-A^{(I)}_i-\mathds{1}_{\{i\in u\}}}
\end{align}
for all $\bs{y}=(y_1,\ldots,y_{d_1})\in(0,1)^{d_1}$ and all $u\subseteq\{1,\ldots,{d_1}\}$. Proposition~\ref{cor:he} then provides the following bound for the first term in~\eqref{eq:var.ell.split}:
\begin{align}
    4\bb{V}\left[\frac{1}{N}\sum_{n=1}^{N}\varphi_{h,M}^{(I)}(\bs{y}^{(n)})\right]{}&\leq4\left(C_{\epsilon,d_2}M^{-1+\epsilon}+C_{\eta_{{\rm{s}},d_1}}h^{\eta_{{\rm{s}},d_1}}\right)^2(\tilde{b}^{(I)})^2B^2_{A^{(I)}}C_{\epsilon,d_1}^2N^{-2+2\epsilon+2A^{(I)}_{\rm{max}}},\nonumber\\
    {}&\leq8\left(C_{\epsilon,d_2}^2M^{-2+2\epsilon}+C_{\eta_{{\rm{s}},d_1}}^2h^{2\eta_{{\rm{s}},d_1}}\right)(\tilde{b}^{(I)})^2B^2_{A^{(I)}}C_{\epsilon,d_1}^2N^{-2+2\epsilon+2A^{(I)}_{\rm{max}}},
\end{align}
for all $\epsilon>0$, where $B_{A^{(I)}}\to\infty$ as $\min_{1\leq i\leq d_1}A^{(I)}_i\to 0$ and $C_{\epsilon,d_1}, C_{\epsilon,d_2}\to\infty$ as $\epsilon\to 0$. Taking
\begin{equation}
    \varphi_{2h,M/2}^{(I)}(\cdot)\coloneqq f'\left(\bar{g}(\cdot)\right)\left(\frac{1}{\frac{M}{2}}\sum_{m=1}^{\frac{M}{2}}g_{2h}(\cdot,\bs{x}^{(m)})-\bar{g}(\cdot)\right),
\end{equation}
a similar bound follows for the second term in~\eqref{eq:var.ell.split}:
\begin{align}
    4\bb{V}\left[\frac{1}{N}\sum_{n=1}^{N}\varphi_{2h,M/2}^{(I)}(\bs{y}^{(n)})\right]{}&\leq8\left(C_{\epsilon,d_2}^22^{2-2\epsilon}M^{-2+2\epsilon}+C_{\eta_{{\rm{s}},d_1}}^22^{2\eta_{{\rm{s}},d_1}}h^{2\eta_{{\rm{s}},d_1}}\right)(\tilde{b}^{(I)})^2B^2_{A^{(I)}}C_{\epsilon,d_1}^2N^{-2+2\epsilon+2A^{(I)}_{\rm{max}}}.
\end{align}
Next, taking
\begin{equation}
    \varphi_{h,M}^{(II)}(\cdot)\coloneqq\left(\frac{1}{M}\sum_{m=1}^Mg_h(\cdot,\bs{x}^{(m)})-\bar{g}(\cdot)\right)^2\int_0^1f''\left(\bar{g}(\cdot)+s \left(\frac{1}{M}\sum_{m=1}^Mg_h(\cdot,\bs{x}^{(m)})-\bar{g}(\cdot)\right)\right)(1-s)\di{}s,
\end{equation}
it follows from Assumption~\ref{asu:boundary.ell} that
\begin{equation}
    \left|\left(\prod_{j\in u}\frac{\partial^{u}}{\partial y_j}\right)\varphi_{h,M}^{(II)}(\bs{y})\right|\leq \left(C_{\epsilon,d_2}^2M^{-2+2\epsilon}+C_{\eta_{{\rm{s}},d_1}}^2h^{2\eta_{{\rm{s}},d_1}}\right)\tilde{b}^{(II)}\prod_{i=1}^{d_1}\min (y_i,1-y_i)^{-A^{(II)}}_i-\mathds{1}_{\{i\in u\}}
\end{equation}
for all $\bs{y}=(y_1,\ldots,y_{d_1})\in(0,1)^{d_1}$ and all $u\subseteq\{1,\ldots,{d_1}\}$. Proposition~\ref{cor:he} then provides the following bound for the third term in~\eqref{eq:var.ell.split}:
\begin{align}
    4\bb{V}\left[\frac{1}{N}\sum_{n=1}^{N}\varphi_{h,M}^{(II)}(\bs{y}^{(n)})\right]\leq8\left(C_{\epsilon,d_2}^4M^{-4+4\epsilon}+C_{\eta_{{\rm{s}},d_1}}^4h^{4\eta_{{\rm{s}},d_1}}\right)(\tilde{b}^{(II)})^2B^2_{A^{(II)}}C_{\epsilon,d_1}^2N^{-2+2\epsilon+2A^{(II)}_{\rm{max}}},
\end{align}
and finally, taking
\begin{equation}
    \varphi_{2h,M/2}^{(II)}(\cdot)\coloneqq\left(\frac{1}{\frac{M}{2}}\sum_{m=1}^{\frac{M}{2}}g_{2h}(\cdot,\bs{x}^{(m)})-\bar{g}(\cdot)\right)^2\int_0^1f''\left(\bar{g}(\cdot)+s \left(\frac{1}{\frac{M}{2}}\sum_{m=1}^{\frac{M}{2}}g_{2h}(\cdot,\bs{x}^{(m)})-\bar{g}(\cdot)\right)\right)(1-s)\di{}s,
\end{equation}
a similar bound applies to the fourth term in~\eqref{eq:var.ell.split}:
\begin{align}
    4\bb{V}\left[\frac{1}{N}\sum_{n=1}^{N}\varphi_{2h,M/2}^{(II)}(\bs{y}^{(n)})\right]\leq8\left(C_{\epsilon,d_2}^42^{4-4\epsilon}M^{-4+4\epsilon}+C_{\eta_{{\rm{s}},d_1}}^42^{4\eta_{{\rm{s}},d_1}}h^{4\eta_{{\rm{s}},d_1}}\right)(\tilde{b}^{(II)})^2B^2_{A^{(II)}}C_{\epsilon,d_1}^2N^{-2+2\epsilon+2A^{(II)}_{\rm{max}}},
    \end{align}
    for all $\epsilon>0$, where $B_{A^{(II)}}\to\infty$ as $\min_{1\leq i \leq d_1}A^{(II)}_i\to 0$, and $C_{\epsilon,d_1}, C_{\epsilon,d_2}\to\infty$ as $\epsilon\to 0$.
    Combining the above results yields the following:
    \begin{align}
        \bb{V}\left[\frac{1}{N}\sum_{n=1}^{N}f\left(\frac{1}{M}\sum_{m=1}^{M}g_{h}\left(\bs{y}^{(n)},\bs{x}^{(m)}\right)\right)-\right.&{}\left.f\left(\frac{1}{\frac{M}{2}}\sum_{m=1}^{\frac{M}{2}}g_{2h}\left(\bs{y}^{(n)},\bs{x}^{(m)}\right)\right)\right]\nonumber\\
        \leq{}&\frac{\left(\frac{C_{\epsilon,d_2}^2}{M^{2-2\epsilon}}+C_{\eta_{{\rm{s}},d_1}}^2h^{2\eta_{{\rm{s}},d_1}}\right)(b^{(I)})^2B^2_{A^{(I)}}C_{\epsilon,d_1}^2}{N^{2-2\epsilon-2A^{(I)}_{\rm{max}}}}\nonumber\\
        {}&+\frac{\left(\frac{C_{\epsilon,d_2}^4}{M^{4-4\epsilon}}+C_{\eta_{{\rm{s}},d_1}}^4h^{4\eta_{{\rm{s}},d_1}}\right)(b^{(II)})^2B^2_{A^{(II)}}C_{\epsilon,d_1}^2}{N^{2-2\epsilon-2A^{(II)}_{\rm{max}}}},
    \end{align}
    where
    \begin{equation}
        (b^{(I)})^2\coloneqq \frac{8(\tilde{b}^{(I)})^2}{\max\{2^{2-2\epsilon},2^{2\eta_{{\rm{s}},d_1}}\}+1}
    \end{equation}
    and
    \begin{equation}
        (b^{(II)})^2\coloneqq \frac{8(\tilde{b}^{(II)})^2}{\max\{2^{4-4\epsilon},2^{4\eta_{{\rm{s}},d_1}}\}+1},
    \end{equation}
    for all $\epsilon>0$, where $B_{A^{(I)}}\to\infty$ as $\min_{1\leq i\leq d_1}A^{(I)}_i\to 0$, $B_{A^{(II)}}\to\infty$ as $\min_{1\leq i\leq d_1}A^{(II)}_i\to 0$, and $C_{\epsilon,d_1}, C_{\epsilon,d_2}\to\infty$ as $\epsilon\to 0$. 
\end{proof}
A base other than two in~\eqref{eq:M.ell.shape} and~\eqref{eq:h.ell.shape} results in different multiplicative constants $b^{(I)}$ and $b^{(II)}$.
\begin{rmk}[Difference between terms in the variance bound]
     If it holds that $A^{(I)}_{\rm{max}}\leq A^{(II)}_{\rm{max}}$, then the second term in the variance bound~\eqref{eq:variance.ell.h.det} is of higher order compared to the first term.
\end{rmk}
For the bound of the variance of the MLDLQMC estimator at the final level $L$, the following additional assumptions are made regarding the variance from the inner randomization:
\begin{asu}\label{asu:boundary.L.l.o.t.}
    Let $\varphi_{h_{L},M_{L}}^{(III)}$ be an auxiliary function, where
    \begin{equation}\label{eq:varphi.h.M.L.l.o.t.}
        \varphi_{h_{L},M_{L}}^{(III)}(\cdot)\coloneqq f'\left(\bar{g}_{h_{L}}(\cdot)\right)\left(\frac{1}{M_{L}}\sum_{m=1}^{M_{L}}g_{h_{L}}(\cdot,\bs{x}^{(m)})-\bar{g}_h(\cdot)\right).
    \end{equation}
    Moreover, let  $h_{L}>0$, $M_{L}\geq 1$,
    \begin{equation}\label{eq:g.h.bar.a}
        \bar{g}_{h_{L}}(\cdot)\coloneqq\int_{[0,1]^{d_2}}g_{h_{L}}(\cdot,\bs{x})\di{}\bs{x},
    \end{equation}
    and $\bs{x}^{(m)}$ be as in~\eqref{eq:x.ell.h.det}, where $1\leq m\leq M_{L}$. Then, assume that there exists $0<b^{(III)}<\infty$ and $A_i^{(III)}>0$ for $1\leq i\leq d_1$, where
\begin{equation}\label{eq:A.max.III}
A_{\rm{max}}^{(III)}\coloneqq\max_{1\leq i\leq d_1}A_i^{(III)}<\frac{1}{2},
\end{equation}
such that
    \begin{equation}
    \left|\varphi_{h_{L},M_{L}}^{(III)}(\bs{y})\right|\leq C_{\epsilon,d_2}M_{L}^{-1+\epsilon}b^{(III)}\prod_{i=1}^{d_1}\min (y_i,1-y_i)^{-A^{(III)}_i}
\end{equation}
for all $\bs{y}=(y_1,\ldots,y_{d_1})\in(0,1)^{d_1}$,
 and all $u\subseteq\{1,\ldots,{d_1}\}$, where $b^{(III)}$ implicitly depends on $h_L$, and $\eta_{{\rm{s}},d_1}$.
\end{asu}
\begin{asu}\label{asu:boundary.L}
    Let $\varphi_{h_{L},M_{L}}^{(IV)}$ be an auxiliary function, where
    \begin{equation}\label{eq:varphi.h.M.L}
        \varphi_{h_{L},M_{L}}^{(IV)}(\cdot)\coloneqq\left(\frac{1}{M_{L}}\sum_{m=1}^{M_{L}}g_{h_{L}}(\cdot,\bs{x}^{(m)})-\bar{g}_{h_{L}}(\cdot)\right)^2\int_0^1f''\left(\bar{g}_{h_{L}}(\cdot)+s \left(\frac{1}{M_{L}}\sum_{m=1}^{M_{L}}g_{h_{L}}(\cdot,\bs{x}^{(m)})-\bar{g}_{h_{L}}(\cdot)\right)\right)(1-s)\di{}s.
    \end{equation}
    Moreover, let $h_{L}>0$, $M_{L}\geq 1$, $\bar{g}_{h_{L}}$ as in~\eqref{eq:g.h.bar.a}, and $\bs{x}^{(m)}$ be as in~\eqref{eq:x.ell.h.det}, where $1\leq m\leq M_{L}$. Then, assume that there exists $0<b^{(IV)}<\infty$ and $A_i^{(IV)}>0$ for $1\leq i\leq d_1$, where
\begin{equation}\label{eq:A.max.IV}
A_{\rm{max}}^{(IV)}\coloneqq\max_{1\leq i\leq d_1}A_i^{(IV)}<\frac{1}{2},
\end{equation}
such that
    \begin{equation}
    \left|\varphi_{h_{L},M_{L}}^{(IV)}(\bs{y})\right|\leq C_{\epsilon,d_2}^2M_{L}^{-2+2\epsilon}b^{(IV)}\prod_{i=1}^{d_1}\min (y_i,1-y_i)^{-A^{(IV)}_i}
\end{equation}
for all $\bs{y}=(y_1,\ldots,y_{d_1})\in(0,1)^{d_1}$ and all $u\subseteq\{1,\ldots,{d_1}\}$, where $b^{(IV)}$ implicitly depends on $h_{L}$ and $\eta_{{\rm{s}},d_1}$.
\end{asu}
\begin{prop}[Variance of level $L$ of the MLDLQMC estimator with inexact sampling]\label{prop:V.ell.MLQ.h.L}
Let $\varphi_{h_{L-1},M_{L-1}}^{(I)}$ and $\varphi_{h_{L},M_{L}}^{(I)}$ as in \eqref{eq:varphi.h.M.l.o.t.} satisfy Assumption~\ref{asu:boundary.ell.l.o.t.},  $\varphi_{h_{L-1},M_{L-1}}^{(II)}$ and $\varphi_{h_{L},M_{L}}^{(II)}$ as in \eqref{eq:varphi.h.M} satisfy Assumption~\ref{asu:boundary.ell}, $\varphi_{h_L,M_L}^{(III)}$ as in \eqref{eq:varphi.h.M.L.l.o.t.} satisfy Assumption~\ref{asu:boundary.L.l.o.t.}, and $\varphi_{h_L,M_L}^{(IV)}$ as in \eqref{eq:varphi.h.M.L} satisfy Assumption~\ref{asu:boundary.L}.
The variance contribution of the $L$th level of the MLDLQMC estimator \eqref{eq:MLrQ.h.det} with inexact sampling then has the following upper bound:
\begin{align}\label{eq:variance.L.h.det}
        \bb{V}\left[I_{{\rm{MLDLQ}}_{L}^{L}}-I_{{\rm{MLDLQ}}_{L}^{L-1}}\right]\leq{}&\frac{\left(\frac{C_{\epsilon,d_2}^2}{M_{L}^{2-2\epsilon}}+C_{\eta_{{\rm{s}},d_1}}^2h_{L}^{2\eta_{{\rm{s}},d_1}}\right)(b^{(I)})^2B_{A^{(I)}}^2C_{\epsilon,d_1}^2}{N_{L}^{2-2\epsilon-2A^{(I)}_{\rm{max}}}}\nonumber\\
        {}&+\frac{\left(\frac{C_{\epsilon,d_2}^4}{M_{L}^{4-4\epsilon}}+C_{\eta_{{\rm{s}},d_1}}^4h_{L}^{4\eta_{{\rm{s}},d_1}}\right)(b^{(II)})^2B_{A^{(II)}}^2C_{\epsilon,d_1}^2}{N_{L}^{2-2\epsilon-2A^{(II)}_{\rm{max}}}}\nonumber\\
        {}&+\frac{C_{\epsilon,d_2}^2(b^{(III)})^2B_{A^{(III)}}^2}{M_L^{2-2\epsilon}}+\frac{C_{\epsilon,d_2}^4(b^{(IV)})^2B_{A^{(IV)}}^2}{M_L^{4-4\epsilon}},
    \end{align}
    for any $\epsilon>0$, where $B_{A^{(I)}}\to\infty$ as $\min_{1\leq i\leq d_1}A^{(I)}_i\to 0$, $B_{A^{(II)}}\to\infty$ as $\min_{1\leq i\leq d_1}A^{(II)}_i\to 0$, $B_{A^{(III)}}\to\infty$ as $\min_{1\leq i\leq d_1}A^{(III)}_i\to 0$, $B_{A^{(IV)}}\to\infty$ as $\min_{1\leq i\leq d_1}A^{(IV)}_i\to 0$, and $C_{\epsilon,d_1}, C_{\epsilon,d_2}\to\infty$ as $\epsilon\to 0$.
\end{prop}

\begin{proof}
The following derivations were performed for the specific choices of $M_{L-1}$, $M_{L}$, $h_{L-1}$, and $h_{L}$ detailed in~\eqref{eq:M.ell.shape} and~\eqref{eq:h.ell.shape}. We omitted the index $L$ for simplicity. The law of total variance is applied to obtain
\begin{align}\label{eq:V.L.total}
    \bb{V}{}&\left[\frac{1}{N}\sum_{n=1}^{N}\left(f\left(\frac{1}{M}\sum_{m=1}^{M}g_{h}\left(\bs{\tau}(\bs{u}^{(n)}),\bs{\rho}(\bs{t}^{(m)})\right)\right)-f\left(\frac{1}{\frac{M}{2}}\sum_{m=1}^{\frac{M}{2}}g_{2h}\left(\bs{\tau}(\bs{u}^{(n)}),\bs{t}^{(m)}\right)\right)\right)\right]\nonumber\\
    {}&\leq \bb{V}\left[\bb{E}\left[\frac{1}{N}\sum_{n=1}^{N}\left(f\left(\frac{1}{M}\sum_{m=1}^{M}g_{h}\left(\bs{\tau}(\bs{u}^{(n)}),\bs{\rho}(\bs{t}^{(m)})\right)\right)-f\left(\frac{1}{\frac{M}{2}}\sum_{m=1}^{\frac{M}{2}}g_{2h}\left(\bs{\tau}(\bs{u}^{(n)}),\bs{t}^{(m)}\right)\right)\right)\Bigg|\bs{\rho}\right]\right]\nonumber\\
    {}&\quad+\bb{E}\left[\bb{V}\left[\frac{1}{N}\sum_{n=1}^{N}\left(f\left(\frac{1}{M}\sum_{m=1}^{M}g_{h}\left(\bs{\tau}(\bs{u}^{(n)}),\bs{\rho}(\bs{t}^{(m)})\right)\right)-f\left(\frac{1}{\frac{M}{2}}\sum_{m=1}^{\frac{M}{2}}g_{2h}\left(\bs{\tau}(\bs{u}^{(n)}),\bs{t}^{(m)}\right)\right)\right)\Bigg|\bs{\rho}\right]\right],
\end{align}
where the randomization $\bs{\tau}$ maps the low-discrepancy points $\bs{u}^{(n)}$, where $1\leq n\leq N$, to a low-discrepancy point set with probability one, and
the randomization $\bs{\rho}$ maps the low-discrepancy points $\bs{t}^{(m)}$, where $1\leq m\leq M$, to a low-discrepancy point set with probability one; thus, the variance bound of the second term in~\eqref{eq:V.L.total} is analogous to the bound for level $\ell$ where $1\leq\ell\leq L-1$. For the first term in~\eqref{eq:V.L.total}, it holds that
\begin{align}
    \bb{V}{}&\left[\bb{E}\left[\frac{1}{N}\sum_{n=1}^{N}\left(f\left(\frac{1}{M}\sum_{m=1}^{M}g_{h}\left(\bs{\tau}(\bs{u}^{(n)}),\bs{\rho}(\bs{t}^{(m)})\right)\right)-f\left(\frac{1}{\frac{M}{2}}\sum_{m=1}^{\frac{M}{2}}g_{2h}\left(\bs{\tau}(\bs{u}^{(n)}),\bs{t}^{(m)}\right)\right)\right)\Bigg|\bs{\rho}\right]\right]\nonumber\\
    {}&=\bb{V}\left[\int_{[0,1]^{d_1}}\left(f\left(\frac{1}{M}\sum_{m=1}^{M}g_{h}\left(\bs{y},\bs{\rho}(\bs{t}^{(m)})\right)\right)-f\left(\frac{1}{\frac{M}{2}}\sum_{m=1}^{\frac{M}{2}}g_{2h}\left(\bs{y},\bs{t}^{(m)}\right)\right)\right)\di{}\bs{y}\right],\nonumber\\
    {}&=\bb{V}\left[\int_{[0,1]^{d_1}}f\left(\frac{1}{M}\sum_{m=1}^{M}g_{h}\left(\bs{y},\bs{\rho}(\bs{t}^{(m)})\right)\right)\di{}\bs{y}\right].
    \end{align}
    The following Taylor expansion is considered:
\begin{align}\label{eq:Taylor.IM.M}
  f\left(\frac{1}{M}\sum_{m=1}^Mg_h(\bs{y},\bs{x}^{(m)})\right)={}& f\left(\bar{g}_h(\bs{y})\right)+ f'\left(\bar{g}_h(\bs{y})\right)\left(\frac{1}{M}\sum_{m=1}^Mg_h(\bs{y},\bs{x}^{(m)})-\bar{g}_h(\bs{y})\right)\nonumber\\
  {}&+ \left(\frac{1}{M}\sum_{m=1}^Mg_h(\bs{y},\bs{x}^{(m)})-\bar{g}_h(\bs{y})\right)^2\nonumber\\
  {}&\quad\times\int_0^1f''\left(\bar{g}_h(\bs{y})+s \left(\frac{1}{M}\sum_{m=1}^Mg_h(\bs{y},\bs{x}^{(m)})-\bar{g}_h(\bs{y})\right)\right)(1-s)\di{}s,
\end{align} 
to obtain
    \begin{align}\label{eq:V.L.M}
    {}&\bb{V}\left[\int_{[0,1]^{d_1}}f\left(\frac{1}{M}\sum_{m=1}^{M}g_{h}\left(\bs{y},\bs{x}^{(m)}\right)\right)\di{}\bs{y}\right]\nonumber\\
    {}&=\bb{V}\left[\int_{[0,1]^{d_1}}f\left(\bar{g}_h(\bs{y})\right)+f'\left(\bar{g}_h(\bs{y})\right)\left(\frac{1}{M}\sum_{m=1}^{M}g_{h}\left(\bs{y},\bs{x}^{(m)}\right)-\bar{g}_h(\bs{y})\right) \right.\nonumber\\
    {}&\left.\quad\quad\quad+\left(\frac{1}{M}\sum_{m=1}^Mg_h(\bs{y},\bs{x}^{(m)}))-\bar{g}_h(\bs{y})\right)^2\int_0^1f''\left(\bar{g}_h(\bs{y})+s \left(\frac{1}{M}\sum_{m=1}^Mg_h(\bs{y},\bs{x}^{(m)})-\bar{g}_h(\bs{y})\right)\right)(1-s)\di{}s\di{}\bs{y}\right],\nonumber\\
    {}&\leq 2\bb{V}\left[\int_{[0,1]^{d_1}}f'\left(\bar{g}_h(\bs{y})\right)\left(\frac{1}{M}\sum_{m=1}^{M}g_{h}\left(\bs{y},\bs{x}^{(m)}\right)-\bar{g}_h(\bs{y})\right) \di{}\bs{y}\right]\nonumber\\
    {}&\quad+2\bb{V}\left[\int_{[0,1]^{d_1}}\left(\frac{1}{M}\sum_{m=1}^Mg_h(\bs{y},\bs{x}^{(m)})-\bar{g}_h(\bs{y})\right)^2\right.\nonumber\\
    {}&\left.\quad\quad\quad\quad\times\int_0^1f''\left(\bar{g}_h(\bs{y})+s \left(\frac{1}{M}\sum_{m=1}^Mg_h(\bs{y},\bs{x}^{(m)})-\bar{g}_h(\bs{y})\right)\right)(1-s)\di{}s\di{}\bs{y}\right].
\end{align}
Taking
\begin{equation}
    \varphi_{h,M}^{(III)}(\cdot)\coloneqq f'\left(\bar{g}_h(\cdot)\right)\left(\frac{1}{M}\sum_{m=1}^{M}g_{h}\left(\cdot,\bs{x}^{(m)}\right)-\bar{g}_h(\cdot)\right),
\end{equation}
it follows from the Cauchy--Schwarz inequality and Assumption~\ref{asu:boundary.L.l.o.t.} that the first term in \eqref{eq:V.L.M} is bounded as follows:
\begin{align}
    2\bb{V}\left[\int_{[0,1]^{d_1}}\varphi_{h,M}^{(III)}(\bs{y}) \di{}\bs{y}\right]{}&\leq2\bb{E}\left[\left|\int_{[0,1]^{d_1}}\varphi_{h,M}^{(III)}(\bs{y}) \di{}\bs{y}\right|^2\right],\nonumber\\
    {}&\leq 2\left|\int_{[0,1]^{d_1}}C_{\epsilon,d_2}M^{-1+\epsilon}\tilde{b}^{(III)}\prod_{i=1}^{d_1}\min (y_i,1-y_i)^{-A^{(III)}_i} \di{}\bs{y}\right|^2,\nonumber\\
    {}&= 2C_{\epsilon,d_2}^2M^{-2+2\epsilon}(\tilde{b}^{(III)})^2\left|\int_{[0,1]^{d_1}}\prod_{i=1}^{d_1}\min (y_i,1-y_i)^{-A^{(III)}_i} \di{}\bs{y}\right|^2,\nonumber\\
    {}&= \frac{C_{\epsilon,d_2}^2(b^{(III)})^2B_{A^{(III)}}^2}{M^{2-2\epsilon}},
\end{align}
where
\begin{equation}
    (b^{(III)})^2\coloneqq 2(\tilde{b}^{(III)})^2,
\end{equation}
and
\begin{equation}
    B_{A^{(III)}}\coloneqq \int_{[0,1]^{d_1}}\prod_{i=1}^{d_1}\min (y_i,1-y_i)^{-A^{(III)}_i} \di{}\bs{y},
\end{equation}
for any $\epsilon>0$, where $C_{\epsilon,d_2}\to\infty$ as $\epsilon\to 0$ and $B_{A^{(III)}}\to\infty$ as $\min_{1\leq i\leq d_1}A^{(III)}_i\to 0$. Finally, taking
\begin{equation}
    \varphi_{h,M}^{(IV)}(\cdot)\coloneqq\left(\frac{1}{M}\sum_{m=1}^Mg_h(\cdot,\bs{x}^{(m)})-\bar{g}_h(\cdot)\right)^2\int_0^1f''\left(\bar{g}_h(\cdot)+s \left(\frac{1}{M}\sum_{m=1}^Mg_h(\cdot,\bs{x}^{(m)})-\bar{g}_h(\cdot)\right)\right)(1-s)\di{}s,
\end{equation}
it follows from the Cauchy--Schwarz inequality and Assumption~\ref{asu:boundary.L} that the second term in \eqref{eq:V.L.M} is bounded as follows:
\begin{align}
    2\bb{V}\left[\int_{[0,1]^{d_1}}\varphi_{h,M}^{(IV)}(\bs{y})\di{}\bs{y}\right]{}&\leq 2\bb{E}\left[\left|\int_{[0,1]^{d_1}}\varphi_{h,M}^{(IV)}(\bs{y})\di{}\bs{y}\right|^2\right],\nonumber\\
    {}&\leq 2\left|\int_{[0,1]^{d_1}}C_{\epsilon,d_2}^2M^{-2+2\epsilon}\tilde{b}^{(IV)}\prod_{i=1}^{d_1}\min (y_i,1-y_i)^{-A^{(IV)}_i}\di{}\bs{y}\right|^2,\nonumber\\
    {}&= 2C_{\epsilon,d_2}^4M^{-4+4\epsilon}(\tilde{b}^{(IV)})^2\left|\int_{[0,1]^{d_1}}\prod_{i=1}^{d_1}\min (y_i,1-y_i)^{-A^{(IV)}_i}\di{}\bs{y}\right|^2,\nonumber\\
    {}&= \frac{C_{\epsilon,d_2}^4(b^{(IV)})^2B_{A^{(IV)}}^2}{M^{4-4\epsilon}},
\end{align}
where
\begin{equation}
    (b^{(IV)})^2\coloneqq 2(\tilde{b}^{(IV)})^2,
\end{equation}
and
\begin{equation}
    B_{A^{(IV)}}\coloneqq \int_{[0,1]^{d_1}}\prod_{i=1}^{d_1}\min (y_i,1-y_i)^{-A^{(IV)}_i} \di{}\bs{y},
\end{equation}
for any $\epsilon>0$, where $C_{\epsilon,d_2}\to\infty$ as $\epsilon\to 0$ and $B_{A^{(IV)}}\to\infty$ as $\min_{1\leq i\leq d_1}A^{(IV)}_i\to 0$.
\end{proof}
A base other than two in~\eqref{eq:M.ell.shape} and~\eqref{eq:h.ell.shape} results in different multiplicative constants $b^{(I)}$, $b^{(II)}$, $b^{(III)}$, and $b^{(IV)}$.
\begin{rmk}[Optimal choices of $M_{\ell}$ and $h_{\ell}$]\label{rmk:optimal.M.h}
    It follows from Propositions~\ref{prop:B.MLrQ-Q.h}, \ref{prop:V.ell.MLQ.h.det},  and~\ref{prop:V.ell.MLQ.h.L} that the choices (see~\eqref{eq:M.ell.shape} and~\eqref{eq:h.ell.shape})
    \begin{equation}\label{M.ell.2}
    M_{\ell}\coloneqq M_02^{\ell}, \quad 1\leq \ell\leq L,
\end{equation}
and
\begin{equation}
    h_{\ell}\coloneqq h_02^{-\ell}, \quad 1\leq \ell\leq L,
\end{equation}
may be suboptimal if $\eta_{\rm{w}}$ differs from $1-\epsilon$. In that case, we select
    \begin{equation}
    M_{\ell}\coloneqq  M_0\left(2^{\frac{\eta_{\rm{w}}}{1-\epsilon}}\right)^{\ell}, \quad 1\leq \ell\leq L.
\end{equation}
It must be ensured that $M_{\ell}$, where $1\leq\ell\leq L$, are powers of two by rounding up appropriately to maintain the properties of the Sobol' sequence. Alternatively, it is possible to select
    \begin{equation}
    h_{\ell}\coloneqq  h_0\left(2^{\frac{1-\epsilon}{\eta_{\rm{w}}}}\right)^{-\ell}, \quad 1\leq \ell\leq L,
\end{equation}
and $M_{\ell}=M_02^{\ell}$. Both of these options lead to the same overall computational cost asymptotically.
\end{rmk}
\begin{rmk}[Extension of the rDLQMC estimator to the multilevel setting]
    It becomes apparent from Propositions~\ref{prop:V.MLrQ-Q.0}--\ref{prop:V.ell.MLQ.h.L} that the particular structure of the MLDLQMC estimator~\eqref{eq:MLrQ.h.det} (deterministic inner samples at levels $0\leq \ell\leq L-1$, and randomized inner samples at level $L$) is crucial. Adapting the rDLQMC estimator \cite{Bar23} to the multilevel setting would deteriorate the convergence of the variances at levels $0\leq \ell\leq L-1$ to rate $N$ rather than rate $N^{2-2\epsilon-2A_{\rm{max}}}$. Instead, randomizing the inner samples in the last level does not affect the overall estimator complexity and provides a practical error estimate for the bias.
\end{rmk}
\begin{rmk}[Antithetic sampling]
    The coupling in multilevel estimators can be increased by an antithetic construction as in the works \cite{Gil14,God20}; however, such an antithetic construction is not possible for the MLDLQMC estimator because of the combination of deterministic and randomized inner samples.
\end{rmk}
\subsection{Optimal work analysis}
\begin{prop}[Optimal work of the MLDLQMC estimator with inexact sampling]\label{prop:Work.MLQ.h.det}
Given the assumptions of Propositions~\ref{prop:B.MLrQ-Q.h}--\ref{prop:V.ell.MLQ.h.L}, the total computational work of the MLDLQMC estimator~\eqref{eq:MLrQ.h.det} with inexact sampling is given by
\begin{align}\label{eq:work.h.det}
  W_{\rm{MLDLQ}}^{\ast}={}&\cl{O}\left(TOL^{-\left(\frac{1}{1-\epsilon}+\frac{\gamma}{\eta_{\rm{w}}}\right)}\right)\nonumber\\
  {}&+\begin{cases}\cl{O}\left(TOL^{-\frac{1}{1-\epsilon-A_{\rm{max}}}}\right), & \text{if}\quad \frac{\eta_{{\rm{s}},d_1}}{\eta_{\rm{w}}(1-\epsilon-A_{\rm{max}})}>\left(\frac{1}{1-\epsilon}+\frac{\gamma}{\eta_{\rm{w}}}\right),\\
  \cl{O}\left(TOL^{-\frac{1}{1-\epsilon-A_{\rm{max}}}} \log_{2}\left(TOL^{-\max\left\{ \frac{1}{\eta_{\rm{w}}},\frac{1}{1-\epsilon}\right\} } \right)\right), & \text{if}\quad \frac{\eta_{{\rm{s}},d_1}}{\eta_{\rm{w}}(1-\epsilon-A_{\rm{max}})}=\left(\frac{1}{1-\epsilon}+\frac{\gamma}{\eta_{\rm{w}}}\right),\\
  \cl{O}\left(TOL^{-\frac{1}{1-\epsilon-A_{\rm{max}}}+\frac{\eta_{{\rm{s}},d_1}}{\eta_{\rm{w}}(1-\epsilon-A_{\rm{max}})}-\left(\frac{1}{1-\epsilon}+\frac{\gamma}{\eta_{\rm{w}}}\right)}\right), & \text{if}\quad \frac{\eta_{{\rm{s}},d_1}}{\eta_{\rm{w}}(1-\epsilon-A_{\rm{max}})}<\left(\frac{1}{1-\epsilon}+\frac{\gamma}{\eta_{\rm{w}}}\right),\end{cases}
\end{align}
where
\begin{equation}\label{def:A.max}
    A_{\rm{max}}\coloneqq \max\left\{A^{(0)}_{\rm{max}},A^{(I)}_{\rm{max}},A^{(II)}_{\rm{max}} \right\},
\end{equation}
with $A^{(0)}_{\rm{max}}$ as in \eqref{eq:A.max.0}, $A^{(I)}_{\rm{max}}$ as in \eqref{eq:A.max.I}, and $A^{(II)}_{\rm{max}}$ as in \eqref{eq:A.max.II},
for any $\epsilon>0$ as $TOL\to 0$, where the multiplicative constant in the $\cl{O}$ notation approaches infinity as $\epsilon\to 0$.
\end{prop}
\begin{proof}
We optimize the total computational work of the MLDLQMC estimator, that is,
\begin{equation}
    W=\sum_{\ell=0}^LN_{\ell}W_{\ell}+\sum_{\ell=0}^LW_{\ell},
\end{equation}
where
\begin{equation}\label{eq:W.ell}
    W_{\ell}\coloneqq M_{\ell}h_{\ell}^{-\gamma},
\end{equation}
under the constraint that
\begin{equation}
    |I_{\rm{MLDLQ}}-I|\leq TOL,
\end{equation}
where $TOL>0$ is a given error tolerance.
From~\eqref{eq:MLQMC.error} and~\eqref{eq:statistical.error}, the constraint on the total error is split as follows:
\begin{equation}\label{eq:bias.constr}
    |\mathbb{E}[I_{\rm{MLDLQ}}]-I|\leq\frac{1}{2}TOL,
\end{equation}
and
\begin{equation}
    \sqrt{C_{\alpha}\bb{V}\left[I_{\rm{MLDLQ}}\right]}\leq \frac{1}{2}TOL,
\end{equation}
which we rewrite as
\begin{equation}
    C_{\alpha}^2\bb{V}\left[I_{\rm{MLDLQ}}\right]\leq \frac{1}{4}TOL^2.
\end{equation}
The bias constraint in~\eqref{eq:bias.constr} is further split into the discretization bias and the inner sampling bias and we require that each of them is bounded by $TOL/4$.
We substitute the previously derived bounds on the bias and variance and arrive at the optimization problem
\begin{align}\label{eq:Lagrangian.h.det.simplified}
  {}&(L^{\ast},N_0^{\ast},\ldots,N_{L^{\ast}})=\argmin_{(L,N_0,\ldots,N_{L})}\sum_{\ell=0}^L(N_{\ell}+1)W_{\ell}\quad \text{subject to} \begin{cases}\sum_{\ell=0}^{L}\frac{D_{\ell}}{N_{\ell}^{a}}\leq  TOL^2,\\[4pt]
  \frac{D_{\rm{iv}}}{M_L^{c}}\leq TOL^2,\\[4pt]
  D_{\rm{w}}h_L^{\eta_{\rm{w}}}\leq TOL,\\[4pt]
  \frac{D_{\rm{bi}}}{M_L^{c}}\leq TOL,\end{cases}
\end{align}
where we introduced the generic exponents $a$ and $c$ for $N_{\ell}$ and $M_{\ell}$, respectively, where $0\leq \ell\leq L$. Let
\begin{align}\label{eq:Dell}
  D_{\ell}\coloneqq \begin{cases}\tilde{D}, & \ell=0\\
  \frac{\tilde{D}_{\rm{iv}}^{(1)}}{M_{\ell}^{b}}+\tilde{D}_{\rm{iv}}^{(2)}h_{\ell}^{2\eta_{{\rm{s}},d_1}}, & 1\leq\ell\leq L
  \end{cases},
\end{align}
where $\tilde{D}$, $\tilde{D}_{\rm{iv}}^{(1)}$, and $\tilde{D}_{\rm{iv}}^{(2)}$ are independent of $L$ and $N_0,\ldots,N_{L}$.
Moreover, $D_{\rm{bi}}, D_{\rm{iv}}, D_{\rm{w}},$ and $D_{\ell}$ for $1\leq \ell,\leq L-1$ are independent of $L$ and $N_0,\ldots,N_{L}$; only $D_L$ depends on L.
 It follows from Propositions~\ref{prop:B.MLrQ-Q.h} to~\ref{prop:V.ell.MLQ.h.L} that
\begin{equation}
    a=2-2\epsilon-2A_{\rm{max}},
\end{equation}
and
\begin{equation}
    b=c=2-2\epsilon.
\end{equation}
We revisit this proof later under stricter regularity assumptions and hence make this distinction here.
It immediately follows from~\eqref{eq:Lagrangian.h.det.simplified} that the constraint
\begin{equation}\label{eq:constr.B}
    \frac{D_{\rm{bi}}}{M_L^{c}}\leq TOL
\end{equation}
is satisfied whenever the constraint
\begin{equation}\label{eq:constr.M.L}
    \frac{D_{\rm{iv}}}{M_L^{c}}\leq TOL^2
\end{equation}
is satisfied as $TOL\to 0$ and the former is thus ignored in the following derivation. We set
\begin{equation}\label{eq:M.ell.eta.w}
    M_{\ell}\coloneqq M_0\left(2^{\frac{2\eta_{\rm{w}}}{c}}\right)^{\ell},
\end{equation}
as described in Remark~\ref{rmk:optimal.M.h}. The Lagrangian
\begin{align}
    \cl{L}(L,N_0,\ldots,N_{L},\lambda,\mu,\nu)={}&\sum_{\ell=0}^{L}(N_{\ell}+1)W_{\ell}+\lambda \left(\sum_{\ell=0}^{L}\frac{D_{\ell}}{N_{\ell}^{a}}- TOL^2\right)+\mu\left(\frac{D_{\rm{iv}}}{M_L^{c}}- TOL^2\right)\nonumber\\
    {}&+\nu\left(D_{\rm{w}}h_L^{\eta_{\rm{w}}}-TOL\right),
\end{align}
is introduced. Taking the derivatives with respect to $\mu$ and $\nu$ and setting them to zero yields the optimal $L^{\ast}$ as follows:
\begin{equation}\label{eq:L.ast}
    L^{\ast}=\max\left\{\log_2\left(\left(\frac{h_0^{\eta_{\rm{w}}}C_{\rm{w}}}{TOL}\right)^{\frac{1}{\eta_{\rm{w}}}}\right),\log_2\left(\left(\frac{D_{\rm{iv}}^{\frac{1}{2}} }{M_0^{\frac{c}{2}}TOL}\right)^{\frac{1}{\eta_{\rm{w}}}}\right)\right\}.
\end{equation}
Next, taking the derivative with respect to $N_0$ yields
\begin{equation}
    \lambda^{\ast}=N_0^{a+1}\frac{W_0}{aD_0}.
\end{equation}
Taking the derivatives with respect to $N_{\ell}$, where $1\leq \ell\leq L^{\ast}$, yields:
\begin{align}\label{eq:N.L-1}
    N_\ell^{a+1}={}&\lambda^{\ast}\frac{aD_{\ell}}{W_{\ell}},\nonumber\\
    ={}&\frac{D_{\ell}W_0}{D_0W_{\ell}}N_0^{a+1}, \quad 1\leq \ell\leq L^{\ast}.
\end{align}
Taking the derivative with respect to $\lambda$ results in
\begin{align}
   TOL^2={}& \sum_{\ell=0}^{L^{\ast}}\frac{D_{\ell}}{N_{\ell}^{a}},\nonumber\\
   ={}&\sum_{\ell=0}^{L^{\ast}}\frac{D_{\ell}}{\left(\frac{D_{\ell}W_0}{D_0W_{\ell}}\right)^{\frac{a}{a+1}}N_{0}^{a}},\nonumber\\
   ={}&\left(\frac{D_0}{W_0}\right)^{\frac{a}{a+1}}\sum_{\ell=0}^{L^{\ast}}\frac{\left(D_{\ell}W_{\ell}^{a}\right)^{\frac{1}{a+1}}}{N_{0}^{a}}.
\end{align}
Rearranging yields
\begin{equation}
    N_0^{\ast a}= TOL^{-2}\left(\frac{D_0}{W_0}\right)^{\frac{a}{a+1}}\sum_{\ell=0}^{L^{\ast}}\left(D_{\ell}W_{\ell}^{a}\right)^{\frac{1}{a+1}},
\end{equation}
from which it follows that
\begin{equation}
    N_0^{\ast}= TOL^{-\frac{2}{a}}\left(\frac{D_0}{W_0}\right)^{\frac{1}{a+1}}\left(\sum_{\ell=0}^{L^{\ast}}\left(D_{\ell}W_{\ell}^{a}\right)^{\frac{1}{a+1}}\right)^{\frac{1}{a}}.
\end{equation}
Substituting this into \eqref{eq:N.L-1} yields
\begin{equation}\label{eq:N.ast}
    N_{\ell}^{\ast}=TOL^{-\frac{2}{a}}\left(\frac{D_{\ell}}{W_{\ell}}\right)^{\frac{1}{a+1}}\left(\sum_{k=0}^{L^{\ast}}\left(D_{k}W_{k}^{a}\right)^{\frac{1}{a+1}}\right)^{\frac{1}{a}}, \quad 0 \leq \ell\leq L^{\ast}.
\end{equation}
The first term in the total work is given as
\begin{align}\label{eq:opt.work.MLDLQ}
    \sum_{\ell=0}^{L^{\ast}}N_{\ell}^{\ast}W_{\ell}{}&=TOL^{-\frac{2}{a}}\sum_{\ell=0}^{L^{\ast}}W_{\ell}\left(\frac{D_{\ell}}{W_{\ell}}\right)^{\frac{1}{a+1}}\left(\sum_{k=0}^{L^{\ast}}\left(D_{k}W_{k}^{a}\right)^{\frac{1}{a+1}}\right)^{\frac{1}{a}},\nonumber\\
    {}&=TOL^{-\frac{2}{a}}\left(\sum_{\ell=0}^{L^{\ast}}\left(D_{\ell}W_{\ell}^{a}\right)^{\frac{1}{a+1}}\right)^{\frac{a+1}{a}}.
\end{align}
It follows from \eqref{eq:Dell} and \eqref{eq:L.ast} that
\begin{align}
    D_{L^{\ast}}{}&=\cl{O}\left( TOL^{2\min\left\{\frac{\eta_{{\rm{s}},d_1}}{\eta_{\rm{w}}},\frac{b}{c}\right\}}\right),
\end{align}
and from \eqref{eq:W.ell} that
\begin{equation}
    W_{L^{\ast}}^{a}=\cl{O} \left(TOL^{-a\left(\frac{2}{c}+\frac{\gamma}{\eta_{\rm{w}}}\right)}\right).
\end{equation}
Moreover, it follows that
\begin{equation}\label{eq:D.L.W.L}
    D_{L^{\ast}}W_{L^{\ast}}^{a}= \cl{O}\left(TOL^{2\min\left\{\frac{\eta_{{\rm{s}},d_1}}{\eta_{\rm{w}}},\frac{b}{c}\right\}-a\left(\frac{2}{c}+\frac{\gamma}{\eta_{\rm{w}}}\right)}\right).
\end{equation}
If the exponent in \eqref{eq:D.L.W.L} is positive, and thus the sum in \eqref{eq:opt.work.MLDLQ} converges, it follows that
\begin{align}
    \sum_{\ell=0}^{L^{\ast}}N_{\ell}^{\ast}W_{\ell}{}&=\cl{O}\left(TOL^{-\frac{2}{a}}\right),\nonumber\\
    {}&=\cl{O}\left(TOL^{-\frac{1}{1-\epsilon-A_{\rm{max}}}}\right).
\end{align}
If the exponent in \eqref{eq:D.L.W.L} is equal to zero, then the sum in \eqref{eq:opt.work.MLDLQ} is proportional to $L^{\ast}$, and it follows that
\begin{align}
    \sum_{\ell=0}^{L^{\ast}}N_{\ell}^{\ast}W_{\ell}{}&=\cl{O}\left(TOL^{-\frac{2}{a}}\log_2\left(TOL^{-1}\right)\right),\nonumber\\
    {}&=\cl{O}\left(TOL^{-\frac{1}{1-\epsilon-A_{\rm{max}}}}\log_2\left(TOL^{-1}\right)\right).
\end{align}
If the exponent in \eqref{eq:D.L.W.L} is negative, and thus the sum in \eqref{eq:opt.work.MLDLQ} diverges, it can be bounded by the largest term, that is, $(D_{L}W_L^a)^{1/a}$, and thus,
\begin{align}
    \sum_{\ell=0}^{L^{\ast}}N_{\ell}^{\ast}W_{\ell}{}&=\cl{O}\left(TOL^{-\frac{2}{a}\left(1-\min\left\{\frac{\eta_{{\rm{s}},d_1}}{\eta_{\rm{w}}},\frac{b}{c}\right\}\right)-\left(\frac{2}{c}+\frac{\gamma}{\eta_{\rm{w}}}\right)}\right).
\end{align}
It follows from $\eta_{{\rm{s}},d_1}\leq \eta_{\rm{w}}$ and $a=b$ that
\begin{align}
    \sum_{\ell=0}^{L^{\ast}}N_{\ell}^{\ast}W_{\ell}{}&=\cl{O}\left(TOL^{-\frac{2}{a}\left(1-\frac{\eta_{{\rm{s}},d_1}}{\eta_{\rm{w}}}\right)-\left(\frac{2}{c}+\frac{\gamma}{\eta_{\rm{w}}}\right)}\right),\nonumber\\
    {}&=\cl{O}\left(TOL^{-\frac{1}{1-\epsilon-A_{\rm{max}}}+\frac{\eta_{{\rm{s}},d_1}}{\eta_{\rm{w}}(1-\epsilon-A_{\rm{max}})}-\left(\frac{1}{1-\epsilon}+\frac{\gamma}{\eta_{\rm{w}}}\right)}\right).
\end{align}
 Finally, the second term in the total work is bounded as follows:
 \begin{align}
    \sum_{\ell=0}^{L^{\ast}}W_{\ell}{}&=\cl{O}\left(TOL^{-\left(\frac{2}{c}+\frac{\gamma}{\eta_{\rm{w}}}\right)}\right),\nonumber\\
    {}&=\cl{O}\left(TOL^{-\left(\frac{1}{1-\epsilon}+\frac{\gamma}{\eta_{\rm{w}}}\right)}\right).
\end{align}
\end{proof}
Instead of choosing $h_{\ell}\propto 2^{-\ell}$ and $M_{\ell}\propto (2^{2\eta_{\rm{w}}/c})^{\ell}$ in the above proof, we may also choose $M_{\ell}\propto 2^{\ell}$ and $h_{\ell}\propto (2^{c/2\eta_{\rm{w}}})^{-\ell}$ for $1\leq\ell\leq L$, see Remark~\ref{rmk:optimal.M.h}. This leads to $L^{\ast}\propto\log_2(TOL^{-2/c})$ in~\eqref{eq:L.ast}. However, $D_{L^{\ast}}$ and $W_{L^{\ast}}$ in~\eqref{eq:D.L.W.L} still exhibit similar asymptotic behavior.
\begin{rmk}[Choice of splitting parameters in the optimal work derivation]
    Several error splitting parameters equal to $1/2$ are introduced in the derivation of the optimal computational work of the MLDLQMC estimator in the proof of Proposition~\ref{prop:Work.MLQ.h.det}. Deriving optimal expressions for these error splitting parameters in MLMC is outside the scope of this study; and we refer to the work~\cite{Haj14} for a thorough investigation.
\end{rmk}
\begin{rmk}[Impact of the dimension on the estimator work]\label{rmk:dimension.dependence}
    The integration dimensions of the inner and outer integrands may affect the multiplicative error terms implied by the $\cl{O}$ notation in the above results. Investigating this dependence is not the focus of the present study; however, \cite{Bar23} successfully applied a nested rQMC estimator for EIG estimation up to a dimension of 15, and \cite{Kaa24} estimated the EIG of an experiment at an even higher dimension using rQMC methods and proved independence of the inner integral estimation of the dimension.
\end{rmk}
\begin{cor}[Optimal work of the MLDLQMC estimator with exact sampling]\label{prop:Work.MLQ.det}
Under the same assumptions as Proposition~\ref{prop:Work.MLQ.h.det}, the total computational work of the MLDLQMC estimator \eqref{eq:MLrQ.h.det} with exact sampling is given as follows:
\begin{align}
  W_{\rm{MLDLQ}}^{\ast}={}&\cl{O}\left(TOL^{-\frac{1}{1-\epsilon-A_{\rm{max}}}}\right)
\end{align}
for any $\epsilon>0$ as $TOL\to 0$.
\end{cor}
\begin{proof}
    It holds that
   \begin{equation}
       \frac{1}{1-\epsilon-A_{\rm{max}}}>\frac{1}{1-\epsilon}
   \end{equation} 
for any $\epsilon>0$, where $0<A_{\rm{max}}<1/2$. The result then follows immediately from Proposition~\ref{prop:Work.MLQ.h.det} in the limit as $h_{\ell}\to0$, where $0\leq \ell\leq L$, as $\gamma=0$ and $\eta_{\rm{w}}=\eta_{{\rm{s}},d_1}=1$ without loss of generality.
\end{proof}
The above results in Proposition~\ref{prop:Work.MLQ.h.det} and Corollary~\ref{prop:Work.MLQ.det} apply to rQMC estimator constructions based on several different low-discrepancy points and randomization methods. For the particular choice of the Sobol' sequence and Owen's scrambling, the stronger result in Proposition~\ref{prop:GHK} may be applied, given that the integrands satisfy stricter regularity criteria imposed by the generalized Hardy--Krause variation of order one.
\begin{cor}[Optimal work of the MLDLQMC estimator for scrambled Sobol' sequence]\label{cor:opt.case}
Let the assumptions of Propositions~\ref{prop:B.MLrQ-Q.h}--\ref{prop:V.ell.MLQ.h.L} hold with $V_{\rm{HK}}$ replaced by $V_{\rm{GHK}}$. Moreover, let $f\circ \bar{g}$ have bounded generalized Hardy--Krause variation $V_{\rm{GHK}}$ of order one (see~\cite[Definition 6.24]{Dic13}). Then, the total computational work of the MLDLQMC estimator~\eqref{eq:MLrQ.h.det} with inexact sampling is given by
\begin{align}
  W_{\rm{MLDLQ}}^{\ast}=
  \cl{O}\left(TOL^{-\frac{2}{3-2\epsilon} \left(1-\min\left\{\frac{\eta_{{\rm{s}},d_1}}{\eta_{\rm{w}}},\frac{2-2\epsilon}{3-2\epsilon}\right\}\right)-\left(\frac{2}{3-2\epsilon}+\frac{\gamma}{\eta_{\rm{w}}}\right)}\right),
\end{align}
 where the multiplicative constant in the $\cl{O}$ notation approaches infinity as $\epsilon\to 0$. Moreover, the total computational work of the MLDLQMC estimator~\eqref{eq:MLrQ.h.det} with exact sampling is given by
\begin{align}
  W_{\rm{MLDLQ}}^{\ast}=
  \cl{O}\left(TOL^{-\frac{8-4\epsilon}{(3-2\epsilon)^2}}\right).
\end{align}
\end{cor}
\begin{proof}
    We replace the exponents $a$, $b$, and $c$ in~\eqref{eq:Lagrangian.h.det.simplified} and~\eqref{eq:Dell} in the proof of Proposition~\ref{prop:Work.MLQ.h.det} by
    \begin{equation}
        a=c=3-2\epsilon,
    \end{equation}
    and
        \begin{equation}
        b=2-2\epsilon.
    \end{equation}
  This follows from the proof of~\cite[Corollary~1]{Bar23}. The exponent in~\eqref{eq:D.L.W.L} then becomes
    \begin{align}
        2\min\left\{\frac{\eta_{{\rm{s}},d_1}}{\eta_{\rm{w}}},\frac{b}{c}\right\}-a\left(\frac{2}{c}+\frac{\gamma}{\eta_{\rm{w}}}\right){}&=2\min\left\{\frac{\eta_{{\rm{s}},d_1}}{\eta_{\rm{w}}},\frac{2-2\epsilon}{3-2\epsilon}\right\}-2-(3-2\epsilon)\frac{\gamma}{\eta_{\rm{w}}},\nonumber\\
        {}&<0.
    \end{align}
    From this, it follows that the sum in~\eqref{eq:opt.work.MLDLQ} diverges, and thus,
\begin{align}\label{eq:Nl.Wl.GHK}
    \sum_{\ell=0}^{L^{\ast}}N_{\ell}^{\ast}W_{\ell}{}&=\cl{O}\left(TOL^{-\frac{2}{a}\left(1-\min\left\{\frac{\eta_{{\rm{s}},d_1}}{\eta_{\rm{w}}},\frac{b}{c}\right\}\right)-\left(\frac{2}{c}+\frac{\gamma}{\eta_{\rm{w}}}\right)}\right),\nonumber\\
    {}&=\cl{O}\left(TOL^{-\frac{2}{3-2\epsilon} \left(1-\min\left\{\frac{\eta_{{\rm{s}},d_1}}{\eta_{\rm{w}}},\frac{2-2\epsilon}{3-2\epsilon}\right\}\right)-\left(\frac{2}{3-2\epsilon}+\frac{\gamma}{\eta_{\rm{w}}}\right)}\right).
\end{align}
and finally,
 \begin{align}
    \sum_{\ell=0}^{L^{\ast}}W_{\ell}{}&=\cl{O}\left(TOL^{-\left(\frac{2}{c}+\frac{\gamma}{\eta_{\rm{w}}}\right)}\right),\nonumber\\
    {}&=\cl{O}\left(TOL^{-\left(\frac{
    2}{3-2\epsilon}+\frac{\gamma}{\eta_{\rm{w}}}\right)}\right),
\end{align}
which is asymptotically negligible compared to the term in~\eqref{eq:Nl.Wl.GHK}. The result for exact sampling follows immediately.
\end{proof}
From the result in Corollary~\ref{cor:opt.case}, it follows that the computational cost of the MLDLQMC estimator for exact sampling increases at a rate of almost $\cl{O}(TOL^{-8/9})$ in the optimal case where both the inner and the outer integrand are highly regular. We demonstrate this rate for a simple polynomial example below. For the EIG application, we do not expect this high degree of regularity. Applying a fixed truncation to the Gaussian observation noise (see~\cite{Bar23}) ensures such regularity; however, as the region of truncation extends towards the boundary of the integration domain, this effect is negated. Only the inner integrand exhibits bounded generalized Hardy--Krause variation of order one for this application. We derive the total computational work for this setting below.
\begin{cor}[Optimal work of the MLDLQMC estimator for scrambled Sobol' sequence for EIG application]\label{cor:EIG.case}
Let the assumptions of Propositions~\ref{prop:B.MLrQ-Q.h}--\ref{prop:V.ell.MLQ.h.L} hold with $V_{\rm{HK}}$ replaced by $V_{\rm{GHK}}$. Then, the total computational work of the MLDLQMC estimator~\eqref{eq:MLrQ.h.det} with inexact sampling is given by
\begin{align}
  W_{\rm{MLDLQ}}^{\ast}={}&\begin{cases}\cl{O}\left(TOL^{-\frac{1}{1-\epsilon-A_{\rm{max}}}}\right), & \text{if}\quad \frac{\min\left\{\frac{\eta_{{\rm{s}},d_1}}{\eta_{\rm{w}}},\frac{2-2\epsilon}{3-2\epsilon}\right\}}{1-\epsilon-A_{\rm{max}}}>\left(\frac{2}{3-2\epsilon}+\frac{\gamma}{\eta_{\rm{w}}}\right),\\
  \cl{O}\left(TOL^{-\frac{1}{1-\epsilon-A_{\rm{max}}}} \log_{2}\left(TOL^{-1} \right)\right), & \text{if}\quad \frac{\min\left\{\frac{\eta_{{\rm{s}},d_1}}{\eta_{\rm{w}}},\frac{2-2\epsilon}{3-2\epsilon}\right\}}{1-\epsilon-A_{\rm{max}}}=\left(\frac{2}{3-2\epsilon}+\frac{\gamma}{\eta_{\rm{w}}}\right),\\
  \cl{O}\left(TOL^{-\frac{1}{1-\epsilon-A_{\rm{max}}}\left(1-\min\left\{\frac{\eta_{{\rm{s}},d_1}}{\eta_{\rm{w}}},\frac{2-2\epsilon}{3-2\epsilon}\right\}\right)-\left(\frac{2}{3-2\epsilon}+\frac{\gamma}{\eta_{\rm{w}}}\right)}\right), & \text{if}\quad \frac{\min\left\{\frac{\eta_{{\rm{s}},d_1}}{\eta_{\rm{w}}},\frac{2-2\epsilon}{3-2\epsilon}\right\}}{1-\epsilon-A_{\rm{max}}}<\left(\frac{2}{3-2\epsilon}+\frac{\gamma}{\eta_{\rm{w}}}\right),\end{cases}
\end{align}
 where the multiplicative constant in the $\cl{O}$ notation approaches infinity as $\epsilon\to 0$. Moreover, the total computational work of the MLDLQMC estimator~\eqref{eq:MLrQ.h.det} with exact sampling is given by
\begin{align}
  W_{\rm{MLDLQ}}^{\ast}=
  \cl{O}\left(TOL^{-\frac{1}{1-\epsilon-A_{\rm{max}}}}\right).
\end{align}
\end{cor}
\begin{proof}
    In this setting, we replace the exponents $a$, $b$, and $c$ in~\eqref{eq:Lagrangian.h.det.simplified} and~\eqref{eq:Dell} in the proof of Proposition~\ref{prop:Work.MLQ.h.det} by
    \begin{equation}
        a=2-2\epsilon-2A_{\rm{max}},
    \end{equation}
        \begin{equation}
        b=2-2\epsilon,
    \end{equation}
    and
        \begin{equation}
        c=3-2\epsilon.
    \end{equation}
 as a consequence of the proof of~\cite[Corollary~1]{Bar23}. The exponent in~\eqref{eq:D.L.W.L} then becomes
    \begin{align}
        2\min\left\{\frac{\eta_{{\rm{s}},d_1}}{\eta_{\rm{w}}},\frac{b}{c}\right\}-a\left(\frac{2}{c}+\frac{\gamma}{\eta_{\rm{w}}}\right){}&=2\min\left\{\frac{\eta_{{\rm{s}},d_1}}{\eta_{\rm{w}}},\frac{2-2\epsilon}{3-2\epsilon}\right\}-(2-2\epsilon-2A_{\rm{max}})\left(\frac{2}{3-2\epsilon}+\frac{\gamma}{\eta_{\rm{w}}}\right).
    \end{align}
    If this exponent is positive, it follows that
\begin{align}\label{eq:Nl.Wl.EIG.conv}
    \sum_{\ell=0}^{L^{\ast}}N_{\ell}^{\ast}W_{\ell}{}&=\cl{O}\left(TOL^{-\frac{1}{1-\epsilon-A_{\rm{max}}}}\right).
\end{align}
If it is zero, then,
\begin{align}\label{eq:Nl.Wl.EIG.const}
    \sum_{\ell=0}^{L^{\ast}}N_{\ell}^{\ast}W_{\ell}{}&=\cl{O}\left(TOL^{-\frac{1}{1-\epsilon-A_{\rm{max}}}}\log_2\left(TOL^{-1}\right)\right).
\end{align}
And finally, if it is negative, then
\begin{align}\label{eq:Nl.Wl.EIG.div}
    \sum_{\ell=0}^{L^{\ast}}N_{\ell}^{\ast}W_{\ell}{}&=\cl{O}\left(TOL^{-\frac{2}{a}\left(1-\min\left\{\frac{\eta_{{\rm{s}},d_1}}{\eta_{\rm{w}}},\frac{b}{c}\right\}\right)-\left(\frac{2}{c}+\frac{\gamma}{\eta_{\rm{w}}}\right)}\right),\nonumber\\
    {}&=\cl{O}\left(TOL^{-\frac{1}{1-\epsilon-A_{\rm{max}}}\left(1-\min\left\{\frac{\eta_{{\rm{s}},d_1}}{\eta_{\rm{w}}},\frac{2-2\epsilon}{3-2\epsilon}\right\}\right)-\left(\frac{2}{3-2\epsilon}+\frac{\gamma}{\eta_{\rm{w}}}\right)}\right).
\end{align}
In the exact sampling case, the exponent is always positive because of the $2A_{\rm{max}}>0$ term.
\end{proof}
Comparing the results in Corollary~\ref{cor:opt.case} and~\ref{cor:EIG.case} for exact sampling, we see that the improved convergence rate in the inner integral approximation does not affect the overall computational work asymptotically. However, there may still be a reduction in the multiplicative terms implied in the $\cl{O}$ notation.
\subsection{Example 1: Polynomial example with exact sampling}\label{sec:polynomial.example}
For the practical implementation of the MLDLQMC estimator applied to the EIG, the CLT confidence parameter $C_{\alpha}=\Phi^{-1}(1-\alpha/2)$ was selected for $\alpha=0.05$ (i.e., $C_{\alpha}=1.96$). This choice proved sufficient, and the more conservative Chebyshev's inequality was not needed (see Figure~\ref{fig:polynomial.evt.work} (C), Figure~\ref{fig:linear.evt.work} (C), Figure~\ref{fig:poisson.evt.work} (C), and Figure~\ref{fig:duffing.evt.work} (C)). Pilot runs were employed to estimate the required constants and convergence rates using several randomizations $S,R$, independent of the error tolerance $TOL$. Thereafter, the optimal number of levels $L^{\ast}$ and outer samples $N_{\ell}^{\ast}$, where $0\leq\ell\leq L^{\ast}$, were derived for a specified error tolerance $TOL>0$ according to the expressions~\eqref{eq:L.ast} and~\eqref{eq:N.ast}.
\subsubsection{Setting and notation}
We apply the MLDLQMC estimator to the polynomial example from~\cite{Bar23} before we introduce the EIG formulation in the following section. Here, both the inner integrand $g$ and the outer integrand $f$ are polynomials over $[0,1]^{d_1}$ and $[0,1]^{d_2}$, respectively, where $d_1=d_2=30$. The outer integrand is given as follows:
\begin{equation}\label{eq:pol.outer}
    f(z)\coloneqq z^2+z, \quad \text{where } z\in\mathbb{R},
\end{equation}
and the inner integrand is given as
\begin{equation}\label{eq:pol.inner}
    g(\bs{y},\bs{x})\coloneqq \sum_{i=1}^{30}(y_i + x_i +  y_ix_i), \quad \text{where } \bs{y},\bs{x}\in[0,1]^{30}.
\end{equation}
It follows immediately that the assumptions of Corollary~\ref{cor:opt.case} hold for the exact sampling case. The closed-form solution to this nested integration problem is presented in~\cite[Subsection 3.6.]{Bar23}.

\subsubsection{Pilot runs}
\begin{figure}[ht]
	\subfloat[The observed convergence rate of the variance of the inner samples at the final level $L$ of 3.]{%
		\includegraphics[width=0.45\textwidth]{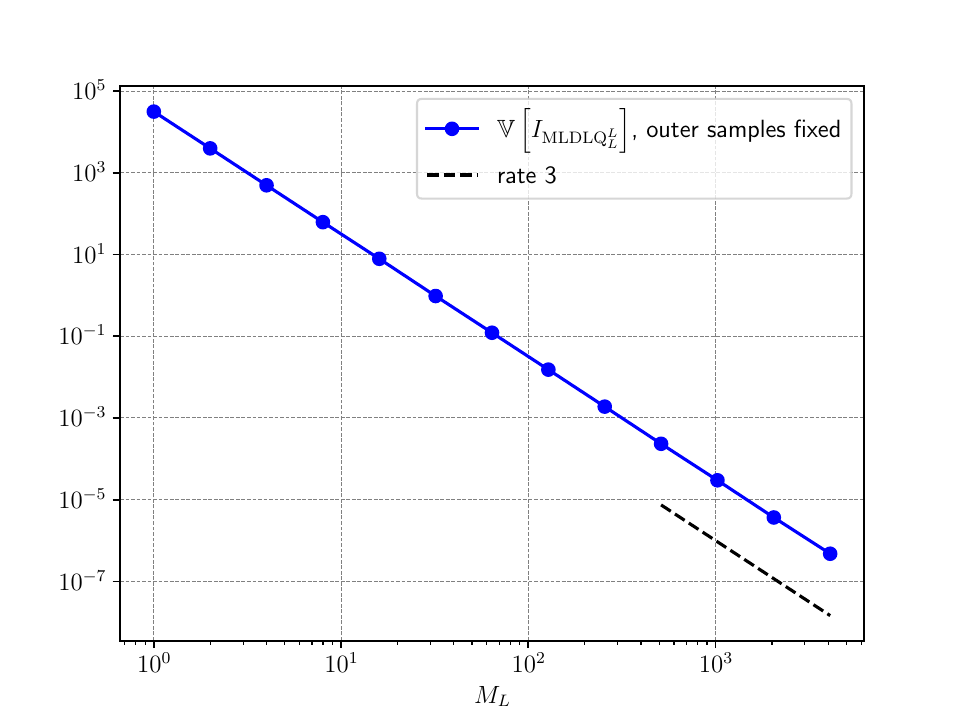}
	}
	\subfloat[The observed convergence rate of the variance at a rate of 3 with a reduced multiplicative factor for the difference estimator.]{%
		\includegraphics[width=0.45\textwidth]{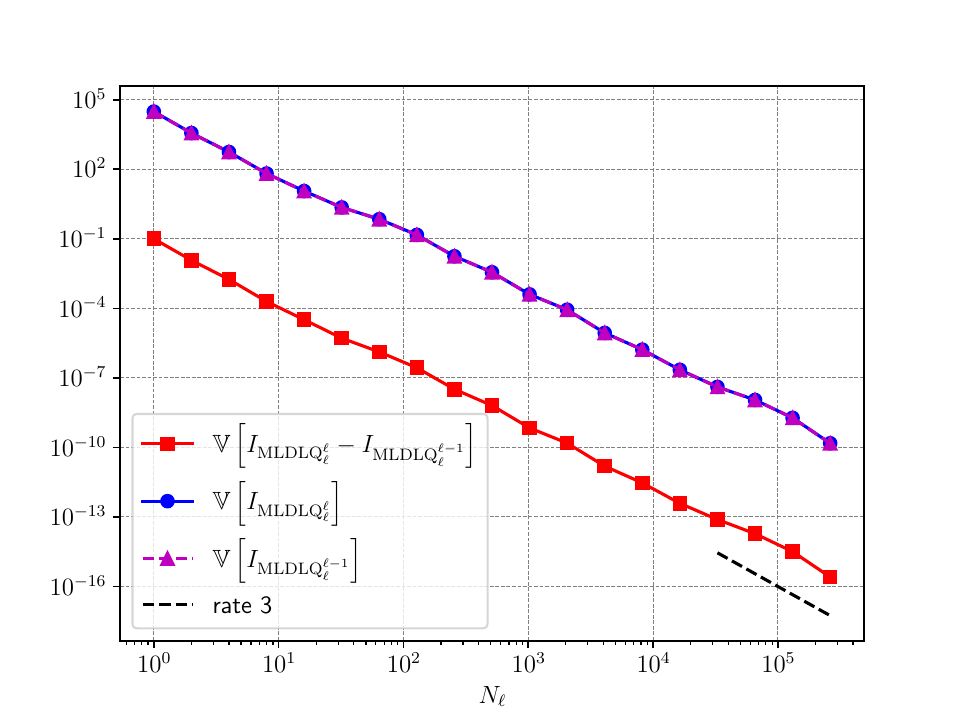}
	}\\
 \subfloat[The variance of the level difference converged at a rate of 2 as a function of $M_{\ell}$ for $N_{\ell}=1$ fixed.]{%
		\includegraphics[width=0.45\textwidth]{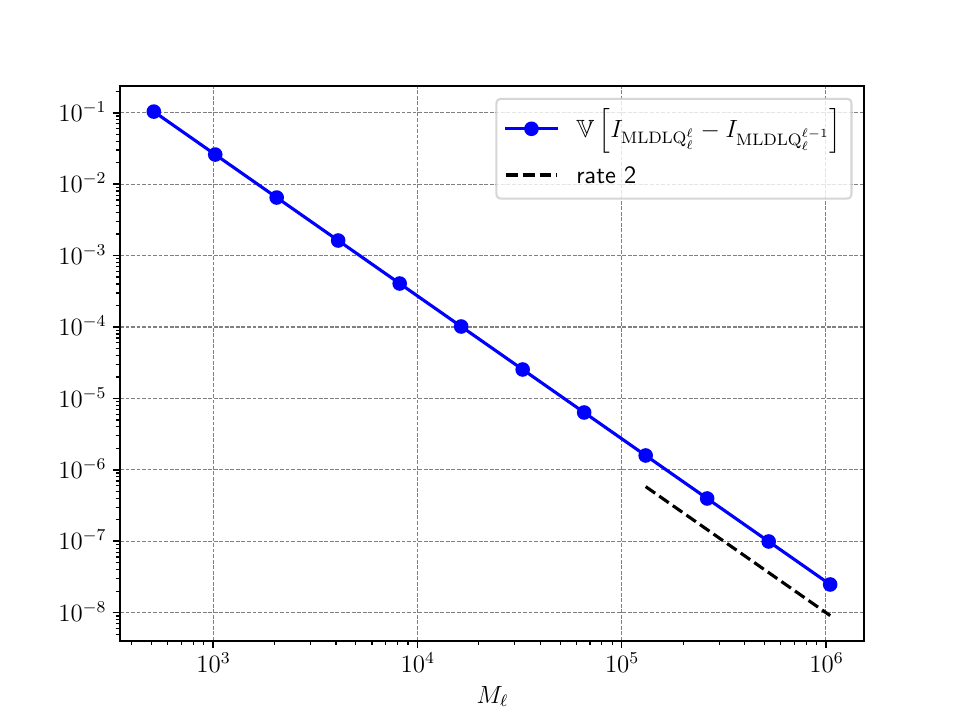}
	}
	\caption{Example 1: Pilot runs for the MLDLQMC estimator with exact sampling. Each observed variance contribution converged at a rate consistent with the theoretical analysis.}
	\label{fig:polynomial.pilot}
\end{figure}

We used pilot runs to verify that the convergence rates indicated by our theoretical results can be achieved in practice and to estimate the multiplicative constants appearing in our error bounds. First, we used $N_L=4$ low-discrepancy points and $S=1$, $R=2^{15}$ randomizations to see that the variance contribution from the randomized inner samples at level $L$ converged at rate 3 (Figure~\ref{fig:polynomial.pilot}~(A)). Next, we used $M_{\ell-1}=2^8$, $M_{\ell}=2^9$ samples and $S=2^7$ randomizations to demonstrate that the variance of the level differences converged at rate 3 as a function of $N_{\ell}$ (Figure~\ref{fig:polynomial.pilot}~(B)). The same rate was also observed for the single-level estimators with a larger multiplicative constant, demonstrating the effect of the multilevel method. Lastly, we used $S=2^{10}$ randomizations to show that the variance of the level difference as a function of $M_{\ell}$ converged at a rate of 2 for $N_{\ell}=1$ fixed (Figure~\ref{fig:polynomial.pilot}~(C)). The lower rate compared to the one observed in Panel (A) is because the inner samples of the intermediate levels are not randomized. Introducing randomness to those samples would complicate the estimator structure and lead to additional variance terms.

\subsubsection{Consistency and work analysis}
\begin{figure}[ht]
\subfloat[The optimal number of levels $L^{\ast}$ increased at a rate of $\log_2(TOL^{-2/3})$ as a function of the error tolerance $TOL>0$.]{%
		\includegraphics[width=0.45\textwidth]{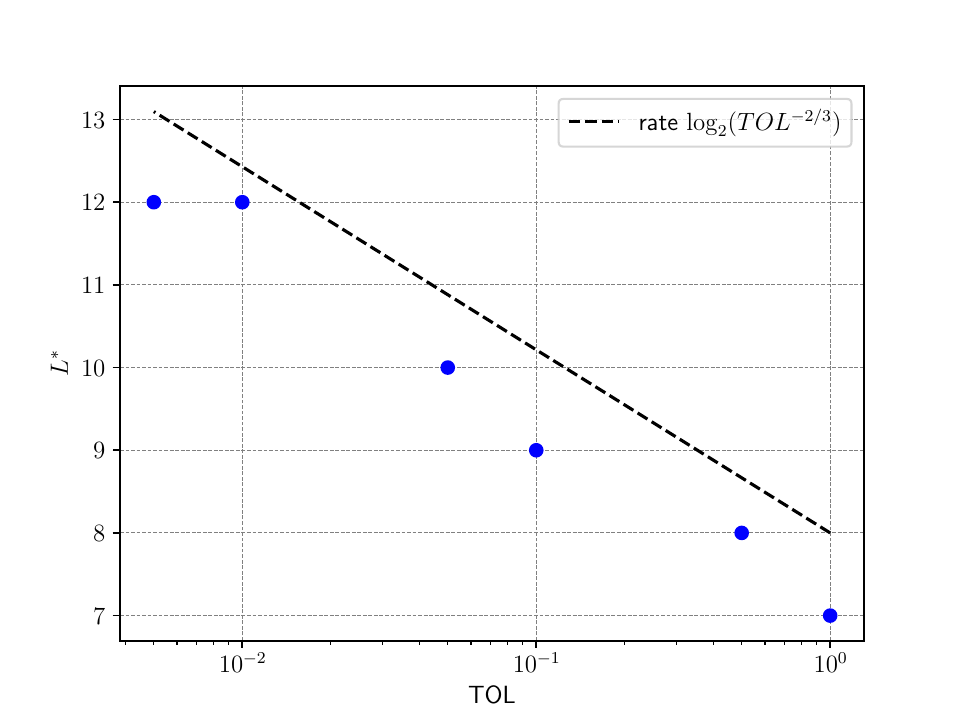}
	}
 \subfloat[Total computational work of the MLDLQMC estimator as a function of the error tolerance $TOL$.]{%
		\includegraphics[width=0.45\textwidth]{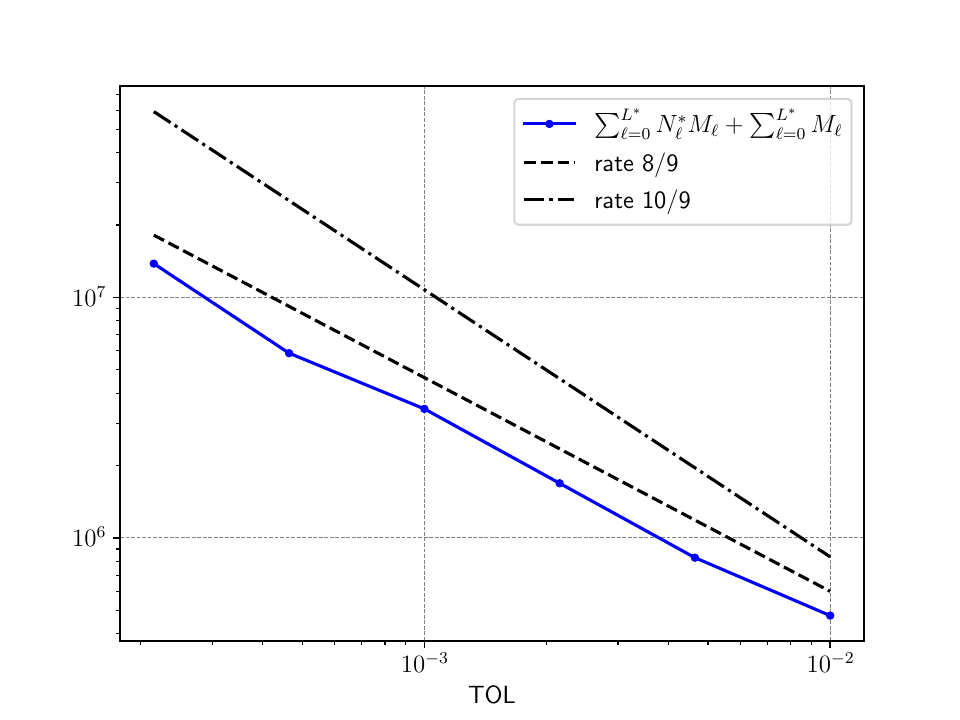}
	}\\
	\subfloat[Error vs.~tolerance $(TOL)$ for 100 runs of the MLDLQMC estimator. Each run resulted in an error less than the tolerance projected by the central limit theorem.]{%
		\includegraphics[width=0.45\textwidth]{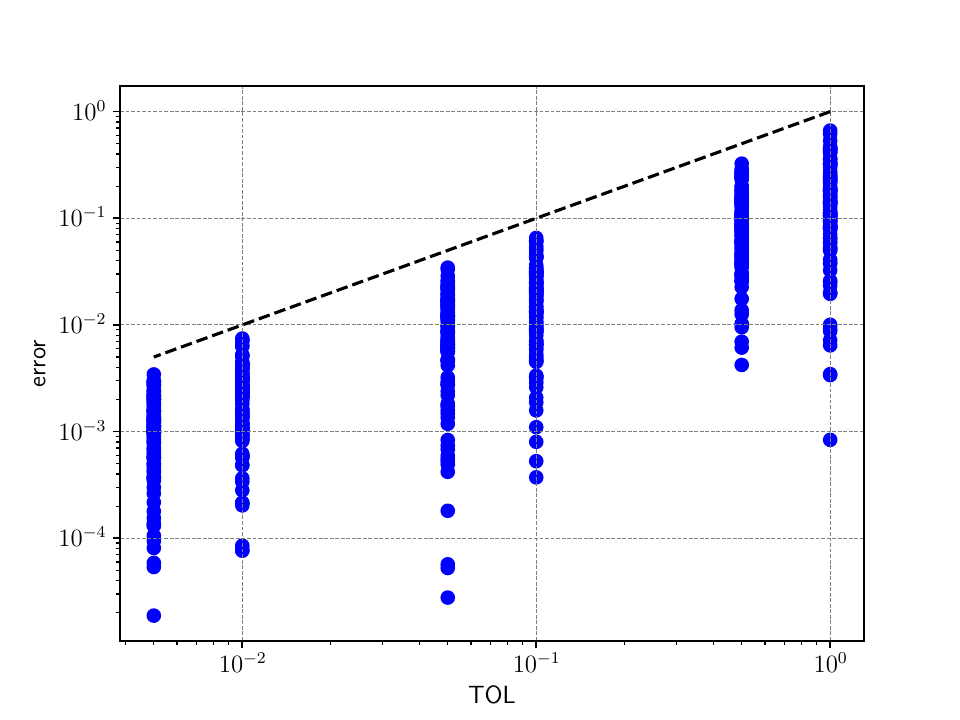}
	}
	
	\caption{Example 1: The MLDLQMC estimator with exact sampling resulted in error less than the projected tolerance $TOL$ (C) at a cost increasing at the projected rate of $8/9$ asymptotically as $TOL\to 0$ (B). This is an improvement compared to the rate of $10/9$ achieved by the single-level estimator in~\cite{Bar23}.}
	\label{fig:polynomial.evt.work}
\end{figure}

Given the estimates from the pilot runs, we determined the optimal final level $L^{\ast}$ as a function of the tolerance $TOL>0$ (Figure~\ref{fig:polynomial.evt.work}~(A)). We observed an increase at a rate of $\log_2(TOL^{-2/3})$. The optimal work of the MLDLQMC estimator with exact sampling for the highly regular polynomial integrands in~\eqref{eq:pol.outer} and~\eqref{eq:pol.inner} as a function of $N_{\ell}^{\ast}$, where $1\leq \ell\leq L^{\ast}$, increased at a rate of $TOL^{-8/9}$, as demonstrated in Figure~\ref{fig:polynomial.evt.work}~(B). This is an improvement of a factor of $2/9$ compared to the single-level estimator from~\cite{Bar23} applied to the same example, where a total computational work of $TOL^{-10/9}$ was observed. Finally, we compared the results of 100 runs of the MLDLQMC estimator with $S=R=1$ randomizations at tolerances $TOL=\{10^0,5\cdot10^{-1},10^{-1},5\cdot10^{-2},10^{-2},5\cdot10^{-3}\}$ against the exact solution of this nested integration problem and found that each run remained below the prescribed tolerance (Figure~\ref{fig:polynomial.evt.work}~(C)). This indicates that the estimates of the multiplicative constants from the pilot runs was slightly conservative, as we would only expect approximately 95 out of 100 runs to remain below the tolerance for the specified confidence level. The value of the nested integration problem given by~\eqref{eq:pol.outer} and~\eqref{eq:pol.inner} in dimension 30 is equal to $I=1449.375$. The smallest relative tolerance depicted in Panel (C) is therefore of order $\cl{O}(10^{-6})$.

\section{Expected information gain estimation}\label{sec:EIG.estimation}
We now demonstrate how to apply the MLDLQMC estimator to estimate of the EIG of an experiment. This is an application from Bayesian experimental design, where uncertain parameters of interest $\bs{\theta}$ are modeled as random variables with a prior distribution $\pi(\bs{\theta})$. The selection of the prior reflects the a priori knowledge of the experimenter. The observation data is modeled as follows:
\begin{equation}\label{eq:data.model}
    \bs{Y}(\bs{\xi})=\bs{G}(\bs{\theta}_t,\bs{\xi})+\bs{\varepsilon},
\end{equation}
where $\bs{Y}\in\bb{R}^{d_{y}}$. The design $\bs{\xi}$ is omitted for conciseness, as it is not essential for the integral approximation that is the focus of the current work, and finding optimal designs is outside the scope of this study. Furthermore, $\bs{G}:\bb{R}^{d_{\theta}}\times\bb{R}^{d_{\xi}}\to \bb{R}^{d_y}$ is an appropriate model for the experiment evaluated at the unknown true parameter vector $\bs{\theta}_t$ and $\bs{\varepsilon}\stackrel{iid}{\sim}\cl{N}(\bs{0},\bs{\Sigma}_{\bs{\varepsilon}})$ is additive Gaussian noise with known covariance matrix $\bs{\Sigma}_{\bs{\varepsilon}}\in\bb{R}^{d_{y}\times d_{y}}$ (i.e., symmetric and strictly positive definite). The noise $\bs{\varepsilon}$ should not be confused with the factor $\epsilon$ in the convergence rate of rQMC methods used throughout the remainder of this study.
The EIG quantifies the amount of information that is expected to be gained from the data observation regarding the parameters of interest. It is defined as
\begin{equation}\label{eq:EIG.posterior}
    EIG\coloneqq \int_{\cl{Y}}\int_{\bs{\Theta}}\log\left(\frac{\pi(\bs{\theta}|\bs{Y})}{\pi(\bs{\theta})}\right)\pi(\bs{\theta}|\bs{Y})\di{}\bs{\theta}p(\bs{Y})\di{}\bs{Y},
\end{equation}
where $\pi(\bs{\theta}|\bs{Y})$ is the posterior distribution of the parameters of interest given the data, and $p(\bs{Y})$ is the evidence, or marginal likelihood. The data model \eqref{eq:data.model} implies the following likelihood:
\begin{equation}\label{eq:likelihood}
    p(\bs{Y}|\bs{\theta})=\frac{1}{\det(2\pi\bs{\Sigma}_{\bs{\varepsilon}})^{\frac{1}{2}}}e^{-\frac{1}{2}\left\lVert\bs{Y}-\bs{G}(\bs{\theta})\right\rVert_{\bs{\Sigma}_{\bs{\varepsilon}}^{-1}}^2},
\end{equation}
where $\lVert\cdot \rVert_{\bs{\Sigma}_{\bs{\varepsilon}}^{-1}}^2\coloneqq\left<\cdot,\cdot\right>_{\bs{\Sigma}_{\bs{\varepsilon}}^{-1}}$, and $\left<\cdot,\cdot\right>_{\bs{\Sigma}_{\bs{\varepsilon}}^{-1}}$ is the ${\bs{\Sigma}_{\bs{\varepsilon}}^{-1}}$-inner product, that is, $\left<\cdot,\cdot\right>_{\bs{\Sigma}_{\bs{\varepsilon}}^{-1}}\coloneqq \left<\cdot,{\bs{\Sigma}_{\bs{\varepsilon}}^{-1}}\cdot\right>$, and $\left<\cdot,\cdot\right>$ is the standard $L^2$ inner product. This likelihood is used together with Bayes' theorem and marginalization to rewrite \eqref{eq:EIG.posterior} as follows:
\begin{align}\label{eq:EIG}
    EIG{}&=\int_{\bs{\Theta}}\int_{\cl{Y}}\left[\log\left(p(\bs{Y}|\bs{\theta})\right)-\log\left(\int_{\bs{\Theta}}p(\bs{Y}|\bs{\vartheta})\pi(\bs{\vartheta})\di{}\bs{\vartheta}\right)\right]p(\bs{Y}|\bs{\theta})\di{}\bs{Y}\pi(\bs{\theta})\di{}\bs{\theta},
\end{align}
where $\bs{\vartheta}$ indicates an auxiliary variable for integration. Substituting \eqref{eq:data.model} into \eqref{eq:likelihood} yields
\begin{align}
    p(\bs{Y}|\bs{\theta}){}&=\frac{1}{\det(2\pi\bs{\Sigma}_{\bs{\varepsilon}})^{\frac{1}{2}}}e^{-\frac{1}{2}\left\lVert\bs{Y}-\bs{G}(\bs{\theta})\right\rVert_{\bs{\Sigma}_{\bs{\varepsilon}}^{-1}}^2},\nonumber\\
    {}&=\frac{1}{\det(2\pi\bs{\Sigma}_{\bs{\varepsilon}})^{\frac{1}{2}}}e^{-\frac{1}{2}\left\lVert\bs{G}(\bs{\theta})+\bs{\varepsilon}-\bs{G}(\bs{\theta})\right\rVert_{\bs{\Sigma}_{\bs{\varepsilon}}^{-1}}^2},\nonumber\\
    {}&=\frac{1}{\det(2\pi\bs{\Sigma}_{\bs{\varepsilon}})^{\frac{1}{2}}}e^{-\frac{1}{2}\left\lVert\bs{\varepsilon}\right\rVert_{\bs{\Sigma}_{\bs{\varepsilon}}^{-1}}^2},
\end{align}
and
\begin{align}
    p(\bs{Y}|\bs{\vartheta}){}&=\frac{1}{\det(2\pi\bs{\Sigma}_{\bs{\varepsilon}})^{\frac{1}{2}}}e^{-\frac{1}{2}\left\lVert\bs{Y}-\bs{G}(\bs{\vartheta})\right\rVert_{\bs{\Sigma}_{\bs{\varepsilon}}^{-1}}^2},\nonumber\\
    {}&=\frac{1}{\det(2\pi\bs{\Sigma}_{\bs{\varepsilon}})^{\frac{1}{2}}}e^{-\frac{1}{2}\left(\left\lVert\bs{G}(\bs{\theta})\right\rVert_{\bs{\Sigma}_{\bs{\varepsilon}}^{-1}}^2+\left\lVert\bs{\varepsilon}\right\rVert_{\bs{\Sigma}_{\bs{\varepsilon}}^{-1}}^2+\left\lVert\bs{G}(\bs{\vartheta})\right\rVert_{\bs{\Sigma}_{\bs{\varepsilon}}^{-1}}^2+2\left<\bs{G}(\bs{\theta}),\bs{\varepsilon}\right>_{\bs{\Sigma}_{\bs{\varepsilon}}^{-1}}-2\left<\bs{G}(\bs{\theta}),\bs{G}(\bs{\vartheta})\right>_{\bs{\Sigma}_{\bs{\varepsilon}}^{-1}}-2\left<\bs{G}(\bs{\vartheta}),\bs{\varepsilon}\right>_{\bs{\Sigma}_{\bs{\varepsilon}}^{-1}}\right)}.
\end{align}
 From this, it follows that
\begin{align}\label{eq:integ}
    EIG{}&=\int_{\bs{\Theta}}\int_{\bb{R}^{d_y}} \left[\frac{1}{2}\left\lVert\bs{G}(\bs{\theta})\right\rVert_{\bs{\Sigma}_{\bs{\varepsilon}}^{-1}}^2+\left<\bs{G}(\bs{\theta}),\bs{\varepsilon}\right>_{\bs{\Sigma}_{\bs{\varepsilon}}^{-1}}\right.\nonumber\\
    {}&\left.\quad-\log\left(\int_{\bs{\Theta}}e^{-\frac{1}{2}\left(\left\lVert\bs{G}(\bs{\vartheta})\right\rVert_{\bs{\Sigma}_{\bs{\varepsilon}}^{-1}}^2-2\left<\bs{G}(\bs{\theta}),\bs{G}(\bs{\vartheta})\right>_{\bs{\Sigma}_{\bs{\varepsilon}}^{-1}}-2\left<\bs{G}(\bs{\vartheta}),\bs{\varepsilon}\right>_{\bs{\Sigma}_{\bs{\varepsilon}}^{-1}}\right)}\pi(\bs{\vartheta})\di{}\bs{\vartheta}\right)\right]\nonumber\\
    {}&\quad\quad\times\frac{1}{\det(2\pi\bs{\Sigma}_{\bs{\varepsilon}})^{\frac{1}{2}}}e^{-\frac{1}{2}\left\lVert\bs{\varepsilon}\right\rVert_{\bs{\Sigma}_{\bs{\varepsilon}}^{-1}}^2}\di{}\bs{\varepsilon}\pi(\bs{\theta})\di{}\bs{\theta}.
\end{align}
The first term in \eqref{eq:integ} is not nested and can thus be estimated using the standard rQMC or similar methods. The second term is symmetric around the origin and thus vanishes as it is integrated against the Gaussian measure for $\bs{\varepsilon}$. The remaining term constitutes a nested integral and is the main subject of this investigation. We introduce the notation $f\equiv\log$, and
\begin{equation}\label{eq:g.simplified}
    g(\bs{y},\bs{x})\coloneqq \exp\left({-\frac{1}{2}\left\lVert\bs{G}(F_{\bs{\theta}}^{-1}(\bs{x}))\right\rVert_{\bs{\Sigma}_{\bs{\varepsilon}}^{-1}}^2+\left<\bs{G}(F_{\bs{\theta}}^{-1}(\bs{y}_1)),\bs{G}(F_{\bs{\theta}}^{-1}(\bs{x}))\right>_{\bs{\Sigma}_{\bs{\varepsilon}}^{-1}}+\left<\bs{G}(F_{\bs{\theta}}^{-1}(\bs{x})),\Phi^{-1}(\bs{y}_2)\right>_{\bs{\tilde{\Sigma}}_{\bs{\varepsilon}}^{-\trans}}}\right),
\end{equation}
where $F_{\bs{\theta}}^{-1}$ is the inverse CDF of the parameters of interest, $\bs{y}=(\bs{y}_1,\bs{y}_2)$, $F_{\bs{\theta}}^{-1}(\bs{y}_1)\coloneqq\bs{\theta}$, $F_{\bs{\theta}}^{-1}(\bs{x})\coloneqq\bs{\vartheta}$, $\bs{\tilde{\Sigma}}_{\bs{\varepsilon}}\Phi^{-1}(\bs{y}_2)\coloneqq \bs{\varepsilon}$, and 
\begin{equation}\label{eq:Sigma.tilde}
\bs{\Sigma}_{\bs{\varepsilon}}=\bs{\tilde{\Sigma}}_{\bs{\varepsilon}}\bs{\tilde{\Sigma}}_{\bs{\varepsilon}}^{\trans}.
\end{equation}
Moreover, we introduce the notation
\begin{equation}
    \bs{\tilde{G}}\coloneqq \bs{G}\circ F_{\bs{\theta}}^{-1}.
\end{equation}
The Gaussian noise present in typical data models together with the logarithm present in the Kullback--Leibler divergence \cite{Kul51, Kul59} introduce singularities at the integration domain boundary, rendering the standard error bound via the Koksma--Hlawka inequality~\eqref{eq:QMC.bound} trivial, as the Hardy--Krause variation for such integrands is infinite. The behavior of these singularities must be carefully controlled to obtain the optimal rate of convergence of the rQMC error.
Moreover, for the case where the experiment model $\bs{\tilde{G}}$ must be approximated as $\bs{\tilde{G}}_h$ by an appropriate method with discretization parameter $h>0$, the inner integrand in~\eqref{eq:integ} is instead defined as follows:
\begin{equation}\label{eq:g.h}
        g_h(\bs{y},\bs{x})\coloneqq \exp\left({-\frac{1}{2}\left\lVert\bs{\tilde{G}}_h(\bs{x})\right\rVert_{\bs{\Sigma}_{\bs{\varepsilon}}^{-1}}^2+\left<\bs{\tilde{G}}_h(\bs{y}_1),\bs{\tilde{G}}_h(\bs{x})\right>_{\bs{\Sigma}_{\bs{\varepsilon}}^{-1}}+\left<\bs{\tilde{G}}_h(\bs{x}),\Phi^{-1}(\bs{y}_2)\right>_{\bs{\tilde{\Sigma}}_{\bs{\varepsilon}}^{-\trans}}}\right).
    \end{equation}
The following assumptions are made in this work:
\begin{asu}[Condition on the approximate experiment model]\label{asu:Lipschitz}
    Let $\bs{\tilde{G}}_h:[0,1]^{d_2}\to\bb{R}^{d_1-d_2}$ be an approximation to $\bs{\tilde{G}}$ and assume that there exists $0\leq q<\infty$ independent of $h$ such that
    \begin{equation}\label{eq:mixed.derivatives}
    \left\lVert\left(\prod_{j\in u}\frac{\partial}{\partial z_j}\right)\bs{\tilde{G}}_h(\bs{z})\right\rVert_{2}\leq q,
\end{equation}
for all $\bs{z}\in[0,1]^{d_{2}}$, all $h>0$,
and all $u\subseteq\{1,\ldots,d_2\}$.
\end{asu}
\begin{asu}[Inverse inequality for the FEM rate]\label{asu:inverse.FEM}
Let $g_h$ be the approximate inner integrand in~\eqref{eq:g.h}, with the inner integral denoted as follows:
    \begin{equation}\label{eq:g.h.bar}
        \bar{g}_h(\bs{y})\coloneqq \int_{[0,1]^{d_2}}\exp\left({-\frac{1}{2}\left\lVert\bs{\tilde{G}}_h(\bs{x})\right\rVert_{\bs{\Sigma}_{\bs{\varepsilon}}^{-1}}^2+\left<\bs{\tilde{G}}_h(\bs{y}_1),\bs{\tilde{G}}_h(\bs{x})\right>_{\bs{\Sigma}_{\bs{\varepsilon}}^{-1}}+\left<\bs{\tilde{G}}_h(\bs{x}),\Phi^{-1}(\bs{y}_2)\right>_{\bs{\tilde{\Sigma}}_{\bs{\varepsilon}}^{-\trans}}}\right)\di{}\bs{x},
    \end{equation}
    where $\bs{\tilde{G}}_h$ as in Assumption~\ref{asu:Lipschitz}, $\bs{y}=(\bs{y}_1,\bs{y}_2)$, $\bs{y}_1, \bs{x}\in[0,1]^{d_2}$, $\bs{y}_2\in[0,1]^{d_1-d_2}$, $\bs{\Sigma_{\bs{\epsilon}}}\in\bb{R}^{(d_1-d_2)\times(d_1-d_2)}$ a symmetric, strictly positive definite matrix, and $\bs{\Sigma}_{\bs{\varepsilon}}=\bs{\tilde{\Sigma}}_{\bs{\varepsilon}}\bs{\tilde{\Sigma}}_{\bs{\varepsilon}}^{\trans}$ as in~\eqref{eq:Sigma.tilde}. Then, assume that
        \begin{equation}
    \sup_{\bs{y}\in[0,1]^{d_1}} \left\lvert \frac{\left(\prod_{j\in u}\frac{\partial}{\partial y_{j}}\right)(\bar{g}_{h}(\bs{y})-\bar{g}(\bs{y}))}{\left(\prod_{j\in u}\frac{\partial}{\partial y_{j}}\right) \bar{g}(\bs{y})} \right\rvert < C_{{\rm{s}},d_1}h^{\eta_{{\rm{s}},d_1}},
\end{equation}
is the relative strong convergence rate depending on the regularity of $\bar{g}_h$ and $\bar{g}$ in the $j$th components, where $j\in u\subseteq \{1,\ldots,d\}$, and $0<C_{{\rm{s}},d_1}<\infty$ is a constant independent of $h$, and $\eta_{{\rm{s}},d_1}>0$.
\end{asu}
The bound on the total error for estimating the EIG using the MLDLQMC estimator for any $A_i>0$, where $1\leq i\leq d_1$ and $d_1=d_{\theta}+d_{y}$, is presented below. 
\begin{lem}[Bound on the integrand at level 0 for the EIG]\label{lem:level.0}
Let $f\equiv\log$ and $g_{h}$ as in~\eqref{eq:g.h} for any $h>0$. Moreover, let $\bs{x}^{(m)}$ be as in~\eqref{eq:x.ell.h.det}, where $1\leq m\leq M$. Then, provided that Assumption~\ref{asu:Lipschitz} holds, there exists $0<b<\infty$ such that
\begin{equation}\label{eq:bound.ell.0}
    \left|\left(\prod_{j\in u}\frac{\partial}{\partial y_j}\right)f\left(\frac{1}{M}\sum_{m=1}^{M}g_{h}\left(\bs{y},\bs{x}^{(m)}\right)\right)\right|\leq b\prod_{i=1}^{d_1}\min (y_i,1-y_i)^{-A_i-\mathds{1}_{\{i\in u\}}}
\end{equation}
for all $\bs{y}=(y_1,\ldots,y_{d_1})\in(0,1)^{d_1}$, all $u\subseteq\{1,\ldots,{d_1}\}$, and any $0<\max_{1\leq i\leq d_1}A_i<1/2$.
\end{lem}
Appendix~\ref{app:ell} presents the proof of this result and Appendix~\ref{app:vcond} provides auxiliary results. A truncation of the Gaussian observation noise in the EIG, directly adapted from~\cite[Lemma 3]{Bar23}, is employed to bound the variance of the levels $1\leq \ell\leq L$:
\begin{lem}[Truncation of the Gaussian observation noise in the EIG]\label{lem:truncation}
    Let the EIG be as in \eqref{eq:integ}, that is,
    \begin{align}\label{eq:EIG.x.y}
    EIG{}&=\int_{[0,1]^{d_2}}\int_{\bb{R}^{d_1-d_2}} \left[\frac{1}{2}\left\lVert\bs{\tilde{G}}(\bs{y}_1)\right\rVert_{\bs{\Sigma}_{\bs{\varepsilon}}^{-1}}^2+\left<\bs{\tilde{G}}(\bs{y}_1),\bs{\varepsilon}\right>_{\bs{\Sigma}_{\bs{\varepsilon}}^{-1}}\right.\nonumber\\
    {}&\left.\quad-\log\left(\int_{[0,1]^{d_2}}e^{-\frac{1}{2}\left(\left\lVert\bs{\tilde{G}}(\bs{x})\right\rVert_{\bs{\Sigma}_{\bs{\varepsilon}}^{-1}}^2-2\left<\bs{\tilde{G}}(\bs{y}_1),\bs{\tilde{G}}(\bs{x})\right>_{\bs{\Sigma}_{\bs{\varepsilon}}^{-1}}-2\left<\bs{\tilde{G}}(\bs{x}),\bs{\varepsilon}\right>_{\bs{\Sigma}_{\bs{\varepsilon}}^{-1}}\right)}\di{}\bs{x}\right)\right]\nonumber\\
    {}&\quad\quad\times\frac{1}{\det(2\pi\bs{\Sigma}_{\bs{\varepsilon}})^{\frac{1}{2}}}e^{-\frac{1}{2}\left\lVert\bs{\varepsilon}\right\rVert_{\bs{\Sigma}_{\bs{\varepsilon}}^{-1}}^2}\di{}\bs{\varepsilon}\di{}\bs{y}_1,
\end{align}
and denote by
    \begin{align}\label{eq:EIG.x.y.tr}
    EIG^{\rm{tr}}{}&=\int_{[0,1]^{d_2}}\int_{[-c(TOL),c(TOL)]^{d_1-d_2}} \left[\frac{1}{2}\left\lVert\bs{\tilde{G}}(\bs{y}_1)\right\rVert_{\bs{\Sigma}_{\bs{\varepsilon}}^{-1}}^2+\left<\bs{\tilde{G}}(\bs{y}_1),\bs{\varepsilon}\right>_{\bs{\Sigma}_{\bs{\varepsilon}}^{-1}}\right.\nonumber\\
    {}&\left.\quad-\log\left(\int_{[0,1]^{d_2}}e^{-\frac{1}{2}\left(\left\lVert\bs{\tilde{G}}(\bs{x})\right\rVert_{\bs{\Sigma}_{\bs{\varepsilon}}^{-1}}^2-2\left<\bs{\tilde{G}}(\bs{y}_1),\bs{\tilde{G}}(\bs{x})\right>_{\bs{\Sigma}_{\bs{\varepsilon}}^{-1}}-2\left<\bs{\tilde{G}}(\bs{x}),\bs{\varepsilon}\right>_{\bs{\Sigma}_{\bs{\varepsilon}}^{-1}}\right)}\di{}\bs{x}\right)\right]\nonumber\\
    {}&\quad\quad\times\frac{1}{\det(2\pi\bs{\Sigma}_{\bs{\varepsilon}})^{\frac{1}{2}}}e^{-\frac{1}{2}\left\lVert\bs{\varepsilon}\right\rVert_{\bs{\Sigma}_{\bs{\varepsilon}}^{-1}}^2}\di{}\bs{\varepsilon}\di{}\bs{y}_1,
\end{align}
the $TOL-$truncated EIG. Here,
\begin{equation}\label{eq:c.TOL.q}
    c(TOL)\coloneqq \tilde{q}(2(1+p))^{\frac{1}{2}}\log(TOL^{-1})^{\frac{1}{2}},
\end{equation}
for some constant $0<\tilde{q}<\infty$, and any $p>0$, $TOL>0$, signifies a region of truncation. Then, it holds that
\begin{equation}
    \left|EIG-EIG^{\rm{tr}}\right|=o(TOL)
\end{equation}
as $TOL\to 0$.
\end{lem}
The result in Lemma~\ref{lem:level.0} still applies to the truncated integrands $g^{\rm{tr}}_h$ and $\bar{g}^{\rm{tr}}_h$ and is not restated for brevity. The following notation is introduced:
\begin{equation}\label{eq:M.tilde}
    \tilde{M}(TOL)\coloneqq\left(\frac{C_{\epsilon,d_2}c(TOL)^{d_2}(1+C_{{\rm{s}},d_1}h^{\eta_{{\rm{s}},d_1}})}{(1-C_{{\rm{s}},d_1}h^{\eta_{{\rm{s}},d_1}})}\right)^{\frac{1}{1-\epsilon}},
\end{equation}
where $c(TOL)$ as in~\eqref{eq:c.TOL.q}, $TOL>0$, $0<h<C_{{\rm{s}},d_1}^{-1/\eta_{{\rm{s}},d_1}}$, $C_{{\rm{s}},d_1}>0$, and $\eta_{{\rm{s}},d_1}>0$.
\begin{lem}[Bound on the higher-order term in the Taylor expansion of the nested integrand]\label{lem:boundary.growth.taylor}
    Let $f\equiv \log$, $g_h^{\rm{tr}}$ as in~\eqref{eq:g.h.truncated}, $0<h<C_{{\rm{s}},d_1}^{-1/\eta_{{\rm{s}},d_1}}$, $M>\tilde{M}(TOL)$, where $\tilde{M}(TOL)$ as in~\eqref{eq:M.tilde}, $TOL>0$, and $\bs{x}^{(m)}$ be as in~\eqref{eq:x.ell.h.det}, where $1\leq m\leq M$. Moreover, let $\varphi_{h,M}$ constitute an auxiliary function, where
    \begin{equation}
        \varphi_{h,M}(\cdot)=\left(\frac{1}{M}\sum_{m=1}^Mg^{\rm{tr}}_h(\cdot,\bs{x}^{(m)})-\bar{g}^{\rm{tr}}(\cdot)\right)^2\int_0^1f''\left(\bar{g}^{\rm{tr}}(\cdot)+s \left(\frac{1}{M}\sum_{m=1}^Mg^{\rm{tr}}_h(\cdot,\bs{x}^{(m)})-\bar{g}^{\rm{tr}}(\cdot)\right)\right)(1-s)\di{}s,
    \end{equation}
    and $\bar{g}^{\rm{tr}}\coloneqq \lim_{h\to 0}\bar{g}^{\rm{tr}}_h$.
Then, given Assumptions~\ref{asu:Lipschitz} and~\ref{asu:inverse.FEM}, there exists $0<b<\infty$ such that
    \begin{equation}\label{eq:boundary.growth.taylor}
    \left|\left(\prod_{j\in u}\frac{\partial^{u}}{\partial y_j}\right)\varphi_{h,M}(\bs{y})\right|\leq \left(C_{\epsilon,d_2}^2M^{-2+2\epsilon}+C_{\eta_{{\rm{s}},d_1}}^2h^{2\eta_{{\rm{s}},d_1}}\right)b\prod_{i=1}^{d_1}\min (y_i,1-y_i)^{-A_i-\mathds{1}_{\{i\in u\}}}
\end{equation}
for all $\bs{y}=(y_1,\ldots,y_{d_1})\in(0,1)^{d_1}$, all $u\subseteq\{1,\ldots,{d_1}\}$, and any $0<A_i<1/2$, where $1\leq i\leq d_1$.
\end{lem}
Appendix~\ref{app:ell} presents the proof of this result. The condition that $M>\tilde{M}(TOL)$ for $TOL>0$ and $h<C_{{\rm{s}},d_1}^{-1/\eta_{{\rm{s}},d_1}}$ in Lemma~\ref{lem:boundary.growth.taylor} is addressed by setting the following in the MLDLQMC estimator:
\begin{equation}
    h_0\coloneqq \delta_{h}C_{{\rm{s}},d_1}^{-\frac{1}{\eta_{{\rm{s}},d_1}}},
\end{equation}
where $0<\delta_{h}<1$, and
\begin{equation}
    M_{0}\coloneqq \left(\frac{C_{\epsilon,d_2}(1+C_{{\rm{s}},d_1}h_0^{\eta_{{\rm{s}},d_1}})}{(1-C_{{\rm{s}},d_1}h_0^{\eta_{{\rm{s}},d_1}})}\right)^{\frac{1}{1-\epsilon}}c(TOL)^{\frac{d_2}{1-\epsilon}}.
\end{equation}
The bound~\eqref{eq:bound.ell.0} in Lemma~\ref{lem:level.0} is independent of $M$ and $h$ and thus unaffected by these conditions. The multiplicative factor $c(TOL)^{d_2/(1-\epsilon)}$ is tracked throughout the remaining error bounds in the following corollary.
\begin{cor}[Total error of the MLDLQMC estimator for EIG estimation]\label{cor:EIG.total.error}
    Given Assumption~\ref{asu:weak.rate}, Assumption~\ref{asu:Lipschitz}, and Assumption~\ref{asu:inverse.FEM}, for $M_{\ell}>\tilde{M}(TOL)$ and $h_{\ell}<C_{{\rm{s}},d_1}^{-1/\eta_{{\rm{s}},d_1}}$, where $\tilde{M}(TOL)$ is as in~\eqref{eq:M.tilde} and $0\leq\ell\leq L$, there exist $0<\delta_{K}<1/2$ such that the total error of the MLDLQMC estimator for the EIG as in \eqref{eq:integ} is bounded as follows:
    \begin{equation}\label{eq:total.err.tr}
        \left|EIG-EIG^{\rm{tr}}\right|=o(TOL),
    \end{equation}
\begin{align}\label{eq:Bias.constraint.h.det.EIG}
  |\mathbb{E}[I_{\rm{MLDLQ}}]-EIG^{\rm{tr}}|\leq C_{\mathrm{w}}h_L^{\eta_{\rm{w}}}{}&+ \frac{\bb{E}[|\bar{g}^{\rm{tr}}_{h}|^2\left|f''(\bar{g}^{\rm{tr}}_{h})\right|]c(TOL)^{2d_2}C_{\epsilon,d_2}^2}{2M_L^{2-2\epsilon}}\nonumber\\
  {}&+\frac{\bb{E}[\left|\bar{g}^{\rm{tr}}_{h}\right|^3\left|f'''(K(TOL)\bar{g}^{\rm{tr}}_{h})\right|]c(TOL)^{3d_2}C_{\epsilon,d_2}^3}{6M_L^{3-3\epsilon}},
\end{align}
\begin{align}\label{eq:total.err.var}
    \bb{V}\left[I_{\rm{MLDLQ}}\right]
    \leq \frac{(b^{(0)})^2B_A^2C_{\epsilon, d_1}^2}{N_0^{2-2\epsilon-2A_{\rm{max}}}} {}&+\sum_{\ell=1}^{L}\left(\frac{\left(\frac{C_{\epsilon,d_2}^2c(TOL)^{2d_2}}{M_{\ell}^{2-2\epsilon}}+C_{\eta_{{\rm{s}},d_1}}^2h_{\ell}^{2\eta_{{\rm{s}},d_1}}\right)(b^{(I)})^2B_A^2C_{\epsilon,d_1}^2}{N_{\ell}^{2-2\epsilon-2A_{\rm{max}}}}\right.\nonumber\\
    {}&\left.\quad\quad +\frac{\left(\frac{C_{\epsilon,d_2}^4c(TOL)^{4d_2}}{M_{\ell}^{4-4\epsilon}}+C_{\eta_{{\rm{s}},d_1}}^4h_{\ell}^{4\eta_{{\rm{s}},d_1}}\right)(b^{(II)})^2B_A^2C_{\epsilon,d_1}^2}{N_{\ell}^{2-2\epsilon-2A_{\rm{max}}}}\right)\nonumber\\
        {}&+\frac{C_{\epsilon,d_2}^2(b^{(III)})^2B_A^2c(TOL)^{2d_2}}{M_L^{2-2\epsilon}}+\frac{C_{\epsilon,d_2}^4(b^{(IV)})^2B_A^2c(TOL)^{4d_2}}{M_L^{4-4\epsilon}},
\end{align}
    for any $\epsilon>0$ and all $K(TOL)\in(\delta_{K},1-\delta_{K})$ for an error tolerance $TOL>0$, where $c(TOL)$ as in \eqref{eq:c.TOL.q}, where $C_{\epsilon,d_1}, C_{\epsilon,d_2}\to\infty$ as $\epsilon\to 0$, $B_A\to\infty$ as $\min_{1\leq i\leq d_1}A_i\to 0$, and any $0<A_{\rm{max}\coloneqq \max_{1\leq i\leq d_1}A_i<1/2}$ where $A_i$ as in~\eqref{eq:bound.ell.0} and~\eqref{eq:boundary.growth.taylor}, where $1\leq i\leq d_1$.
\end{cor}
\begin{proof}
    The result in \eqref{eq:total.err.tr} follows from Lemma~\ref{lem:truncation}. Assumption~\ref{asu:f.3.diff} follows from Corollary~\ref{cor:bound.nested.integrand} and Assumption~\ref{asu:g.L4} follows from Lemma~\ref{cor:inverse.g.h} for $k$ replaced with $c(TOL)^{d_2}$ and Corollary~\ref{cor:bar.g} as $\max_{i}A_i<1/2$. It is verified that
    \begin{equation}
        |f'''(K(TOL)\bar{g}^{\rm{tr}}_{h})|=\left|\frac{2}{(K(TOL))^3(\bar{g}^{\rm{tr}}_{h})^3}\right|
    \end{equation}
 is monotonic for positive $K(TOL)\bar{g}^{\rm{tr}}_{h}$ and $f\equiv\log$. The result in \eqref{eq:Bias.constraint.h.det.EIG} thus follows from Proposition~\ref{prop:B.MLrQ-Q.h}. Assumption~\ref{asu:boundary.0} follows from Lemma~\ref{lem:level.0}. Assumption~\ref{asu:boundary.ell} follows from Lemma~\ref{lem:boundary.growth.taylor} with the additional condition that $M_{\ell}>\tilde{M}(TOL)$ and $h_{\ell}<C_{{\rm{s}},d_1}^{-1/\eta_{{\rm{s}},d_1}}$, where $0\leq\ell\leq L$. Assumption~\ref{asu:boundary.ell.l.o.t.}, Assumption~\ref{asu:boundary.L.l.o.t.}, and Assumption~\ref{asu:boundary.L} also follow from Lemma~\ref{lem:boundary.growth.taylor} with a similar condition on $M_{\ell}$ and $h_{\ell}$. 
 Finally, the result in \eqref{eq:total.err.var} follows from Propositions~\ref{prop:V.MLrQ-Q.0} to \ref{prop:V.ell.MLQ.h.L}, again with the condition that $M_L>\tilde{M}(TOL)$ and $h_L<C_{{\rm{s}},d_1}^{-1/\eta_{{\rm{s}},d_1}}$.
\end{proof}
\begin{rmk}[Total complexity of the MLDLQMC estimator for EIG estimation]
    The total computational complexity of the MLDLQMC estimator applied to the EIG differs from the result stated in Proposition~\ref{prop:Work.MLQ.h.det} by a $\log$ factor as a result of the truncation of the observation noise. 
\end{rmk}
\begin{rmk}[Total complexity of the MLDLQMC estimator for scrambled Sobol' sequence]
    The convergence rate of the inner sampling bias in~\eqref{eq:Bias.constraint.h.det.EIG} and the inner variance at level $L$ in~\eqref{eq:total.err.var} improve for rDLQMC estimators based on scrambled Sobol' sequences. Thus, the exponent $2-2\epsilon$ of $M_L$ in the above expressions is replaced by $3-2\epsilon$ for appropriately adjusted multiplicative constants depending on the generalized Hardy--Krause variation of order one of the inner integrand $g_h^{\rm{tr}}$. The assumptions of Corollary~\ref{cor:EIG.case} relating to the generalized Hardy--Krause variation of order one for the EIG application follow from~\cite[Lemma~4]{Bar23}.
\end{rmk}
\section{Numerical results}\label{sec:numerical.results}
The derivation in~\eqref{eq:g.simplified} appears to be required for technical reasons only (see Appendix~\ref{app:vcond}), and the following inner integrand was applied in the implementation:
\begin{equation}
    g(\bs{y},\bs{x})\coloneqq \exp\left(-\frac{1}{2}\left\lVert\bs{Y}-\bs{G}(\bs{\vartheta})\right\rVert_{\bs{\Sigma}_{\bs{\varepsilon}}^{-1}}^2\right),
\end{equation}
 This integrand displayed preferable numerical behavior to that in~\eqref{eq:g.simplified}. All numerical experiments were performed without truncation of the Gaussian observation noise, indicating that the MLDLQMC estimator is capable of addressing such integrands for the considered tolerances. Moreover, a deterministic shift modulo 1 was applied to the inner points, that is,
\begin{equation}
    \bs{\theta}^{(m)}=\left(\bs{t}^{(m)}+\frac{1}{\sqrt{2}}\right)_{\rm{mod 1}}, \quad 1\leq m\leq M_{\ell}, \, 0\leq \ell\leq L-1,
\end{equation}
to ensure that the initial point $\bs{\theta}^{(1)}$ is not the origin $\bs{0}$, where $(\cdot)_{\rm{mod1}}$ denotes the modulo 1 operator. The choice of $1/\sqrt{2}$ is arbitrary, and any irrational number would serve a similar purpose of moving the Sobol' points away from the origin.

\subsection{Example 2: linear EIG example with exact sampling}
\subsubsection{Setting and notation}
The second example we consider is given as follows:
\begin{equation}\label{eq:model.ex1}
    \bs{G}(\bs{\theta})\coloneqq \bs{A}\bs{\theta},
\end{equation}
where $\bs{\theta}\stackrel{iid}{\sim}\cl{U}([0,1]^{d_{\theta}})$, for $d_{\theta}=4$, and $\bs{A}\coloneqq \bs{I}_{d_{\theta}\times d_{\theta}}$ is the identity matrix in $\bb{R}^{d_{\theta}}$. Moreover,
\begin{equation}
    \bs{Y}=\bs{G}(\bs{\theta}) +\bs{\varepsilon},
\end{equation}
where $\bs{\varepsilon}\stackrel{iid}{\sim}\cl{N}(\bs{0},\sigma^2\bs{I}_{d_y\times d_y})$, for $\sigma^2=0.25$ and $\bs{I}_{d_y\times d_y}$ is the identity matrix in dimension $d_y=4$. Thus, in the nested integration setting, $d_1=d_{\theta}+d_y=8$ and $d_2=d_{\theta}=4$. A design $\bs{\xi}$ is not specified for this example to focus on the demonstration of applicability of the MLDLQMC estimator to the EIG of an experiment. Optimizing the EIG regarding a design is outside the scope of this study. It is easily verified that Assumption~\ref{asu:Lipschitz} holds for $\bs{G}_h\equiv\bs{G}$ in expression~\eqref{eq:model.ex1}. The EIG for this linear example separates into four independent inner integrals in one dimension, allowing for efficient estimation via deterministic quadrature methods. The MLDLQMC estimator with exact sampling was nevertheless able to efficiently approximate the four-dimensional integration problem, as evidenced in the following sections.

\subsubsection{Pilot runs}
\begin{figure}[ht]
	\subfloat[The observed convergence rate of the variance of the inner samples at the final level $L$ of 3.]{%
		\includegraphics[width=0.45\textwidth]{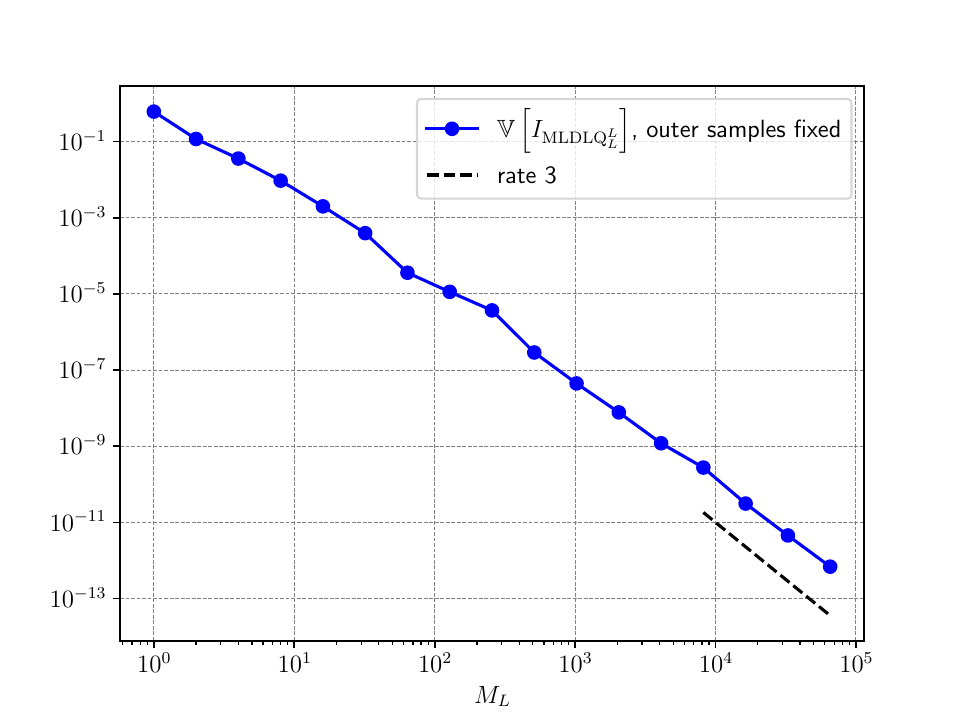}
	}
	\subfloat[The observed convergence rate of the variance at a rate of 2 with a reduced multiplicative factor for the difference estimator.]{%
		\includegraphics[width=0.45\textwidth]{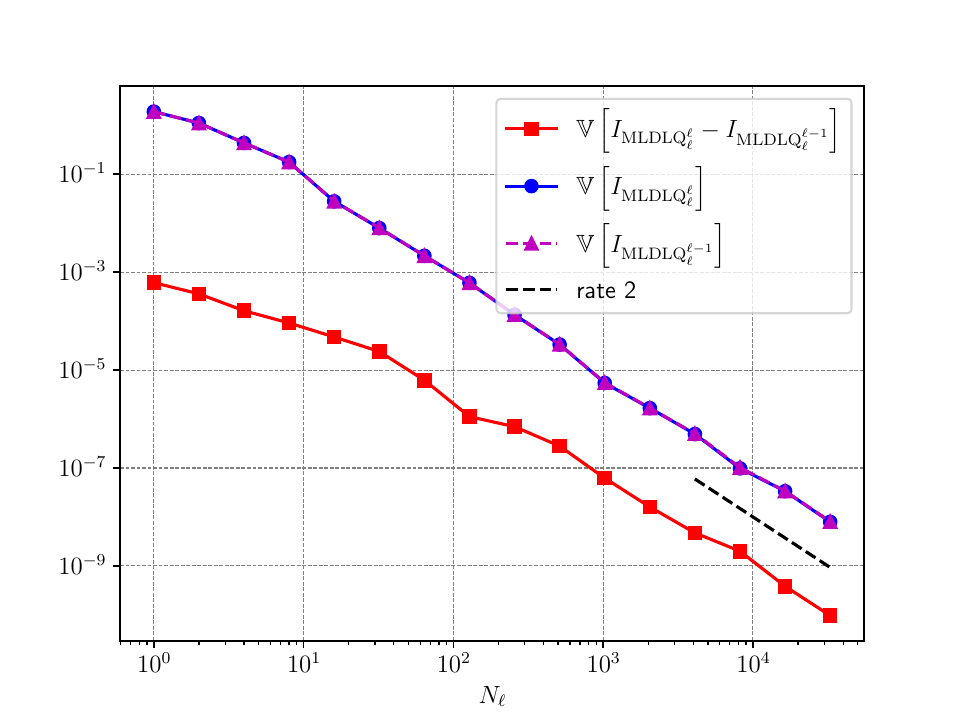}
	}\\
 \subfloat[The variance of the level difference converged at a rate of 2 as a function of $M_{\ell}$ for $N_{\ell}=1$ fixed.]{%
		\includegraphics[width=0.45\textwidth]{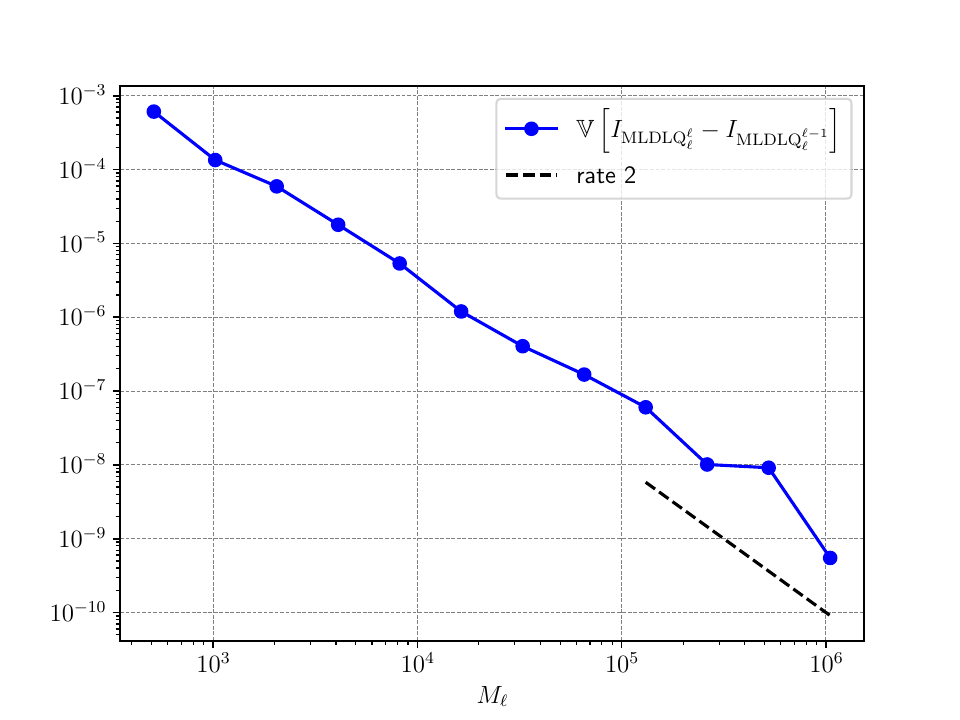}
	}
	\caption{Example 2: Pilot runs for the MLDLQMC estimator with exact sampling. Each observed variance contribution converged at least at the rate predicted by the theoretical analysis.}
	\label{fig:linear.pilot}
\end{figure}

We employed pilot runs to verify that the convergence rates predicted via the theoretical analysis were achieved. This verification also allowed us to estimate the constant terms in the error bounds appearing in Corollary~\ref{cor:EIG.total.error} required for the optimal setting specified in Corollary~\ref{cor:EIG.case}. First, we verified that the variance contribution from the inner samples at the final level $L$ converged at a rate of 3 (Figure~\ref{fig:linear.pilot} (A)). This estimation was carried out using $N_L=1$ outer sample and $S=1$, $R=2^{12}$ randomizations. The convergence of the variance of the level differences as a function of $N_{\ell}$ at the expected rate of 2 was demonstrated in Panel (B), where $M_{\ell-1}=2^{8}$, $M_{\ell}=2^{9}$ samples and $S=2^{10}$ randomizations were used. The variance of the single-level estimators converged at the same rate with a larger multiplicative constant and was independent of $M_{\ell}$, $M_{\ell-1}$. To demonstrate the convergence of the variance of the level difference as a function of $M_{\ell}$, we kept $N_{\ell}=1$ fixed. This variance converged at a rate of 2, as seen in Panel (C) for $S=2^{10}$ randomizations.

\subsubsection{Consistency and work analysis}
\begin{figure}[ht]
\subfloat[The optimal number of levels $L^{\ast}$ increased at a rate of $\log_2(TOL^{-2/3})$ as a function of the error tolerance $TOL>0$.]{%
		\includegraphics[width=0.45\textwidth]{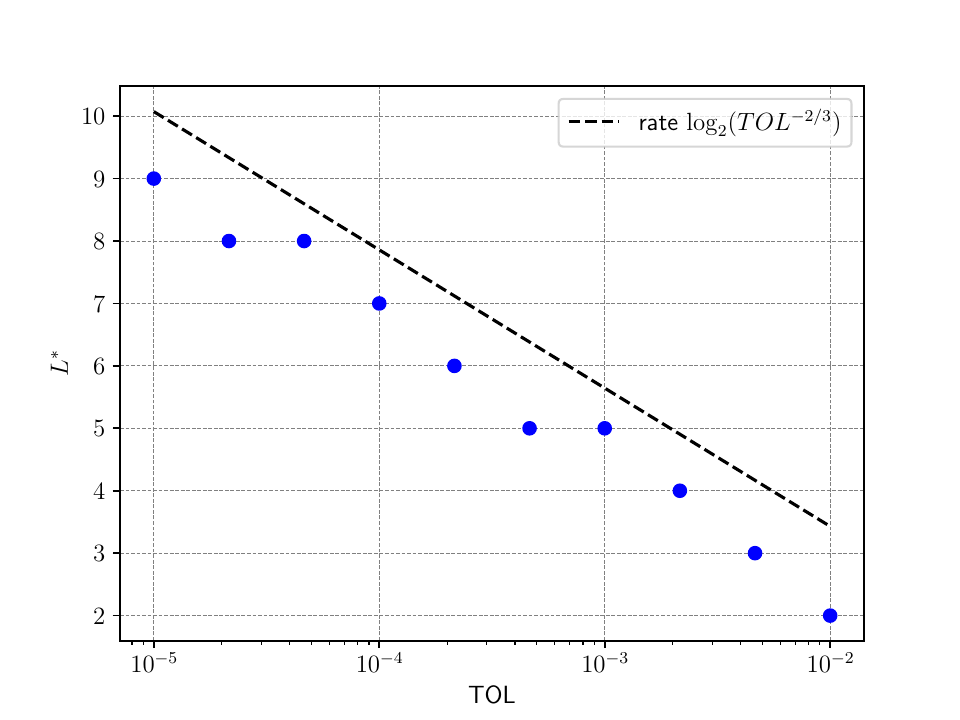}
	}
 \subfloat[Total computational work of the MLDLQMC estimator as a function of the error tolerance $TOL$. The projected rate of approximately 1 was only attained for low tolerances.]{%
		\includegraphics[width=0.45\textwidth]{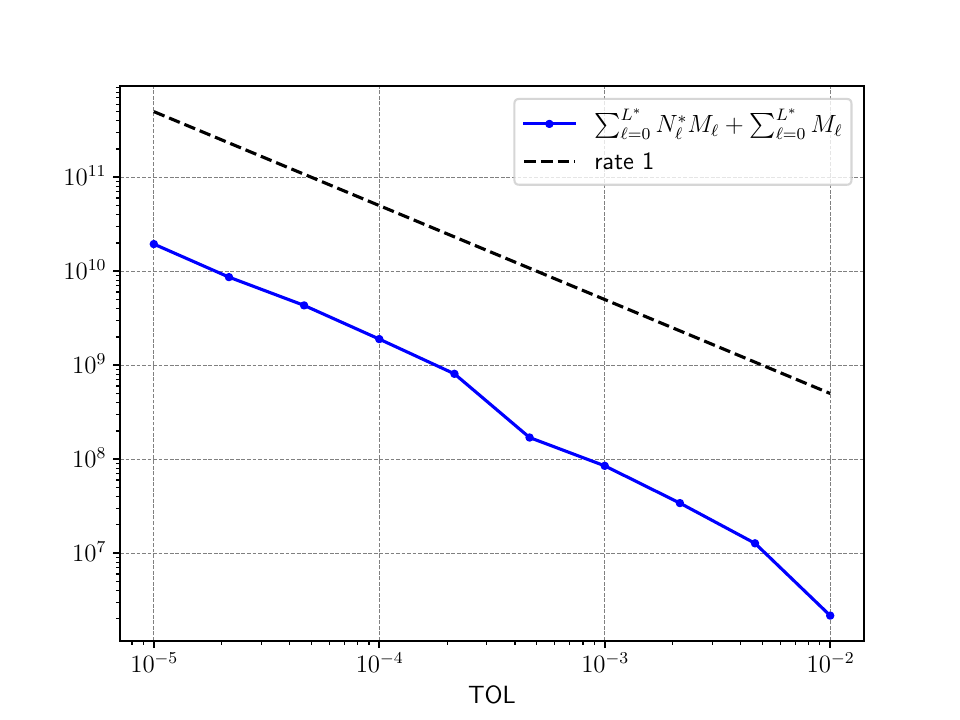}
	}\\
	\subfloat[Error vs.~tolerance $(TOL)$ for 100 runs of the MLDLQMC estimator. Each run resulted in an error less than the tolerance projected by the central limit theorem.]{%
		\includegraphics[width=0.45\textwidth]{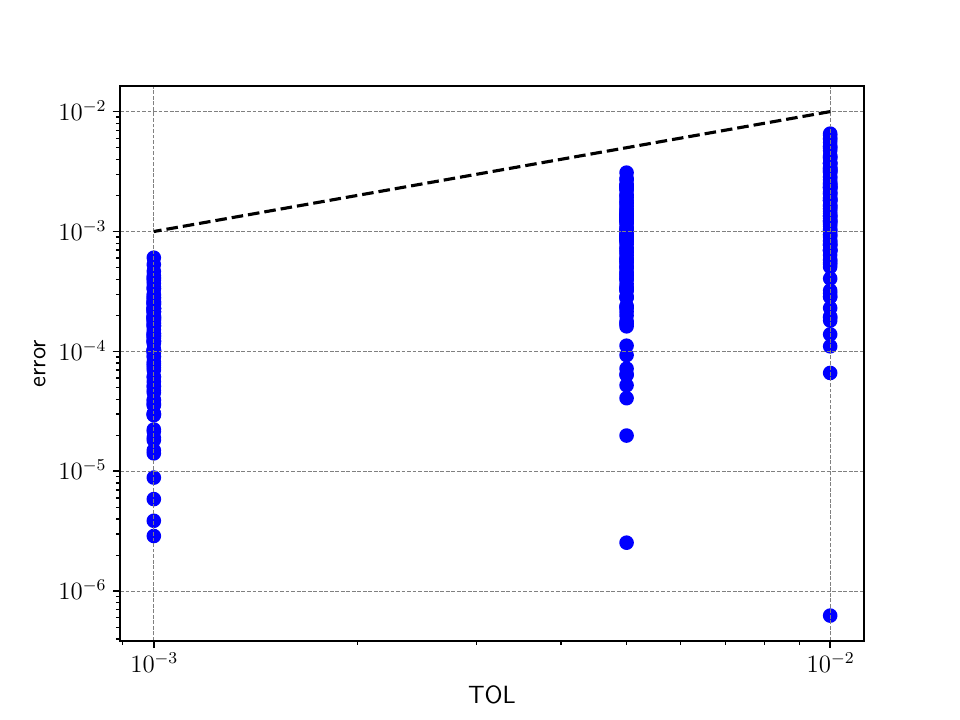}
	}
	
	\caption{Example 2: The MLDLQMC estimator with exact sampling resulted in error less than the projected tolerance $TOL$ (C) at a cost increasing at the projected rate of approximately 1 asymptotically as $TOL\to 0$ (B).}
	\label{fig:linear.evt.work}
\end{figure}

We then determined the optimal final level $L^{\ast}$ of the MLDLQMC estimator with exact sampling for a given tolerance $TOL>0$. The number of levels is displayed in Figure~\ref{fig:linear.evt.work} (A) and the rate of increase is $\log_2(TOL^{-2/3})$, as expected.
Finally, the accuracy and efficiency of the MLDLQMC estimator for optimal outer samples $N_{\ell}^{\ast}$, where $1\leq \ell\leq L^{\ast}$, was analyzed for different error tolerances $TOL$. In Panel (B), the total work of the MLDLQMC estimator, depending on the optimal $L^{\ast}$ and $N_{\ell}^{\ast}$, where $0\leq \ell\leq L^{\ast}$, is presented as a function of the tolerance $TOL$. The projected rate of increase of approximately one was only observed for $TOL<10^{-4}$. In Panel (C), 100 runs of the MLDLQMC estimator with $S=R=1$ randomizations at each level for $TOL=10^{-2}$, $TOL=5\cdot 10^{-3}$, and $TOL=10^{-3}$ are presented, each of them staying below the error tolerance projected by the CLT for the given $L^{\ast}$ and $N_{\ell}^{\ast}$, where $0\leq\ell\leq L^{\ast}$ and $M_0=2^{8}$. An average of 100 runs for $TOL=5\cdot 10^{-4}$ was applied as a reference solution. This reference solution was approximately equal to $I\approx 0.575$ and the smallest relative tolerance depicted in Panel (C) is therefore of order $\cl{O}(10^{-3})$.

\FloatBarrier
\subsection{Example 3: Poisson EIG example with inexact sampling}
\subsubsection{Setting and notation}
The third example we consider is given as follows:
\begin{equation}
    \bs{G}(\bs{\theta},\bs{\xi})\coloneqq \begin{pmatrix}
        u(\xi_1),u(\xi_2)
    \end{pmatrix},
\end{equation}
where $u\in H_0^1(\cl{D})$ is the solution to the Poisson equation in one dimension, that is,
\begin{align}\label{eq:Poisson}
    -u''{}&=\psi\quad \text{in } \cl{D},\nonumber\\
    u{}&=0, \quad \text{on } \{0,1\}.
\end{align}
Here, $\cl{D}=(0,1)$ and
\begin{equation}
    \psi(x,\bs{\theta})\coloneqq\theta_1\sin(4x)+\theta_2\cos(4x)+\theta_3\sin(8x)+\theta_4\cos(8x), \quad x\in\cl{D},
\end{equation}
where $\bs{\theta}=(\theta_1,\theta_2,\theta_3,\theta_4)\stackrel{iid}{\sim}\cl{U}([0,1]^{d_{\theta}})$, for $d_{\theta}=4$. Moreover, $\bs{\xi}\in(0,1)^{d_{\xi}}$, where $d_{\xi}=2$, and
\begin{equation}
    \bs{Y}=\bs{G}(\bs{\theta},\bs{\xi})+\bs{\varepsilon},
\end{equation}
where $\bs{\varepsilon}\stackrel{iid}{\sim}\cl{N}(\bs{0},\sigma^2\bs{I}_{d_y\times d_y})$, for $\sigma^2=10^{-4}$ and $\bs{I}_{d_y\times d_y}$ is the identity matrix in dimension $d_y=2$. Thus, in the nested integration setting, $d_1=d_{\theta}+d_y=6$ and $d_2=d_{\theta}=4$. We consider the discretization
\begin{equation}\label{eq:model.ex2}
    \bs{G}_h(\bs{\theta},\bs{\xi})\coloneqq \begin{pmatrix}
        u_h(\xi_1),u_h(\xi_2)
    \end{pmatrix},
\end{equation}
where $u_h(x):[0,1]\mapsto\bb{R}$ is an FEM approximation of $u$, where $1/h$ specifies the number of elements for $h>0$. Here, we rounded up to the nearest power of two for convenience. Thus, for this example, evaluating the experiment model $\bs{G}_h$ involves the application of an FEM discretization. It is easily verified that Assumption~\ref{asu:Lipschitz} holds for $\bs{G}_h$ in expression~\eqref{eq:model.ex2}. Verifying Assumption~\ref{asu:inverse.FEM}, concerning the strong convergence of mixed derivatives of the inner integrand in a relative sense is contingent on the particular discretization method and will thus be reserved for future research. However, the convergence of the estimator bias as a function of the FEM discretization parameter was demonstrated numerically (in the weak sense) in Figure~\ref{fig:poisson.pilot} (A), and the convergence of the estimator variance (in the strong sense) was indicated in Panel~(E).

\begin{figure}[ht]
\subfloat[Exact solution $u$ of the Poisson equation~\eqref{eq:Poisson} in one dimension as a function of $\bs{\theta}$ and locations of the observations $\bs{\xi}$.]{%
		\includegraphics[width=0.45\textwidth]{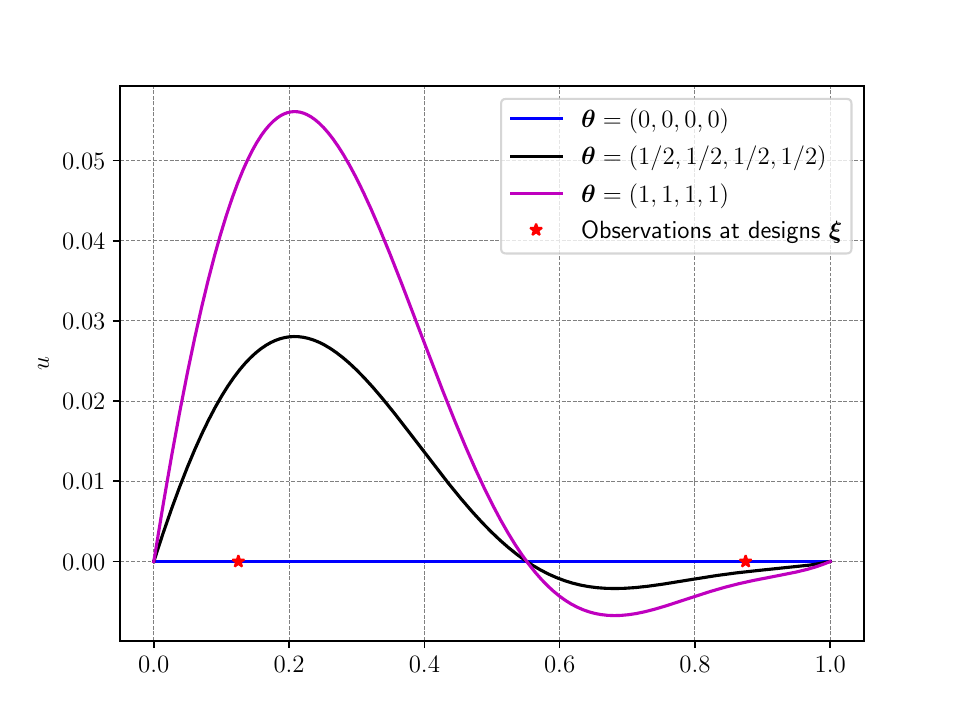}
	}	
 \subfloat[The observed weak convergence rate $\eta_{\rm{w}}=2$ of the FEM solution $u_h$ of the Poisson equation~\eqref{eq:Poisson} in one dimension with linear elements as a function of the discretization parameter $h$.]{%
		\includegraphics[width=0.45\textwidth]{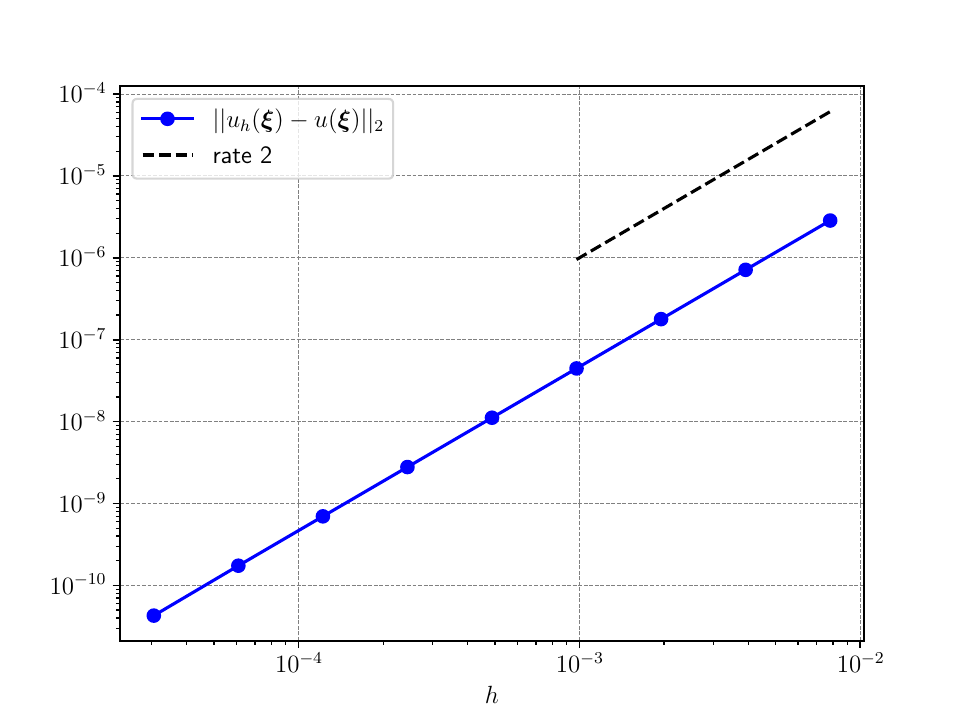}
	}\\
	\subfloat[Runtime in seconds to obtain the FEM solution $u_h$ of the Poisson equation~\eqref{eq:Poisson} yields the work rate $\gamma\approx 1$ as a function of the discretization parameter $h$ (averaged over 20 runs).]{%
		\includegraphics[width=0.45\textwidth]{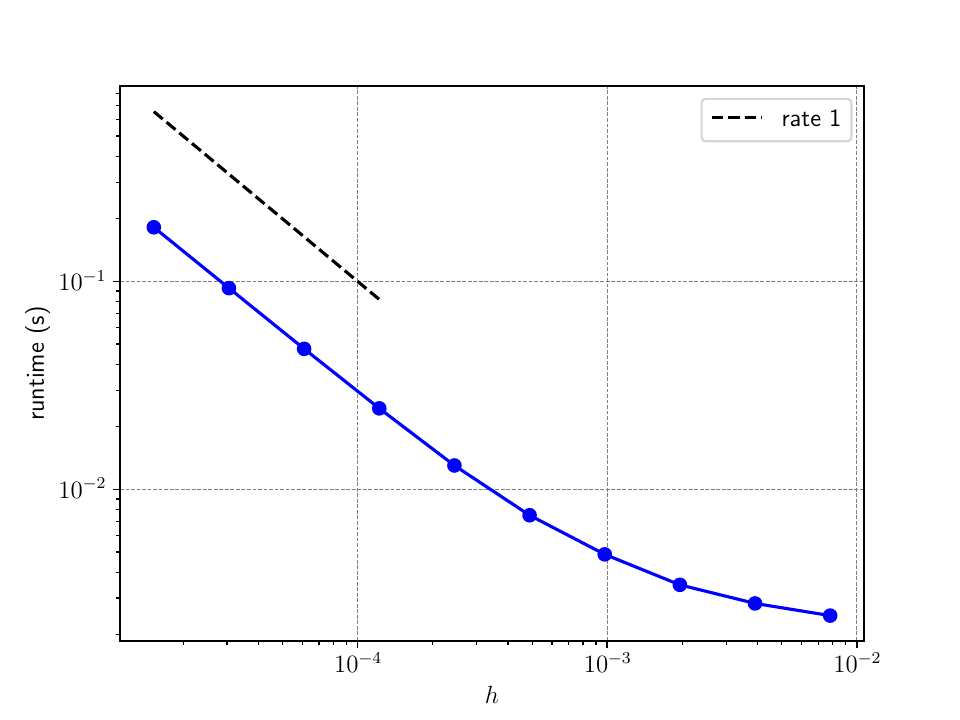}
	}
	\caption{Example 3: FEM approximation of the Poisson equation~\eqref{eq:Poisson} in one dimension with linear elements, where the number of elements is specified as $1/h$.}
	\label{fig:poisson.convergence}
\end{figure}

The Poisson equation~\eqref{eq:Poisson} can be solved exactly, however, for the sake of demonstration, we use a naive FEM implementation. The exact solution $u$ is presented in Figure~\ref{fig:poisson.convergence} (A) for three choices of $\bs{\theta}$. The design is considered fixed at $\bs{\xi}=(1/8, 7/8)$. Using this solution, we were able to verify that the bias of this FEM discretization converged at an observed rate of 2 (Figure~\ref{fig:poisson.convergence} (B)). Moreover, this method has a computational cost of $\cl{O}(h^{-\gamma})$, where $\gamma$ was estimated to be 1 (Figure~\ref{fig:poisson.convergence} (C)).

\subsubsection{Pilot runs}
\begin{figure}[ht]
	\subfloat[The observed weak convergence rate $\eta_{\rm{w}}$ of the discretization bias of the EIG estimation was 2.]{%
		\includegraphics[width=0.45\textwidth]{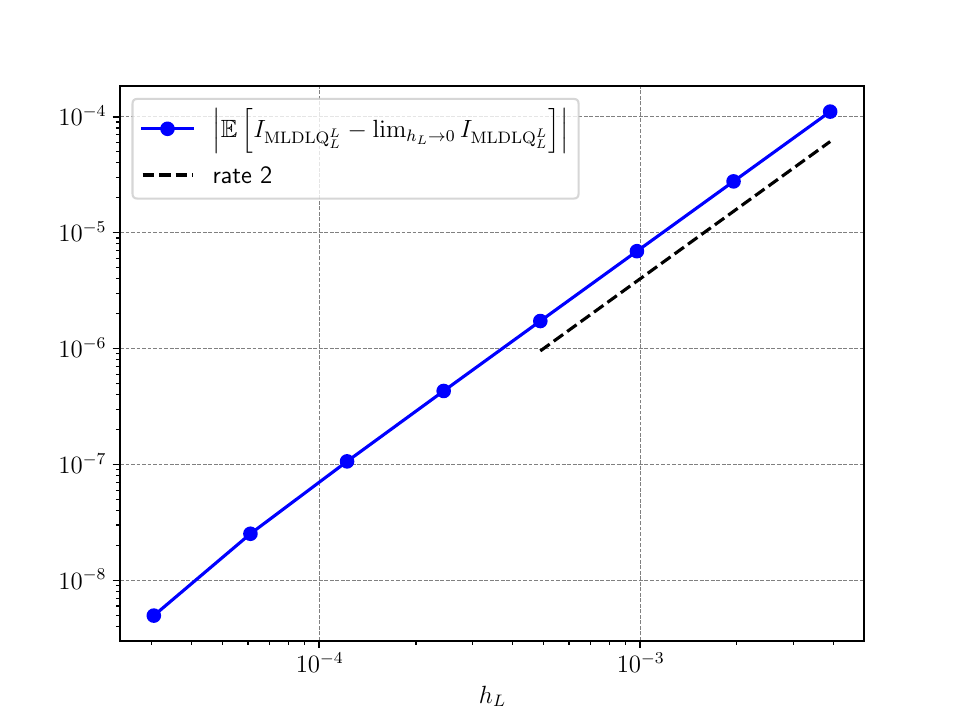}
	}
	\subfloat[The observed convergence rate of the variance of the inner samples at the final level $L$: $3$.]{%
		\includegraphics[width=0.45\textwidth]{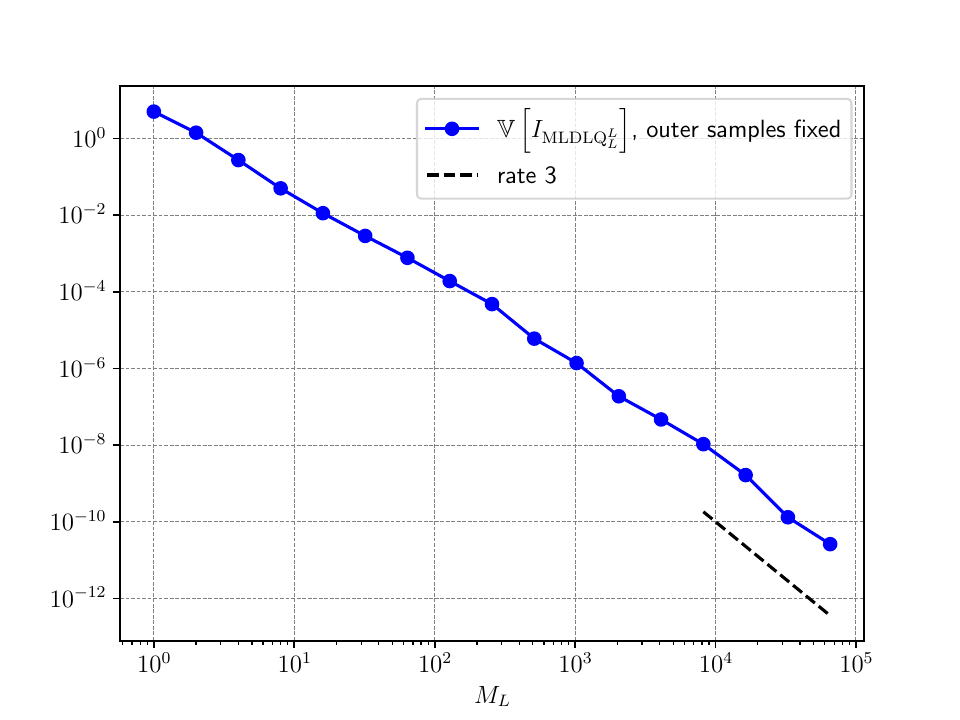}
	}\\
 \subfloat[The observed convergence rate of the variance of 2, with a reduced multiplicative factor for the difference estimator.]{%
		\includegraphics[width=0.45\textwidth]{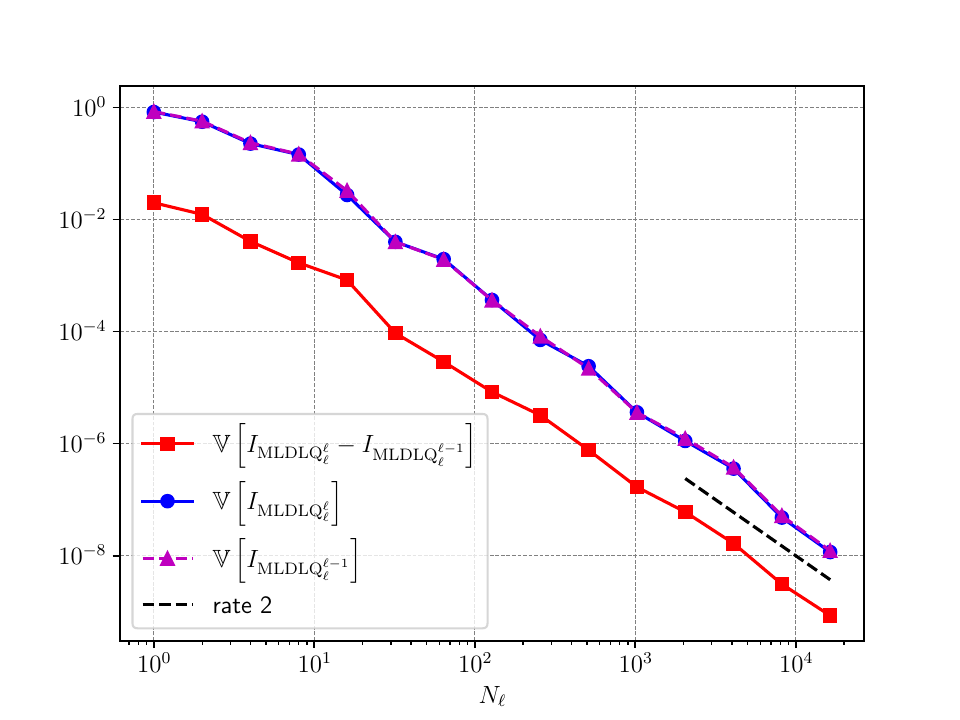}
	}
 \subfloat[Convergence of the level difference as a function of $M_{\ell}$ for $N_{\ell}$ and $h_{\ell}$ fixed at a rate of 2.]{%
		\includegraphics[width=0.45\textwidth]{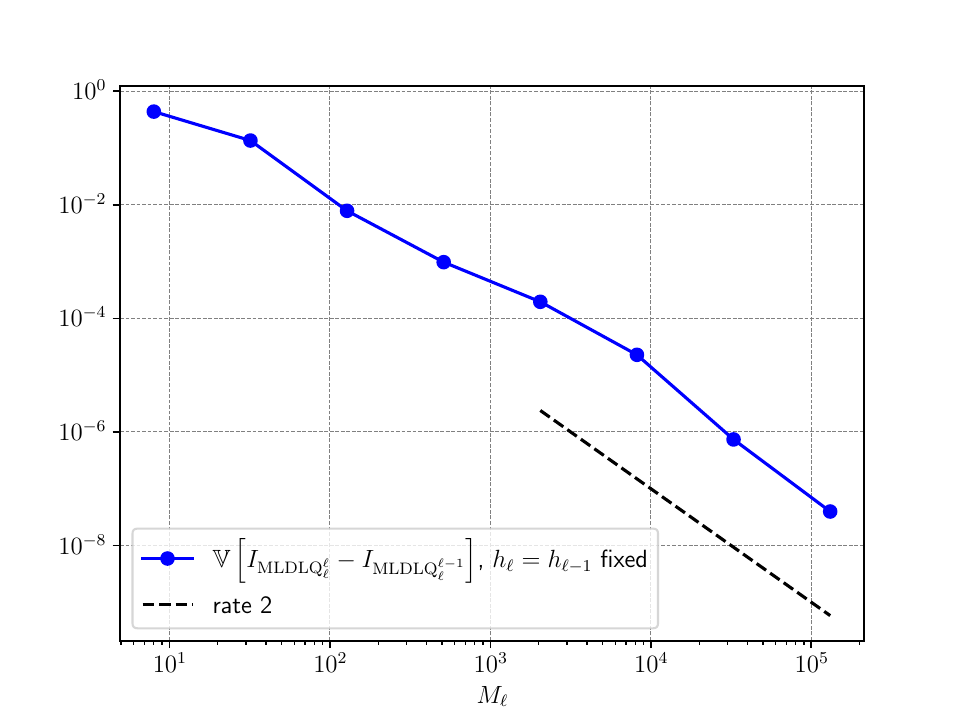}
	}\\
	\subfloat[Convergence of the level difference at a rate of $2\eta_{{\rm{s}},d_1}\approx4$ as a function of $h_{\ell}$ for fixed values of $N_{\ell}$ and $M_{\ell}$.]{%
		\includegraphics[width=0.45\textwidth]{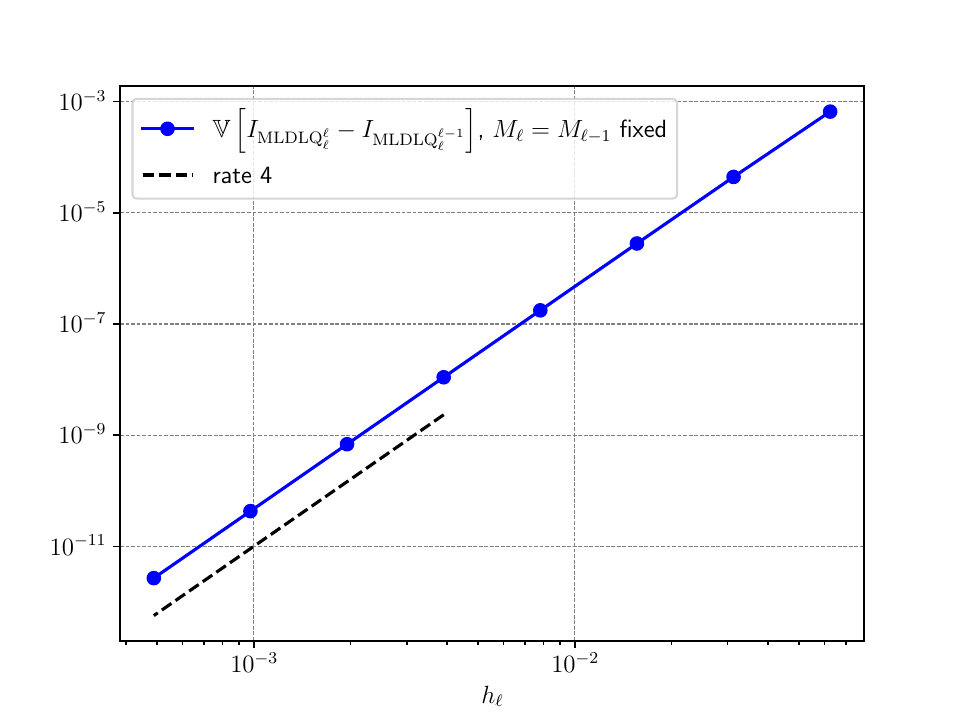}
	}
 \caption{Example 3: Pilot runs for the MLDLQMC estimator with inexact sampling.}
	\label{fig:poisson.pilot}
\end{figure}

We again verified the convergence rates predicted by the theoretical analysis using pilot runs. This verification enabled us to estimate the constant terms in the error bounds appearing in Corollary~\ref{cor:EIG.total.error} required for the optimal setting specified in Corollary~\eqref{cor:EIG.case}. First, we verified that the bias induced by the FEM approximation converged at a rate of $\eta_{\rm{w}}\approx 2$ by fixing the inner number of samples $M_L$ at the final level $L$ as presented in Figure~\ref{fig:poisson.pilot} (A). For this estimation, $N_L=1$, $M_L=2^5$ samples, and $S=100$, $R=1$ randomizations were used, and the bias was estimated using $h_L=2^{-16}$ as a reference solution. To demonstrate generality in our estimations, we preferred this reference solution to the exact solution as the latter may not be available in general. Next, we confirmed that the variance contribution from the inner samples at the final level $L$ converged at a rate of 3 (Panel (B)). This estimation was conducted using $N_L=1$ outer sample, $h_L=2^{-4}$, and $S=2^{7}$, $R=1$ randomizations. The convergence of the variance of the level differences as a function of $N_{\ell}$ at the expected rate of 2 is illustrated in Panel (C), where $M_{\ell-1}=2^{9}$, $M_{\ell}=2^{10}$ samples, $h_{\ell-1}=2^{-4}$, $h_{\ell}=2^{-5}$, and $S=2^{10}$ randomizations were applied. The variances of the single-level estimators converged at the same rate with a larger multiplicative constant and were independent of $M_{\ell-1}$, $M_{\ell}$ and $h_{\ell-1}$, $h_{\ell}$. To demonstrate the convergence of the variance of the level difference as a function of $M_{\ell}$, we kept $N_{\ell}=1$ and $h_{\ell}=h_{\ell-1}=2^{-4}$ fixed. This variance converged at rate 2, as seen in Panel (D) for $S=100$ randomizations. Similarly, we kept $M_{\ell}=M_{\ell-1}=2^{3}$ fixed to demonstrate the convergence of the variance of the level difference as a function of $h_{\ell}$ in Panel (E). This variance converged at rate of $2\eta_{{\rm{s}},d_1}\approx4$, implying that for this example, $\eta_{\rm{w}}=\eta_{{\rm{s}},d_1}\approx2$.

\subsubsection{Consistency and work analysis}
\begin{figure}[ht]
\subfloat[The optimal number of levels $L^{\ast}$ increased at a rate of $\log_2(TOL^{-2/3})$ as a function of the error tolerance $TOL>0$.]{%
		\includegraphics[width=0.45\textwidth]{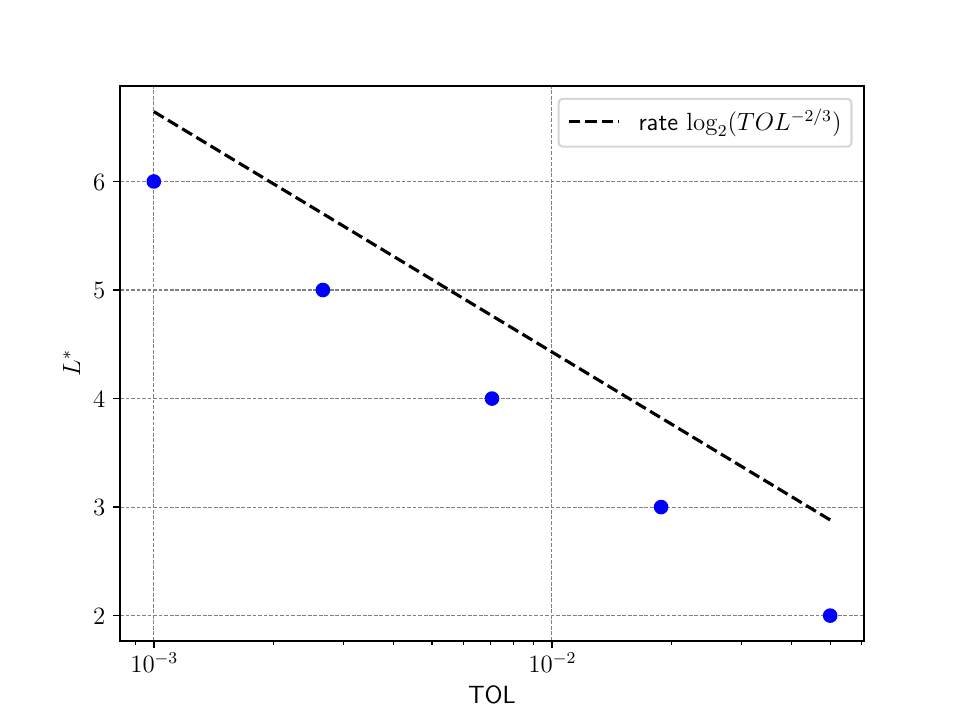}
	}	
 \subfloat[Total computational work of the MLDLQMC estimator as a function of the error tolerance $TOL$. The projected rate of $1.5$ was attained for the considered tolerances.]{%
		\includegraphics[width=0.45\textwidth]{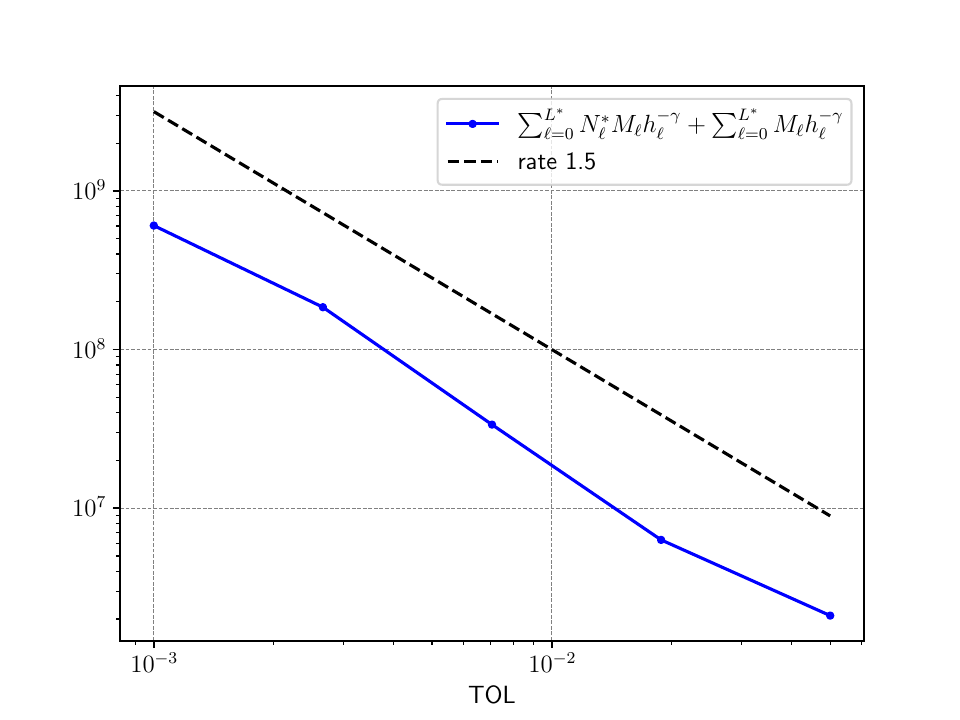}
	}\\
 \subfloat[Error vs.~tolerance $(TOL)$ for 100 runs of the MLDLQMC estimator.]{%
		\includegraphics[width=0.45\textwidth]{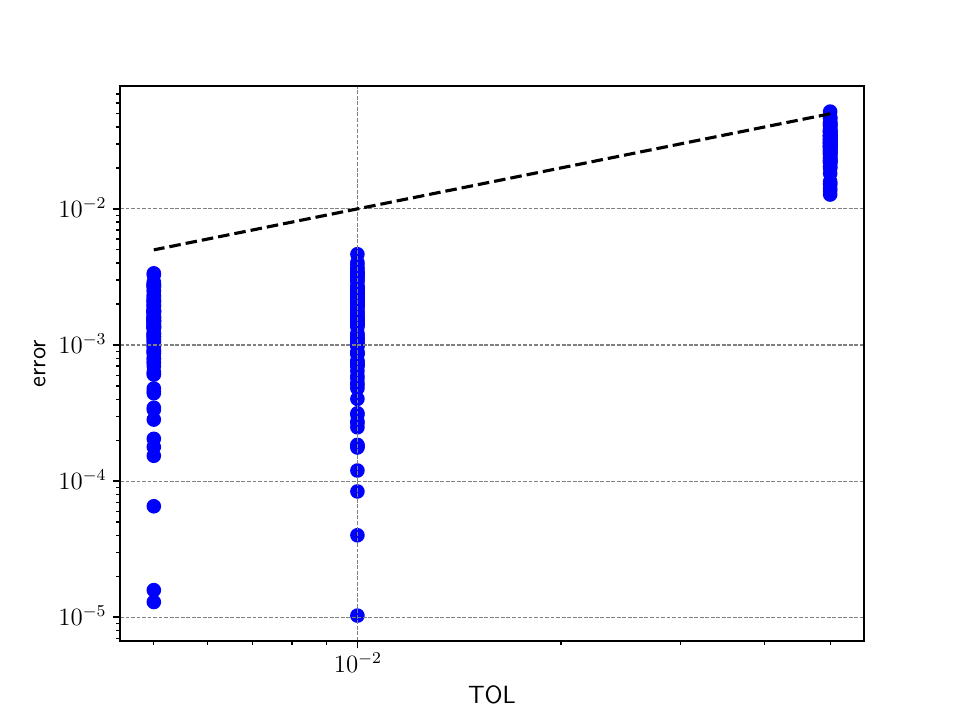}
	}
	
	\caption{Example 3: The MLDLQMC estimator with inexact sampling resulted in error less than the projected tolerance $TOL$ (C) at a cost increasing at the projected rate of 1.5 asymptotically as $TOL\to 0$ (B).}
	\label{fig:poisson.evt.work}
\end{figure}

To balance the contribution of the inner number of samples and mesh discretization parameter, the following choice was selected for this example:
\begin{equation}
    M_{\ell}=M_02^{\ell}, \quad 1\leq \ell\leq L,
\end{equation}
and
\begin{equation}
    h_{\ell}=h_0\left(2^{\frac{3}{4}}\right)^{-\ell}, \quad 1\leq \ell\leq L.
\end{equation}
We determined the optimal final level $L^{\ast}$ of the MLDLQMC estimator with inexact sampling for a given tolerance $TOL>0$. The number of levels is displayed in Figure~\ref{fig:poisson.evt.work} (A) and the rate of increase was $\log_2(TOL^{-2/3})$, as expected. Finally, the accuracy and efficiency of the MLDLQMC estimator for optimal outer samples $N_{\ell}^{\ast}$, where $1\leq \ell\leq L^{\ast}$, was analyzed for various error tolerances $TOL$. In Panel (B), the total work of the MLDLQMC estimator, depending on the optimal $L^{\ast}$ and $N_{\ell}^{\ast}$, where $0\leq \ell\leq L^{\ast}$, is presented as a function of the tolerance $TOL$. The projected rate of increase of approximately 1.5 implied by Corollary~\ref{cor:EIG.case} was observed for this example. In Panel (C), 100 runs of the MLDLQMC estimator with $S=R=1$ randomizations at each level for $TOL=5\cdot 10^{-2}$, $TOL=10^{-2}$, and $TOL=5\cdot 10^{-3}$ are presented, each remaining below the error tolerance projected by the CLT for the given $L^{\ast}$ and $N_{\ell}^{\ast}$, where $0\leq\ell\leq L^{\ast}$ except for one run at the coarsest tolerance. Here, $M_0=2^{7}$ and $h_0\approx 0.42$. An average of 100 runs for $TOL=10^{-3}$ was used as a reference solution. Here, the reference solution was approximately equal to $I\approx 0.779$ and the smallest relative tolerance depicted in Panel (C) is therefore of order $\cl{O}(10^{-2})$.

\FloatBarrier

\subsection{Example 4: Duffing oscillator EIG example with inexact sampling}
\begin{figure}[ht]
    \centering
    \includegraphics[width=0.8\linewidth]{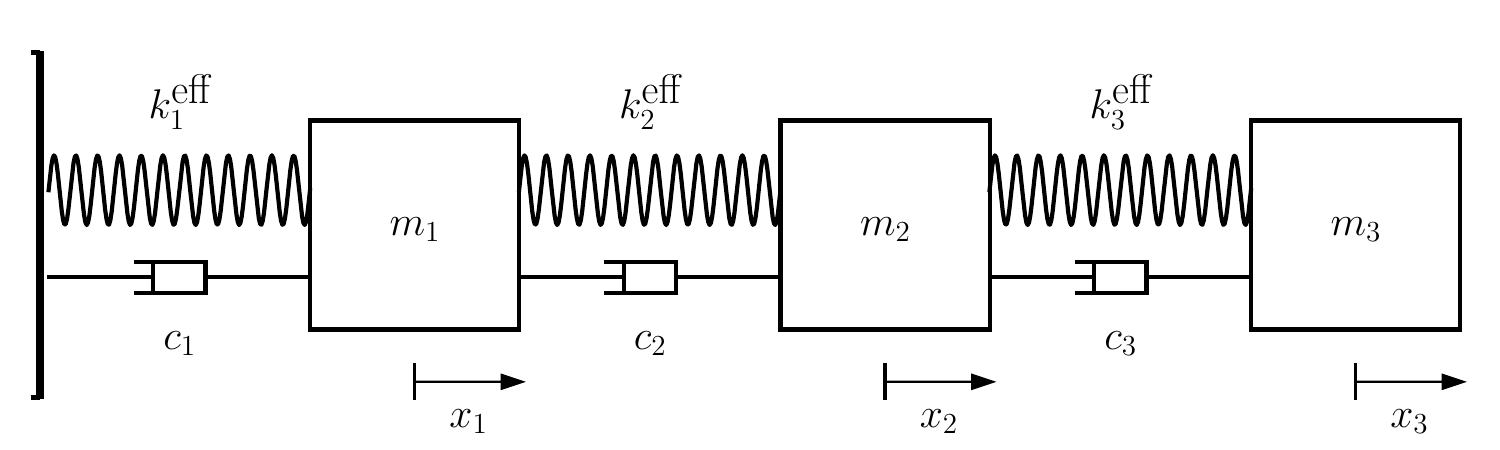}
    \caption{Example 4: diagram representing the spring–mass–damper system. Three masses $(m_1, m_2, m_3)$ are connected by springs with stiffness coefficients $(k_1^\text{eff}, k_2^\text{eff}, k_3^\text{eff})$ and dampers with coefficients $(c_1, c_2, c_3)$.
    Observations of the displacements  $(x_1, x_2, x_3)$ are used to infer the unknown quadratic contributions to the stiffness coefficients.
    }
    \label{fig:duffing.diagram}
\end{figure}

We consider a three-degree-of-freedom spring–mass–damper system in which each spring’s stiffness is augmented by a quadratic term in its relative displacement, yielding a Duffing-type oscillator. In this formulation, the effective stiffness of each spring grows nonlinearly with the square of the displacement difference between its endpoints. In our experiment, we observe the three mass positions at a single time $T$ and use these observations to infer the unknown quadratic stiffness coefficients for each spring. A diagram describing the system is depicted in Figure~\ref{fig:duffing.diagram}.

\subsubsection{Setting and notation}
We consider
\begin{equation}\label{eq:model.ex3}
     \bs{G}(\bs{\theta},\xi)=\bs{u}(\xi),
\end{equation}
where $\bs{u}:[0,T]\mapsto\bb{R}^3$ is the state vector $\bs{u}=\bs{x}=(x_1,x_2,x_3)^{\trans}$ 
and the mechanics of the system are governed by the nonlinear second-order differential equation
\begin{align}\label{eq:Duffing}
    \bs{M} \bs{\ddot{x}} + \bs{C} \bs{\dot{x}} + \bs{K}(\bs{x}) \, \bs{x} {}&= \bs{0}, \quad \text{for all } 0<t\leq T,\\
    \bs{x}(0){}&=\bs{x}_0,\\
    \bs{\dot{x}}(0){}&=\bs{\dot{x}}_0.
\end{align}
For this problem, we assume a diagonal mass matrix
\begin{equation}
      \bs{M} =
  \begin{pmatrix}
    m_1 & 0           & 0 \\
    0            & m_2 & 0\\
    0                           & 0          & m_3
  \end{pmatrix},
\end{equation}
a damping matrix
\begin{equation}
    \bs{C} = \begin{pmatrix}
    c_1 + c_2 & -c_2       & 0 \\
    -c_2      & c_2 + c_3  & -c_3 \\
    0         & -c_3       & c_3
    \end{pmatrix},
\end{equation}
and a state-dependent stiffness matrix
    \begin{equation}
  \bs{K}(\bs{x}) =
  \begin{pmatrix}
    k_1^\text{eff} + k_2^\text{eff} & -k_2^\text{eff}           & 0 \\
    -k_2^\text{eff}              & k_2^\text{eff} + k_3^\text{eff} & -k_3^\text{eff} \\
    0                           & -k_3^\text{eff}           & k_3^\text{eff}
  \end{pmatrix},
\end{equation}
where
    \begin{align}
  k_1^\text{eff} &= k_1 + \theta_1 \, x_1^2, \\
  k_2^\text{eff} &= k_2 + \theta_2 \, (x_2 - x_1)^2, \\
  k_3^\text{eff} &= k_3 + \theta_3 \, (x_3 - x_2)^2.
\end{align}
We moreover assume initial conditions $\bs{x}_0=(100,50,-50)^{\trans}$mm and $\bs{\dot{x}}_0=(0,0,0)^{\trans}$mm/s. The fixed model parameters are presented in Table~\ref{tab:duffing_parameters}.
\begin{table}[ht]
    \centering
    \begin{tabular}{cccc}
        \toprule
         $i$& $m_i$ (kg)  & $k_i$ (N/mm) & $c_i$ (N\,s/mm) \\ \midrule
         1  & 20          & 0.5       & 0.01  \\
         2  & 50          & 0.75       & 0.007 \\
         3  & 30          & 1       & 0.005 \\ \bottomrule
    \end{tabular}
    \caption{Example 4: Duffing oscillator fixed model parameters.}
    \label{tab:duffing_parameters}
\end{table}
The prior distribution of the nonlinear stiffness coefficients is modeled as $\bs{\theta}=(\theta_1,\theta_2,\theta_3)\stackrel{iid}{\sim}\cl{U}([6\cdot10^{-5},7\cdot10^{-5}]^{d_{\theta}})$N/mm$^3$, constituting the parameters of interest; hence $d_{\theta}=3$. We only consider one observation at final time $\xi=T=5$s; thus, the design $\xi$ is in dimension $d_{\xi}=1$. The system~\eqref{eq:Duffing} describes three masses connected by springs and dampers, where $x_i$ is the displacement of the $i$th mass for $1\leq i\leq 3$.  Moreover, the data observations are
\begin{equation}
    \bs{Y}=\bs{G}(\bs{\theta},\bs{\xi})+\bs{\varepsilon},
\end{equation}
where $\bs{\varepsilon}\stackrel{iid}{\sim}\cl{N}(\bs{0},\sigma^2\bs{I}_{d_y\times d_y})$, for $\sigma^2=10^{2}$mm$^2$ and $\bs{I}_{d_y\times d_y}$ is the identity matrix in dimension $d_y=3$, yielding $d_1=d_{\theta}+d_y=6$ and $d_2=d_{\theta}=3$ in the nested integration setting. We consider a discretization
\begin{equation}
     \bs{G}_h(\bs{\theta},\xi)=\bs{u}_h(\xi),
\end{equation}
where $\bs{u}_h$ is obtained via the Runge--Kutta 4 method with weak error convergence rate $\eta_{\rm{w}}=4$ and work rate $\gamma=1$; thus $\gamma/\eta_{\rm{w}}=0.25$. The following choice was selected to balance the contribution of the inner number of samples $M$ and the time stepping parameter $h$:
\begin{equation}
    M_{\ell}=M_02^{\ell}, \quad 1\leq \ell\leq L,
\end{equation}
and
\begin{equation}
    h_{\ell}=h_0\left(2^{\frac{3}{8}}\right)^{-\ell}, \quad 1\leq \ell\leq L.
\end{equation}

\subsubsection{Pilot runs}
\begin{figure}[ht]
	\subfloat[The observed weak convergence rate of the discretization bias of the EIG estimation was $\eta_{\rm{w}}\approx4$.]{%
		\includegraphics[width=0.45\textwidth]{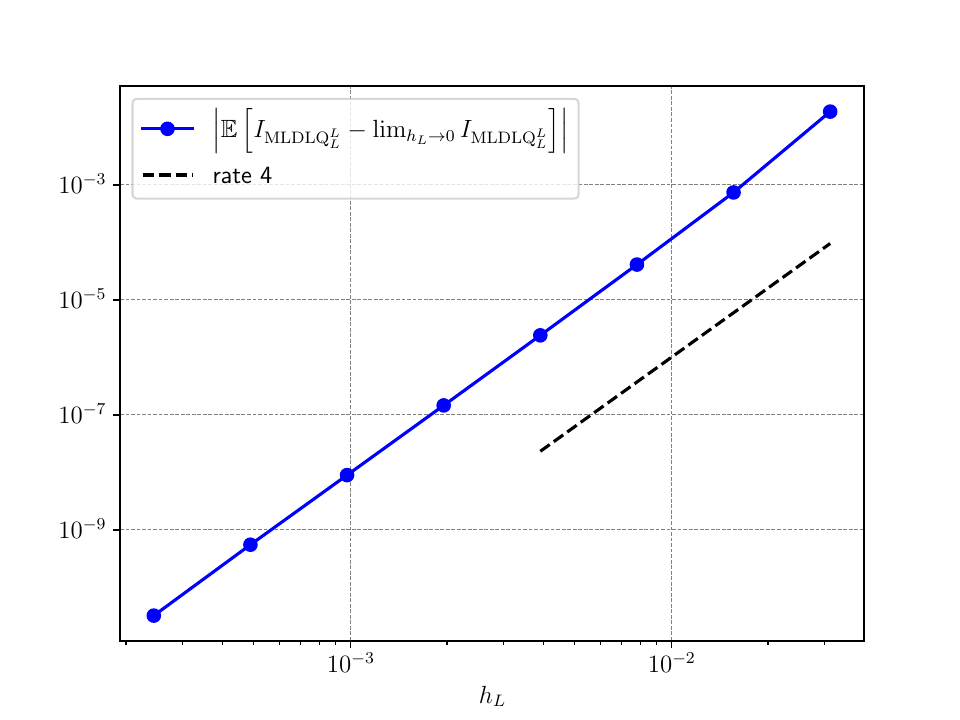}
	}
	\subfloat[The observed convergence rate of approximately 3 of the variance of the inner samples at the final level $L$.]{%
		\includegraphics[width=0.45\textwidth]{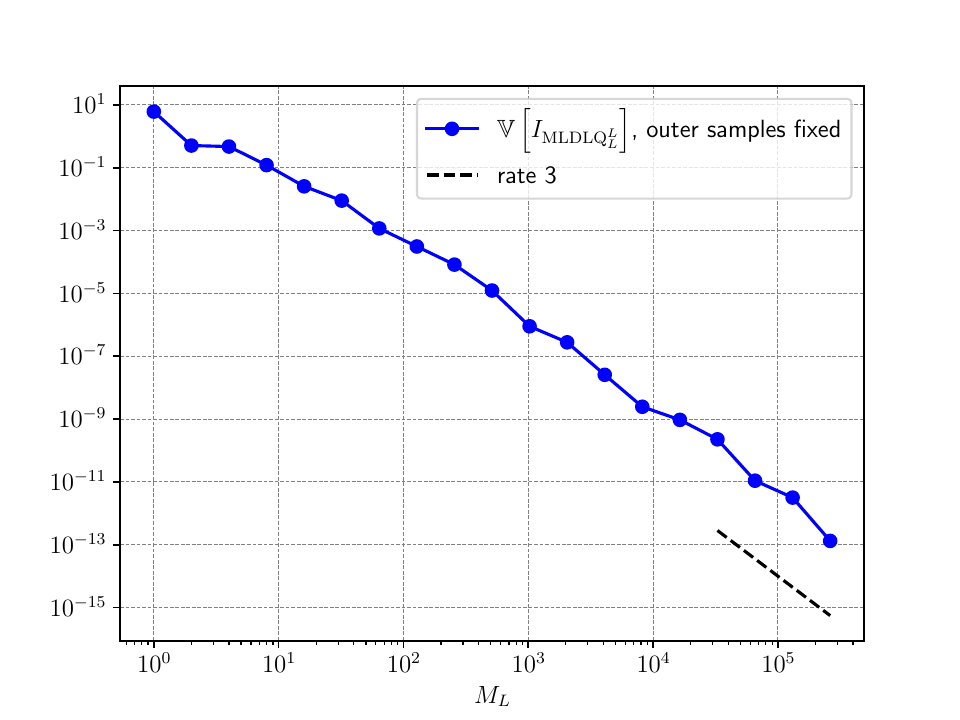}
	}\\
 \subfloat[The observed convergence rate of the variance of approximately 2, with a reduced multiplicative factor for the difference estimator.]{%
		\includegraphics[width=0.45\textwidth]{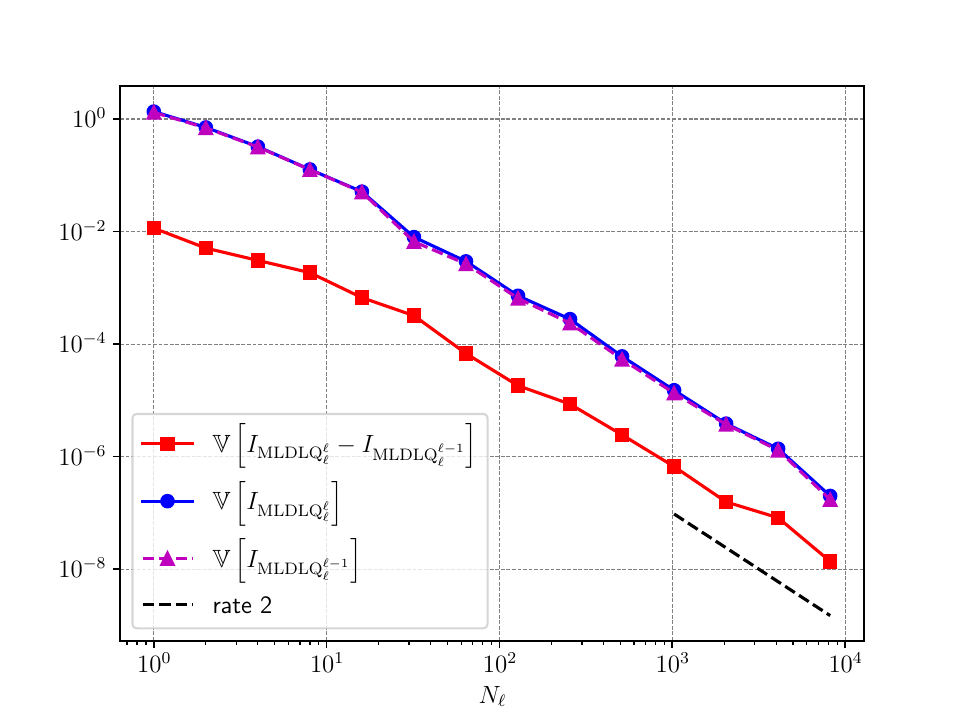}
	}
 \subfloat[Convergence of the level difference as a function of $M_{\ell}$ for $N_{\ell}$ and $h_{\ell}$ fixed at a rate of approximately 2.]{%
		\includegraphics[width=0.45\textwidth]{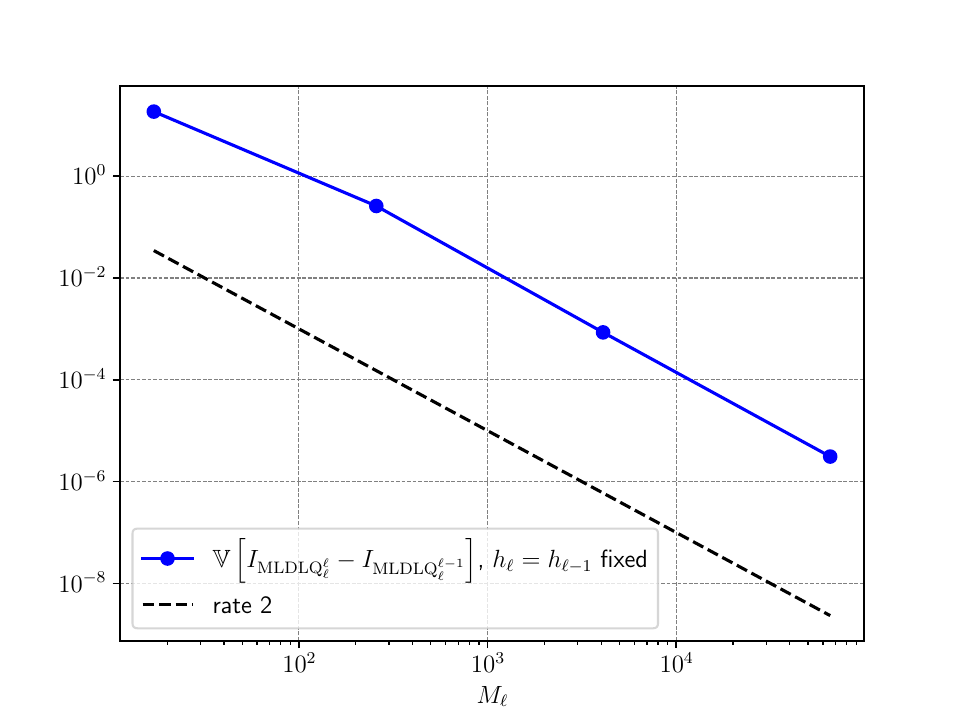}
	}\\
	\subfloat[Convergence of the level difference at a rate of $2\eta_{{\rm{s}},d_1}\approx8$ as a function of $h_{\ell}$ for fixed values of $N_{\ell}$ and $M_{\ell}$.]{%
		\includegraphics[width=0.45\textwidth]{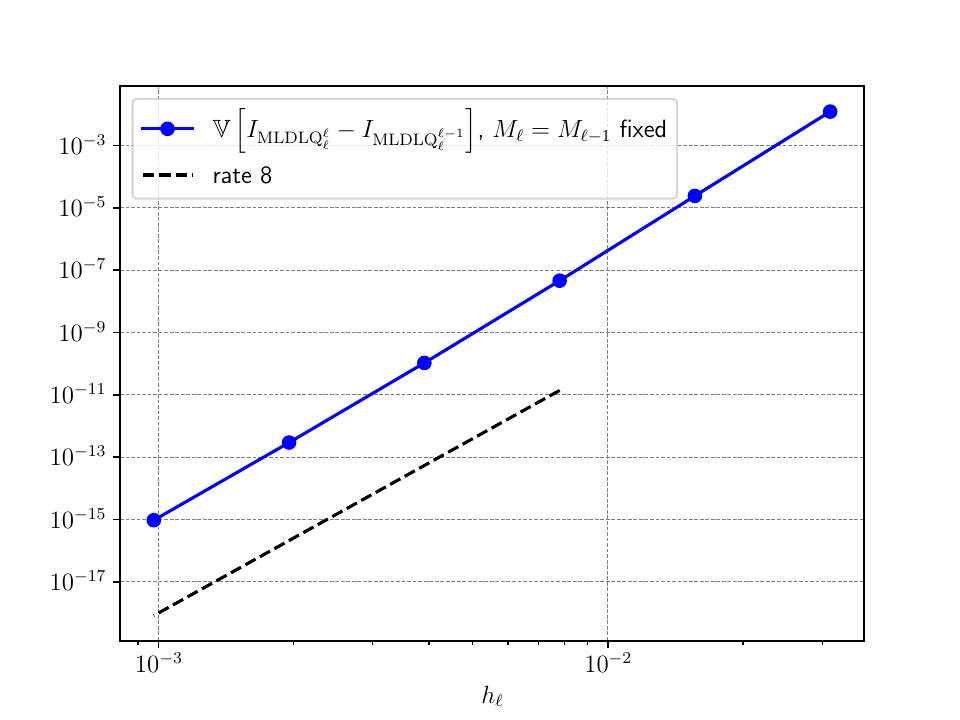}
	}
 \caption{Example 4: Pilot runs for the MLDLQMC estimator with inexact sampling.}
	\label{fig:duffing.pilot}
\end{figure}

 We verified that the bias induced by the Runge--Kutta 4 method converged at a rate of $\eta_{\rm{w}}\approx 4$ by fixing the inner number of samples $M_L$ at the final level $L$ as presented in Figure~\ref{fig:duffing.pilot} (A). For this estimation, $N_L=1$, $M_L=2^{3}$ samples, and $S=100$, $R=1$ randomizations were used, and the bias was estimated using $h_L=2^{-12}$ as a reference solution. Next, we confirmed that the variance contribution from the inner samples at the final level $L$ converged at a rate of 3 (Panel (B)). This estimation was conducted using $N_L=1$ outer sample, $h_L=2^{-8}$, and $S=2^4$, $R=1$ randomizations. The convergence of the variance of the level differences as a function of $N_{\ell}$ at the expected rate of 2 is illustrated in Panel (C), where $M_{\ell-1}=2^{13}$, $M_{\ell}=2^{14}$  samples, $h_{\ell-1}=2^{-4}$, $h_{\ell}=2^{-5}$, and $S=2^{7}$ randomizations were applied. The variances of the single-level estimators converged at the same rate with a larger multiplicative constant. To demonstrate the convergence of the variance of the level difference as a function of $M_{\ell}$, we kept $N_{\ell}=1$ and $h_{\ell}=h_{\ell-1}=2^{-4}$ fixed. This variance converged at rate 2, as seen in Panel (D) for $S=2^{4}$ randomizations. Similarly, we kept $M_{\ell}=M_{\ell-1}=2^{8}$ fixed to demonstrate the convergence of the variance of the level difference as a function of $h_{\ell}$ in Panel (E). This variance converged at rate of $2\eta_{{\rm{s}},d_1}\approx8$, implying that for this example, $\eta_{\rm{w}}=\eta_{{\rm{s}},d_1}\approx4$, as expected.

\subsubsection{Optimal work analysis}
\begin{figure}[ht]
\subfloat[The optimal number of levels $L^{\ast}$ increased at a rate of $\log_2(TOL^{-2/3})$ as a function of the error tolerance $TOL>0$.]{%
		\includegraphics[width=0.45\textwidth]{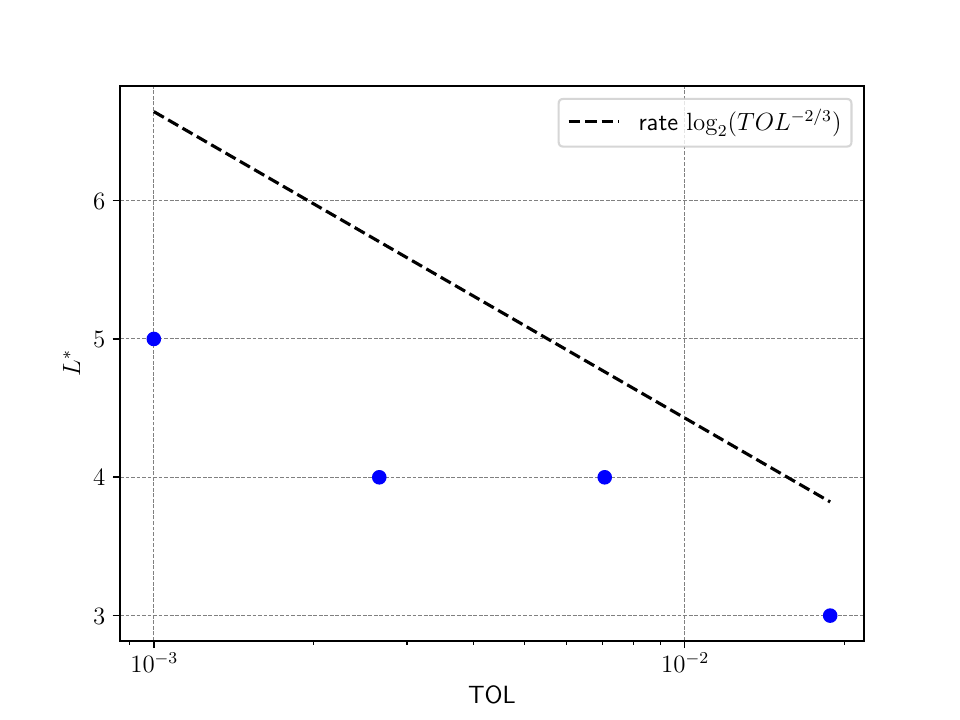}
	}	
 \subfloat[Total computational work of the MLDLQMC estimator as a function of the error tolerance $TOL$.]{%
		\includegraphics[width=0.45\textwidth]{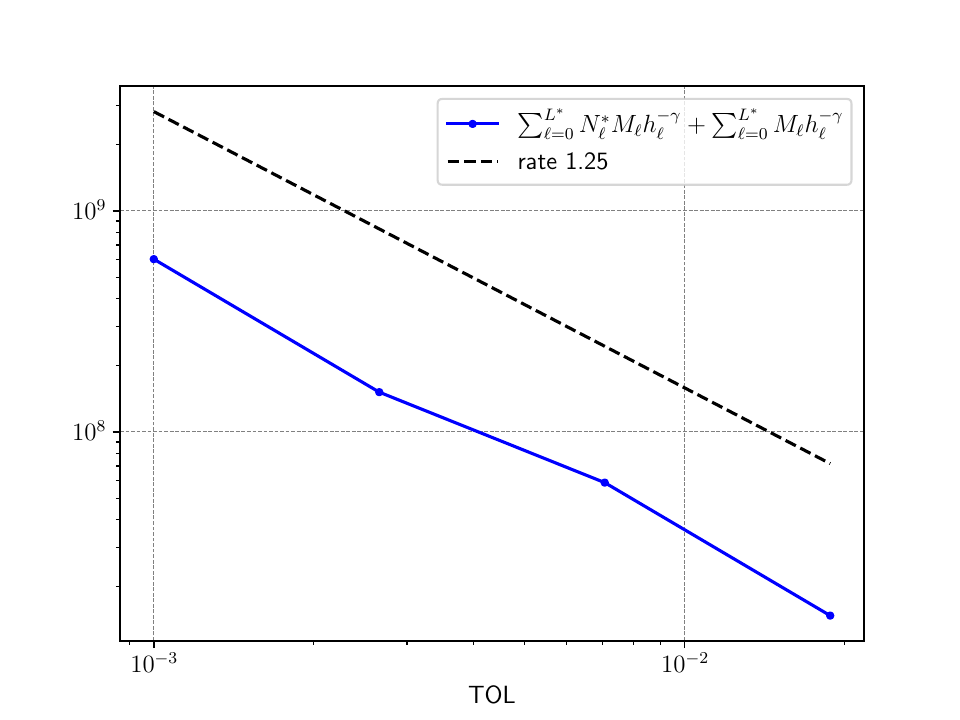}
	}\\
     \subfloat[Error vs.~tolerance $(TOL)$ for 100 runs of the MLDLQMC estimator.]{%
		\includegraphics[width=0.45\textwidth]{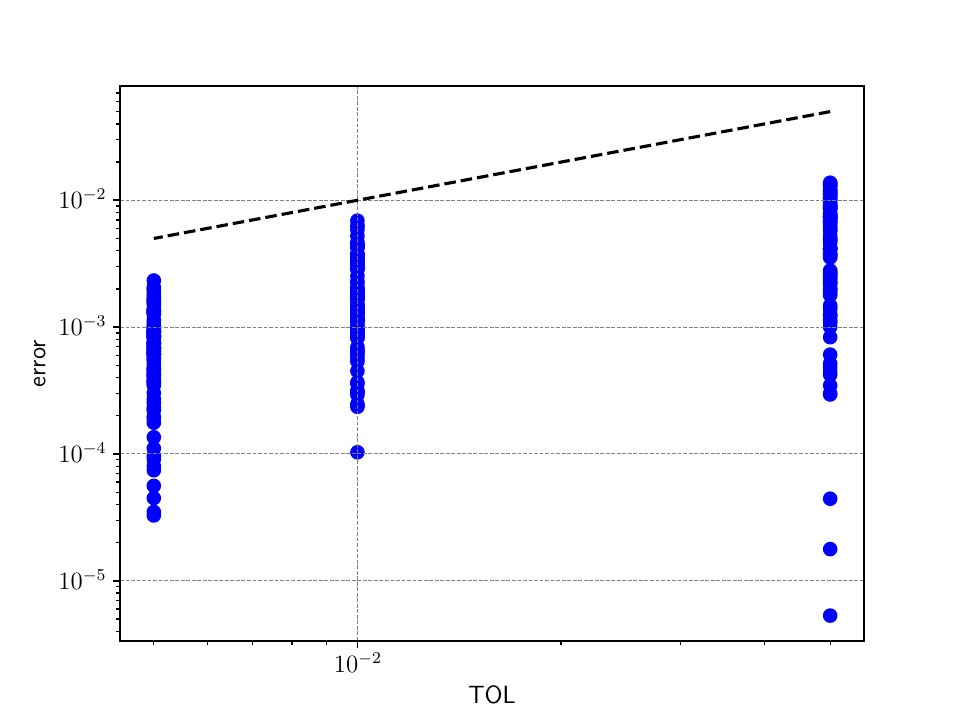}
	}
	\caption{Example 4: The MLDLQMC estimator with inexact sampling resulted in a cost increasing at the projected rate of 1.25 asymptotically as $TOL\to 0$ (B).}
	\label{fig:duffing.evt.work}
\end{figure}

The optimal final level $L^{\ast}$ of the MLDLQMC estimator with inexact sampling for a given tolerance $TOL>0$ is displayed in Figure~\ref{fig:duffing.evt.work} (A), where the rate of increase was $\log_2(TOL^{-2/3})$. In Panel (B), the total work of the MLDLQMC estimator is presented as a function of the tolerance $TOL$, where the projected rate of increase of approximately 1.25 was observed. We used 100 runs of the MLDLQMC estimator with $S=R=1$ randomizations at each level for $TOL=5\cdot 10^{-2}$, $TOL=10^{-2}$, and $TOL=5\cdot 10^{-3}$ and observed that each of them remained below the error tolerance (Panel (C)), where an average of 100 runs for $TOL=10^{-3}$ was used as a reference solution. Here, we used $M_0=2^{7}$ and $h_0\approx 0.08$. The reference solution was approximately equal to $I\approx 0.638$ and the smallest relative tolerance depicted in Panel (C) is therefore of order $\cl{O}(10^{-2})$.

\section{Conclusion}
This work proposes the multilevel double-loop randomized quasi-Monte Carlo estimator for nested integrals with inexact sampling. This estimator combines the multilevel strategy with deterministic and randomized quasi-Monte Carlo (QMC) methods to achieve high approximation accuracy when the inner integrand of a nested integration problem must itself be approximated. The bias and statistical errors are bounded asymptotically, yielding an estimate of the total number of inner and outer samples and the discretization parameter required to achieve a specified error tolerance. Depending on the regularity of the integrands, several different settings are analyzed. In the optimal setting, where both the outer and inner integrands have bounded generalized Hardy--Krause variation of order one, the MLDLQMC estimator achieves a total complexity of almost $\cl{O}(TOL^{-8/9})$ if no additional discretization is necessary, that is, the exact sampling case. A simple example with polynomial integrands demonstrates this optimal complexity numerically. The estimator is then applied to expected information gain (EIG) estimation, where singularities render standard QMC error bounds via the Koksma--Hlawka inequality pointless. A truncation scheme of the Gaussian observation noise in the experiment model with error bounds introduced by Owen~\cite{Owe06} and extended by He et al.~\cite{He23} provide error bounds for this setting, yielding a total complexity of almost $\cl{O}(TOL^{-1-\gamma/\eta_{\rm{w}}})$, where $\gamma$ and $\eta_{\rm{w}}$ relate to the cost and weak convergence rate of the discretization scheme of the inner integrand, respectively. The efficiency of the method for EIG applications is demonstrated on three numerical examples.

\section{Statements and Declarations}
\subsection{Conflict of interest}
The authors have no conflicts to disclose.

\subsection{Data availability}
The data that support the findings of this study are available from the corresponding author upon reasonable request.



\appendix

\section{Auxiliary results for the error bounds to estimate the EIG}\label{app:vcond}
\begin{lem}[Bound on the Hardy--Krause variation of $g_h(\bs{y},\cdot)$]\label{lem:VHK}
     Let $g_h(\bs{y},\cdot)$ be as in \eqref{eq:g.h}. Then, given Assumption~\ref{asu:Lipschitz}, there exists $0<b<\infty$ such that
    \begin{equation}
        V_{\rm{HK}}\left(\left(\prod_{j\in u}\frac{\partial}{\partial y_j}\right)g_h(\bs{y},\cdot)\right)\leq b\prod_{i=1}^{d_1}\min (y_i,1-y_i)^{-A_i-\mathds{1}_{\{i\in u\}}}
    \end{equation}
    for all $\bs{y}\in[0,1]^{d_1}$, all $u\subseteq \{1,\ldots,d_1\}$, and any $A_i>0$, where $1\leq i\leq d_1$.
\end{lem}
\begin{proof}
    For the case where $d_1=2$, $d_2=1$, and
    \begin{equation}
        \bs{\Sigma}_{\bs{\varepsilon}}\coloneqq \sigma^2,
    \end{equation}
    for $\sigma>0$, it follows that
    \begin{equation}\label{eq:g.h.1d}
        g_h(\bs{y},x)= \exp\left({-\frac{1}{2\sigma^2}(\tilde{G}_h(x))^2+\frac{1}{\sigma^2}\tilde{G}_h(y_1)\tilde{G}_h(x)+\frac{1}{\sigma} \tilde{G}_h(x)\Phi^{-1}(y_2)}\right).
    \end{equation}
As a first step, we demonstrate that there exists $0<b<\infty$ such that
\begin{equation}
    V_{\rm{HK}}\left(g_h(\bs{y},\cdot)\right)\leq b\min(y_2, 1-y_2)^{-A_2},
\end{equation}
for any $A_2>0$. The Hardy--Krause variation is bounded as
\begin{align}
    V_{\rm{HK}}\left(g_h(\bs{y},\cdot)\right){}&=\int_{[0,1]}\left|\frac{\partial}{\partial x} \exp\left({-\frac{1}{2\sigma^2}(\tilde{G}_h(x))^2+\frac{1}{\sigma^2}\tilde{G}_h(y_1)\tilde{G}_h(x)+\frac{1}{\sigma} \tilde{G}_h(x)\Phi^{-1}(y_2)}\right)\right| \di{}x,\nonumber\\
    {}&= \int_{[0,1]}\left|\frac{\di{}}{\di{}x}\tilde{G}_h(x)\left(-\frac{1}{\sigma^2}\tilde{G}_h(x)+\frac{1}{\sigma^2}\tilde{G}_h(y_1)+\frac{1}{\sigma}\Phi^{-1}(y_2)\right)\right.\nonumber\\
    {}&\left.\quad\quad \times\exp\left({-\frac{1}{2\sigma^2}(\tilde{G}_h(x))^2+\frac{1}{\sigma^2}\tilde{G}_h(y_1)\tilde{G}_h(x)+\frac{1}{\sigma} \tilde{G}_h(x)\Phi^{-1}(y_2)}\right)\right| \di{}x,\nonumber\\
    {}&\leq \int_{[0,1]}\left|\frac{\di{}}{\di{}x}\tilde{G}_h(x)\right|\left(\frac{1}{\sigma^2}\left|\tilde{G}_h(x)\right|+\frac{1}{\sigma^2}\left|\tilde{G}_h(y_1)\right|+\frac{1}{\sigma}\left|\Phi^{-1}(y_2)\right|\right)\nonumber\\
    {}&\quad\quad \times\exp\left({-\frac{1}{2\sigma^2}(\tilde{G}_h(x))^2+\frac{1}{\sigma^2}|\tilde{G}_h(y_1)||\tilde{G}_h(x)|+\frac{1}{\sigma} |\tilde{G}_h(x)||\Phi^{-1}(y_2)|}\right) \di{}x.
\end{align}
For the case in which $y_2\to0^{+}$, the approximation in \cite{Owe06, Pat96} is employed:
\begin{align}\label{eq:CDF.approx.0}
    \Phi^{-1}(y_2){}&=-\sqrt{2}\sqrt{\log(y_2^{-1})}+o(1),
\end{align}
for which it holds that
\begin{equation}
    \left|\Phi^{-1}(y_2)\right|\leq \sqrt{2}\sqrt{\log(y_2^{-1})}+P,
\end{equation}
for any $P>0$ and all $y_2\leq y_{2,0}$ for some $y_{2,0}\in(0,1)$ as $y_2\to 0^{+}$. Moreover, it follows that
\begin{align}\label{eq:sqrt.A}
    \sqrt{2}\sqrt{\log(y_2^{-1})}{}&\leq \tilde{A}_2\log(y_2^{-1}),\nonumber\\
    {}&= \log(y_2^{-\tilde{A}_2})
\end{align}
for any $\tilde{A}_2>0$ as $y_2\to 0^{+}$, yielding, together with Assumption~\ref{asu:Lipschitz}, that
\begin{align}\label{eq:V.HK.bound}
    {}& \int_{[0,1]}\left|\frac{\di{}}{\di{}x}\tilde{G}_h(x)\right|\left(\frac{1}{\sigma^2}\left|\tilde{G}_h(x)\right|+\frac{1}{\sigma^2}\left|\tilde{G}_h(y_1)\right|+\frac{1}{\sigma}\left|\Phi^{-1}(y_2)\right|\right)\nonumber\\
    {}&\quad\quad \times\exp\left({-\frac{1}{2\sigma^2}(\tilde{G}_h(x))^2+\frac{1}{\sigma^2}|\tilde{G}_h(y_1)||\tilde{G}_h(x)|+\frac{1}{\sigma} |\tilde{G}_h(x)||\Phi^{-1}(y_2)|}\right) \di{}x,\nonumber\\
    {}&\quad\leq \int_{[0,1]}q\left(\frac{2q}{\sigma^2}+\frac{1}{\sigma}\log(y_2^{-\tilde{A}_2})+\frac{1}{\sigma}P\right)\exp\left({\frac{1}{2\sigma^2}q^2+\frac{1}{\sigma} q\log(y_2^{-\tilde{A}_2})+\frac{1}{\sigma}qP}\right) \di{}x,\nonumber\\
    {}&\quad\leq q\left(\frac{2q}{\sigma^2}+\frac{1}{\sigma}\log(y_2^{-\tilde{A}_2})+\frac{1}{\sigma}P\right)y_2^{-\frac{q\tilde{A}_2}{\sigma}}\exp\left({\frac{1}{2\sigma^2}q^2+\frac{1}{\sigma}qP}\right),\nonumber\\
    {}&\quad\leq by_2^{-A_2},
\end{align}
for any $A_2>0$ and an appropriately chosen constant $0<b<\infty$ independent of $y_2$. For the case in which $y_2\to1^{-}$, the approximation in \cite{Owe06, Pat96} is employed:
\begin{align}\label{eq:CDF.approx.1}
    \Phi^{-1}(y_2){}&=\sqrt{2}\sqrt{\log((1-y_2)^{-1})}+o(1),
\end{align}
for which it holds that
\begin{equation}
    \left|\Phi^{-1}(y_2)\right|\leq \sqrt{2}\sqrt{\log((1-y_2)^{-1})}+Q,
\end{equation}
for any $Q>0$ and all $y_2\geq y_{2,1}$ for some $y_{2,1}\in(0,1)$, where $y_{2,1}>y_{2,0}$, as $y_2\to 1^{-}$. It follows immediately that the same bound as in~\eqref{eq:V.HK.bound} applies for any $A_2>0$ and possibly different $0<b<\infty$. The Hardy--Krause variation $V_{\rm{HK}}\left(g_h(\bs{y},\cdot)\right)$ is trivially bounded independently of $y_2$ in between the values $0<y_{2,0}<y_{2,1}<1$, again for an appropriately chosen constant $0<b<\infty$. Thus, we arrive at the following bound:
\begin{equation}
    V_{\rm{HK}}\left(g_h(\bs{y},\cdot)\right)\leq b\min(y_2, 1-y_2)^{-A_2},
\end{equation}
for any $A_2>0$.
As a second step, we demonstrate that there exists $0<b<\infty$ such that
\begin{equation}
    V_{\rm{HK}}\left(\frac{\partial}{\partial y_2}g_h(\bs{y},\cdot)\right)
    \leq b \min(y_2,1-y_2)^{-A_2-1}
\end{equation}
for any $A_2>0$. For this bound, we utilize the following result:
    \begin{align}\label{eq:Phi.prime}
    \frac{\di{}}{\di{} y_2}\Phi^{-1}(y_2){}&=\left(\Phi^{\prime}(\Phi^{-1}(y_2))\right)^{-1},\nonumber\\
{}&=\sqrt{2\pi}e^{\frac{1}{2}(\Phi^{-1}(y_2))^2},\nonumber\\
    {}&\leq\sqrt{2\pi}e^{\left(\log(y_2^{-1})+2P\sqrt{\log(y_2^{-1})}+P^2\right)},\nonumber\\
    {}&\leq\sqrt{2\pi}e^{\left(\log(y_2^{-1})+\tilde{A}_2\log(y_2^{-1})+P^2\right)},\nonumber\\
    {}&=\sqrt{2\pi}e^{P^2}y_2^{-\tilde{A}_2-1},
\end{align}
as $y_2\to0^{+}$ to bound
\begin{align}\label{eq:Vhk.bound.y2}
    V_{\rm{HK}}{}&\left(\frac{\partial}{\partial y_2}g_h(\bs{y},\cdot)\right)\nonumber\\ 
    {}&=\int_{[0,1]}\left|\frac{\partial}{\partial x}\frac{\partial}{\partial y_2}\exp\left({-\frac{1}{2\sigma^2}(\tilde{G}_h(x))^2+\frac{1}{\sigma^2}\tilde{G}_h(y_1)\tilde{G}_h(x)+\frac{1}{\sigma} \tilde{G}_h(x)\Phi^{-1}(y_2)}\right)\right| \di{}x,\nonumber\\
    {}&\leq \int_{[0,1]}\left|\frac{\di{}}{\di{}y_2}\Phi^{-1}(y_2)\frac{1}{\sigma}\frac{\di{}}{\di{}x}\tilde{G}_h(x)\exp\left({-\frac{1}{2\sigma^2}(\tilde{G}_h(x))^2+\frac{1}{\sigma^2}\tilde{G}_h(y_1)\tilde{G}_h(x)+\frac{1}{\sigma} \tilde{G}_h(x)\Phi^{-1}(y_2)}\right)\right|\nonumber\\
    {}&\quad+ \left|\frac{\di{}}{\di{}y_2}\Phi^{-1}(y_2)\frac{1}{\sigma}\tilde{G}_h(x)\frac{\di{}}{\di{}x}\tilde{G}_h(x)\left(-\frac{1}{\sigma^2}\tilde{G}_h(x)+\frac{1}{\sigma^2}\tilde{G}_h(y_1)+\frac{1}{\sigma}\Phi^{-1}(y_2)\right)\right.\nonumber\\
    {}&\left.\quad\quad \times\exp\left({-\frac{1}{2\sigma^2}(\tilde{G}_h(x))^2+\frac{1}{\sigma^2}\tilde{G}_h(y_1)\tilde{G}_h(x)+\frac{1}{\sigma} \tilde{G}_h(x)\Phi^{-1}(y_2)}\right)\right| \di{}x,\nonumber\\
    {}&\leq \tilde{b} y_2^{-\tilde{A}_2-1}
\end{align}
as $y_2\to0^{+}$ for some $0<\tilde{b}<\infty$ and any $\tilde{A}_2>0$. From an analogous argument as $y_2\to0^{-}$, it follows that there exists $0<b<\infty$ such that
\begin{equation}
    V_{\rm{HK}}\left(\frac{\partial}{\partial y_2}g_h(\bs{y},\cdot)\right)
    \leq b \min(y_2,1-y_2)^{-A_2-1}
\end{equation}
for any $A_2>0$. The partial derivative with respect to $y_1$ follows trivially. The general case follows from Fa\`{a} di Bruno's formula.
\end{proof}

\begin{cor}[Bound on $\bar{g}_h(\bs{y})$]\label{cor:bar.g}
Let $\bar{g}_h$ be as in~\eqref{eq:g.h.bar}. Then, given Assumption~\ref{asu:Lipschitz}, there exists $0<b<\infty$ such that
\begin{equation}
\left|\left(\prod_{j\in u}\frac{\partial}{\partial y_j}\right)\bar{g}_h(\bs{y})\right|\leq b\prod_{i=1}^{d_1}\min (y_i,1-y_i)^{-A_i-\mathds{1}_{\{i\in u\}}},
    \end{equation}
    for all $\bs{y}\in[0,1]^{d_1}$, all $u\subseteq \{1,\ldots,d_1\}$, and any $A_i>0$, where $1\leq i\leq d_1$.
\end{cor}
\begin{proof}
    For the case where $d_1=2$, $d_2=1$, and $g_h$ as in~\eqref{eq:g.h.1d}, the result follows from the Leibniz integral rule and a similar argument as in~\eqref{eq:Vhk.bound.y2}. The general case follows from Fa\`{a} di Bruno's formula as well.
\end{proof}
\begin{cor}[Bound on the nested integrand]\label{cor:bound.nested.integrand}
    Given Assumption\ref{asu:Lipschitz}, for $f\equiv\log$ and $g_h$ as in \eqref{eq:g.h}, there exists $0<\delta_{K}<1/2$ and $0<b<\infty$ such that
    \begin{equation}
        \left|f^{(2+d_1)}\left(K\bar{g}_h(\bs{y})\right)\right|\leq b\prod_{i=d_2+1}^{d_1}\min(y_i,1-y_i)^{-A_i}
    \end{equation}
    for any $A_i>0$, where $1\leq i\leq d_1$, and all $K\in(\delta_{K},1-\delta_{K})$. Here, $f^{(d)}$ denotes the $d$th derivative of $f$.
\end{cor}
\begin{proof}
    The following holds for $f\equiv\log$:
    \begin{align}
        \left|f^{(2+d_1)}\left(K\bar{g}_h(\bs{y})\right)\right|{}&=(1+d_1)!K^{-2-d_1}(\bar{g}_h(\bs{y}))^{-2-d_1}.
    \end{align}
    
    In Corollary~\ref{cor:bar.g}, $\bar{g}_h$ was bounded for the cases when $\Phi^{-1}(\bs{y})\to-\infty$ and $\Phi^{-1}(\bs{y})\to+\infty$; thus, similar bounds also apply to $(\bar{g}_h)^{-1}$.
\end{proof}
The following result is a direct consequence of \cite[Lemma 2]{Bar23}.
\begin{lem}[Inverse inequality for the inner integrand with truncated noise]\label{cor:inverse.g.h}
    Let $TOL>0$ and $c(TOL)$ as in~\eqref{eq:c.TOL.q}. Moreover, let
\begin{equation}\label{eq:inverse.CDF}
F_{c(TOL)}^{-1}(y_2)\coloneqq\Phi^{-1}\Big((2\Phi(c(TOL))-1)y_2+(1-\Phi(c(TOL)))\Big)
\end{equation}
denote the inverse of the standard normal distribution truncated to the interval $[-c(TOL),c(TOL)]$. Then, it holds for
\begin{equation}\label{eq:g.h.truncated}
        g_h^{\rm{tr}}(\bs{y},\bs{x})\coloneqq \exp\left({-\frac{1}{2}\left\lVert\bs{\tilde{G}}_h(\bs{x})\right\rVert_{\bs{\Sigma}_{\bs{\varepsilon}}^{-1}}^2+\left<\bs{\tilde{G}}_h(\bs{y}_1),\bs{\tilde{G}}_h(\bs{x})\right>_{\bs{\Sigma}_{\bs{\varepsilon}}^{-1}}+\left<\bs{\tilde{G}}_h(\bs{x}),F_{c(TOL)}^{-1}(\bs{y}_2)\right>_{\bs{\tilde{\Sigma}}_{\bs{\varepsilon}}^{-\trans}}}\right)
    \end{equation}
    and
    \begin{equation}
        \bar{g}_h^{\rm{tr}}(\bs{y})\coloneqq \int_{[0,1]^{d_2}}g_h^{\rm{tr}}(\bs{y},\bs{x})\di{}\bs{x}
    \end{equation}
    that there exists $0<\tilde{L}<\infty$ such that
    \begin{equation}
    \sup_{\bs{y}\in[0,1]^{d_1}} \left\lvert \frac{V_{\rm{HK}}(g_h^{\rm{tr}}(\bs{y},\cdot))}{ \bar{g}^{\rm{tr}}_h(\bs{y})} \right\rvert < \tilde{L}^{d_2}c(TOL)^{d_2},
\end{equation}
with $V_{\rm{HK}}$ as in \eqref{eq:VHK}.
\end{lem}
\section{Proof of Lemma~\ref{lem:level.0}}\label{app:0}
\begin{proof}[Proof of Lemma~\ref{lem:level.0}]
    We again consider the case where $d_1=2$, $d_2=1$, and $g_h$ as in~\eqref{eq:g.h.1d}. As a first step, we demonstrate that there exists $0<b<\infty$ such that
    \begin{equation}
    \left|f\left(\frac{1}{M}\sum_{m=1}^{M}g_{h}\left(\bs{y},x^{(m)}\right)\right)\right|\leq b\min (y_2,1-y_2)^{-A_2}
\end{equation}
for any $A_2>0$.
    For $f\equiv\log$ and $g_{h}$ as in~\eqref{eq:g.h.1d}, it follows that
    \begin{align}
       {}&\left|f\left(\frac{1}{M}\sum_{m=1}^{M}g_{h}\left(\bs{y},x^{(m)}\right)\right)\right|\nonumber\\
       {}&\quad=\left|\log\left(\frac{1}{M}\sum_{m=1}^{M}\exp\left({-\frac{1}{2\sigma^2}(\tilde{G}_h(x^{(m)}))^2+\frac{1}{\sigma^2}\tilde{G}_h(y_1)\tilde{G}_h(x^{(m)})+\frac{1}{\sigma} \tilde{G}_h(x^{(m)})\Phi^{-1}(y_2)}\right)\right)\right|.
    \end{align}
    For now, we consider the case where $\tilde{G}_h$ is nonnegative for all $x^{(m)}$, where $1\leq m\leq M$, that is,
    \begin{equation}
        \min_{1\leq m\leq M}\tilde{G}_h(x^{(m)})\geq 0.
    \end{equation}
    Then, as $y_2\to0^{+}$, it follows from Assumption~\ref{asu:Lipschitz} and the bound in~\eqref{eq:CDF.approx.0} that there exists $0<\tilde{b}<\infty$ such that
    \begin{align}
        \exp{}&\left({-\frac{1}{2\sigma^2}(\tilde{G}_h(x^{(m)}))^2+\frac{1}{\sigma^2}\tilde{G}_h(y_1)\tilde{G}_h(x^{(m)})+\frac{1}{\sigma} \tilde{G}_h(x^{(m)})\Phi^{-1}(y_2)}\right)\nonumber\\
        {}&\geq\exp\left(-\frac{3}{2\sigma^2}q^2\right)\exp\left(-\frac{1}{\sigma} q\left|\Phi^{-1}(y_2)\right|\right),\nonumber\\
        {}&\geq\tilde{b}y_2^{\tilde{A}_2},\quad 1\leq m\leq M,
    \end{align}
    for any $\tilde{A}_2>0$.
    From this result, it follows that
    \begin{equation}
        \left|\log\left(\tilde{b}y_2^{\tilde{A}_2}\right)\right|\leq by_2^{-A_2}
    \end{equation}
    for some appropriately chosen $0<b<\infty$ and any $A_2>0$ as $y_2\to0^{+}$. Now, we consider the case that $\tilde{G}_h$ is negative for at least one $x^{(m)}$, where $1\leq m\leq M$, that is,
    \begin{equation}
        \min_{1\leq m\leq M} \tilde{G}_h(x^{(m)})<0,
    \end{equation}
    and furthermore, let
    \begin{equation}
        m^{\ast}\coloneqq \argmin_{1\leq m\leq M} \tilde{G}_h(x^{(m)}).
    \end{equation}
    Thus, it follows that there exists $0<\tilde{b},\tilde{c}<\infty$ such that
    \begin{align}
        \frac{1}{M}{}&\sum_{m=1}^{M}\exp\left({-\frac{1}{2\sigma^2}(\tilde{G}_h(x^{(m)}))^2+\frac{1}{\sigma^2}\tilde{G}_h(y_1)\tilde{G}_h(x^{(m)})+\frac{1}{\sigma} \tilde{G}_h(x^{(m)})\Phi^{-1}(y_2)}\right)\nonumber\\
        {}&\leq \tilde{c}\exp\left({-\frac{1}{2\sigma^2}(\tilde{G}_h(x^{(m^{\ast})}))^2+\frac{1}{\sigma^2}\tilde{G}_h(y_1)\tilde{G}_h(x^{(m^{\ast})})+\frac{1}{\sigma} \tilde{G}_h(x^{(m^{\ast})})\Phi^{-1}(y_2)}\right),\nonumber\\
        {}&\leq\tilde{c}\exp\left(\frac{1}{2\sigma^2}q^2\right)\exp\left(\frac{1}{\sigma} q\left|\Phi^{-1}(y_2)\right|\right),\nonumber\\
        {}&\leq \tilde{b}y_2^{-\tilde{A}_2},
    \end{align}
    for any $\tilde{A}_2>0$, and hence,
    \begin{equation}
        \left|\log\left(\tilde{b}y_2^{-\tilde{A}_2}\right)\right|\leq by_2^{-A_2}
    \end{equation}
    for some appropriately chosen $0<b<\infty$ and any $A_2>0$ as $y_2\to0^{+}$. The case where $y_2\to1^{-}$ follows by a similar argument and the bound~\eqref{eq:CDF.approx.1}. Thus, it is demonstrated that
    \begin{equation}
    \left|f\left(\frac{1}{M}\sum_{m=1}^{M}g_{h}\left(\bs{y},x^{(m)}\right)\right)\right|\leq b\min (y_2,1-y_2)^{-A_2}
\end{equation}
for any $A_2>0$.
    As a second step, it follows that
    \begin{align}\label{eq:f.0.y2}
        {}&\left|\frac{\partial}{\partial y_2}f\left(\frac{1}{M}\sum_{m=1}^{M}g_{h}\left(\bs{y},\bs{x}^{(m)}\right)\right)\right|\nonumber\\
        {}&= \left|\frac{\frac{1}{M}\sum_{m=1}^{M}\left(\frac{\di{}}{\di{}y_2}\Phi^{-1}(y_2)\frac{1}{\sigma} \tilde{G}_h(x^{(m)})\right)\exp\left({-\frac{1}{2\sigma^2}(\tilde{G}_h(x^{(m)}))^2+\frac{1}{\sigma^2}\tilde{G}_h(y_1)\tilde{G}_h(x^{(m)})+\frac{1}{\sigma} \tilde{G}_h(x^{(m)})\Phi^{-1}(y_2)}\right)}{\frac{1}{M}\sum_{m=1}^{M}\exp\left({-\frac{1}{2\sigma^2}(\tilde{G}_h(x^{(m)}))^2+\frac{1}{\sigma^2}\tilde{G}_h(y_1)\tilde{G}_h(x^{(m)})+\frac{1}{\sigma} \tilde{G}_h(x^{(m)})\Phi^{-1}(y_2)}\right)}\right|.
    \end{align}
    For the numerator in equation~\eqref{eq:f.0.y2}, there exists $0<\tilde{b}<\infty$ such that
\begin{align}
    {}&\left|\frac{1}{M}\sum_{m=1}^{M}\left(\frac{\di{}}{\di{}y_2}\Phi^{-1}(y_2)\frac{1}{\sigma} \tilde{G}_h(x^{(m)})\right)\exp\left({-\frac{1}{2\sigma^2}(\tilde{G}_h(x^{(m)}))^2+\frac{1}{\sigma^2}\tilde{G}_h(y_1)\tilde{G}_h(x^{(m)})+\frac{1}{\sigma} \tilde{G}_h(x^{(m)})\Phi^{-1}(y_2)}\right)\right|\nonumber\\
    {}&\leq \tilde{b} y_2^{-\tilde{A}_2-1},
\end{align}
for any $\tilde{A}_2>0$ by a similar derivation and the bound~\eqref{eq:Phi.prime}. For the denominator in equation~\eqref{eq:f.0.y2}, it holds that
    \begin{equation}
    \left|\frac{1}{M}\sum_{m=1}^{M}\exp\left({-\frac{1}{2\sigma^2}(\tilde{G}_h(x^{(m)}))^2+\frac{1}{\sigma^2}\tilde{G}_h(y_1)\tilde{G}_h(x^{(m)})+\frac{1}{\sigma} \tilde{G}_h(x^{(m)})\Phi^{-1}(y_2)}\right)\right|\geq \tilde{b} y_2^{-\tilde{A}_2},
\end{equation}
again by a similar derivation. Thus, it follows that there exists $0<b<\infty$ such that
\begin{equation}
    \left|\frac{\partial}{\partial y_2}f\left(\frac{1}{M}\sum_{m=1}^{M}g_{h}\left(\bs{y},\bs{x}^{(m)}\right)\right)\right|\leq b y_2^{-A_2-1}
\end{equation}
for any $A_2>0$. The general case follows similarly by Fa\`{a} di Bruno's formula.
\end{proof}
    
\section{Proof of Lemma~\ref{lem:boundary.growth.taylor}}\label{app:ell}
\begin{proof}[Proof of Lemma~\ref{lem:boundary.growth.taylor}]
We introduce the following notations:
\begin{equation}
    \varphi_{h,M}(\bs{y})=\left(\lambda(\bs{y})\right)^2\mu(\bs{y}),
\end{equation}
where
\begin{equation}
    \lambda(\bs{y})\coloneqq \frac{1}{M}\sum_{m=1}^Mg^{\rm{tr}}_h(\bs{y},\bs{x}^{(m)})-\bar{g}^{\rm{tr}}(\bs{y}),
\end{equation}
and
\begin{equation}
    \mu(\bs{y})\coloneqq \int_0^1f''\left(\bar{g}^{\rm{tr}}(\bs{y})+s \left(\frac{1}{M}\sum_{m=1}^Mg^{\rm{tr}}_h(\bs{y},\bs{x}^{(m)})-\bar{g}^{\rm{tr}}(\bs{y})\right)\right)(1-s)\di{}s.
\end{equation}
For the case $d_1=2$, it must be verified that the following inequalities hold:
\begin{align}
    \left|\varphi_{h,M}(y_1,y_2)\right|{}&\leq \left(C_{\epsilon,d_2}^2M^{-2+2\epsilon}+C_{\eta_{{\rm{s}},d_1}}^2h^{2\eta_{{\rm{s}},d_1}}\right)b\min(y_1,1-y_1)^{-A_1}\min(y_2,1-y_2)^{-A_2},\label{eq:f.h.M.0}\\
        \left|\frac{\partial}{\partial y_1}\varphi_{h,M}(y_1,y_2)\right|{}&\leq \left(C_{\epsilon,d_2}^2M^{-2+2\epsilon}+C_{\eta_{{\rm{s}},d_1}}^2h^{2\eta_{{\rm{s}},d_1}}\right)b\min(y_1,1-y_1)^{-A_1-1}\min(y_2,1-y_2)^{-A_2},\label{eq:f.h.M.1}\\
    \left|\frac{\partial}{\partial y_2}\varphi_{h,M}(y_1,y_2)\right|{}&\leq \left(C_{\epsilon,d_2}^2M^{-2+2\epsilon}+C_{\eta_{{\rm{s}},d_1}}^2h^{2\eta_{{\rm{s}},d_1}}\right)b\min(y_1,1-y_1)^{-A_1}\min(y_2,1-y_2)^{-A_2-1},\label{eq:f.h.M.2}\\
    \left|\frac{\partial^2}{\partial y_1\partial y_2}\varphi_{h,M}(y_1,y_2)\right|{}&\leq \left(C_{\epsilon,d_2}^2M^{-2+2\epsilon}+C_{\eta_{{\rm{s}},d_1}}^2h^{2\eta_{{\rm{s}},d_1}}\right)b\min(y_1,1-y_1)^{-A_1-1}\min(y_2,1-y_2)^{-A_2-1}.\label{eq:f.h.M.3}
\end{align}
We start with the Condition~\eqref{eq:f.h.M.0}, for which it follows that
\begin{align}\label{eq:f.h.M.0.split}
    \left|\varphi_{h,M}(y_1,y_2)\right|{}&=\left|\left(\lambda(y_1,y_2)\right)^2\mu(y_1,y_2)\right|\nonumber\\
    {}&=\left(\lambda(y_1,y_2)\right)^2\left|\mu(y_1,y_2)\right|.
\end{align}
The bound for the first term in \eqref{eq:f.h.M.0.split} is
\begin{align}\label{eq:f.h.M.0.split.M}
    \left(\lambda(\bs{y})\right)^2{}&=\left(\frac{1}{M}\sum_{m=1}^Mg^{\rm{tr}}_h(\bs{y},\bs{x}^{(m)})-\bar{g}^{\rm{tr}}(\bs{y})\right)^2\nonumber\\
    {}&=\left(\frac{1}{M}\sum_{m=1}^Mg^{\rm{tr}}_h(\bs{y},\bs{x}^{(m)})-\bar{g}^{\rm{tr}}_h(\bs{y})+\bar{g}^{\rm{tr}}_h(\bs{y})-\bar{g}^{\rm{tr}}(\bs{y})\right)^2,\nonumber\\
    {}&\leq 2\left(\frac{1}{M}\sum_{m=1}^Mg^{\rm{tr}}_h(\bs{y},\bs{x}^{(m)})-\bar{g}^{\rm{tr}}_h(\bs{y})\right)^2+2\left(\bar{g}^{\rm{tr}}_h(\bs{y})-\bar{g}^{\rm{tr}}(\bs{y})\right)^2.
\end{align}
By the Koksma--Hlawka inequality~\eqref{eq:QMC.bound} and Lemma~\ref{lem:VHK}, there exists $0<\tilde{b}<\infty$ such that the first term in \eqref{eq:f.h.M.0.split.M} is bounded as follows:
\begin{align}
    2\left(\frac{1}{M}\sum_{m=1}^Mg^{\rm{tr}}_h(\bs{y},\bs{x}^{(m)})-\bar{g}^{\rm{tr}}_h(\bs{y})\right)^2 {}&\leq2V_{\rm{HK}}(g^{\rm{tr}}_h(\bs{y},\cdot))^2C_{\epsilon,d_2}^2M^{-2+2\epsilon},\nonumber\\
    {}&\leq2C_{\epsilon,d_2}^2M^{-2+2\epsilon}\tilde{b}^2\min (y_1,1-y_1)^{-2\tilde{A}_1}\min (y_2,1-y_2)^{-2\tilde{A}_2}
\end{align}
 for any $\tilde{A}_1,\tilde{A}_2>0$. By Assumption~\ref{asu:inverse.FEM} and Corollary~\ref{cor:bar.g}, there exists $0<\tilde{b}<\infty$ such that the second term in \eqref{eq:f.h.M.0.split.M} is bounded as follows:
\begin{align}
    2\left(\bar{g}^{\rm{tr}}_h(\bs{y})-\bar{g}^{\rm{tr}}(\bs{y})\right)^2{}&\leq 2(\bar{g}^{\rm{tr}}(\bs{y}))^2C_{\eta_{{\rm{s}}}}^2h^{2\eta_{\rm{s}}},\nonumber\\
    {}&\leq2C_{\eta_{{\rm{s}}}}^2h^{2\eta_{\rm{s}}}\tilde{b}^2\min (y_1,1-y_1)^{-2\tilde{A}_1}\min (y_2,1-y_2)^{-2\tilde{A}_2}
\end{align}
 for any $\tilde{A}_1,\tilde{A}_2>0$. Thus, the bound for the first term in~\eqref{eq:f.h.M.0.split} follows as
 \begin{equation}\label{eq:bound.lambda}
     \left(\lambda(\bs{y})\right)^2\leq \left(C_{\epsilon,d_2}^2M^{-2+2\epsilon}+C_{\eta_{{\rm{s}}}}^2h^{2\eta_{\rm{s}}}\right)2\tilde{b}^2\min (y_1,1-y_1)^{-2\tilde{A}_1}\min (y_2,1-y_2)^{-2\tilde{A}_2}.
 \end{equation}
  To bound the second term in \eqref{eq:f.h.M.0.split}, it follows directly from the definition~\eqref{eq:g.h} that
\begin{equation}
    \bar{g}^{\rm{tr}}(\bs{y})+s \left(\frac{1}{M}\sum_{m=1}^Mg^{\rm{tr}}_h(\bs{y},\bs{x}^{(m)})-\bar{g}^{\rm{tr}}(\bs{y})\right)> 0.
\end{equation}
For $f\equiv\log$, it is thus verified that $|f''|$ is monotonically decreasing for positive inputs. It follows that
\begin{align}\label{eq:monotonicity.bound}
\left|\mu(\bs{y})\right|={}&\Bigg|\int_0^1f''\left(\bar{g}^{\rm{tr}}(\bs{y})+s \left(\frac{1}{M}\sum_{m=1}^Mg^{\rm{tr}}_h(\bs{y},\bs{x}^{(m)})-\bar{g}^{\rm{tr}}(\bs{y})\right)\right)(1-s)\di{}s\Bigg|\nonumber\\
\leq{}&\int_0^1\Bigg|f''\left(\bar{g}^{\rm{tr}}(\bs{y})+s \left(\frac{1}{M}\sum_{m=1}^Mg^{\rm{tr}}_h(\bs{y},\bs{x}^{(m)})-\bar{g}^{\rm{tr}}(\bs{y})\right)\right)\Bigg|(1-s)\di{}s,\nonumber\\
\leq{}&\int_0^1\Bigg|f''\left(\bar{g}^{\rm{tr}}(\bs{y})-s\Bigg| \frac{1}{M}\sum_{m=1}^Mg^{\rm{tr}}_h(\bs{y},\bs{x}^{(m)})-\bar{g}^{\rm{tr}}(\bs{y})\Bigg|\right)\Bigg|(1-s)\di{}s,\nonumber\\
\leq{}&\Bigg|f''\left(\bar{g}^{\rm{tr}}(\bs{y})-\Bigg| \frac{1}{M}\sum_{m=1}^Mg^{\rm{tr}}_h(\bs{y},\bs{x}^{(m)})-\bar{g}^{\rm{tr}}(\bs{y})\Bigg|\right)\Bigg|\int_0^1(1-s)\di{}s,\nonumber\\
={}&\frac{1}{2}\Bigg|f''\left(\bar{g}^{\rm{tr}}(\bs{y})-\Bigg| \frac{1}{M}\sum_{m=1}^Mg^{\rm{tr}}_h(\bs{y},\bs{x}^{(m)})-\bar{g}^{\rm{tr}}(\bs{y})\Bigg|\right)\Bigg|.
\end{align}
From the Koksma--Hlawka inequality~\eqref{eq:QMC.bound} and Assumption~\ref{asu:inverse.FEM} it follows that
\begin{align}
\frac{1}{2}{}&\Bigg|f''\left(\bar{g}^{\rm{tr}}(\bs{y})-\Bigg| \frac{1}{M}\sum_{m=1}^Mg^{\rm{tr}}_h(\bs{y},\bs{x}^{(m)})-\bar{g}^{\rm{tr}}(\bs{y})\Bigg|\right)\Bigg|\nonumber\\
={}&\frac{1}{2}\Bigg|f''\left(\bar{g}^{\rm{tr}}(\bs{y})-\Bigg| \frac{1}{M}\sum_{m=1}^Mg^{\rm{tr}}_h(\bs{y},\bs{x}^{(m)})-\bar{g}^{\rm{tr}}_h(\bs{y})+\bar{g}^{\rm{tr}}_h(\bs{y})-\bar{g}^{\rm{tr}}(\bs{y})\Bigg|\right)\Bigg|,\nonumber\\
\leq{}&\frac{1}{2}\Bigg|f''\left(\bar{g}^{\rm{tr}}(\bs{y})-\left(\Bigg| \frac{1}{M}\sum_{m=1}^Mg^{\rm{tr}}_h(\bs{y},\bs{x}^{(m)})-\bar{g}^{\rm{tr}}_h(\bs{y})\Bigg|+\Bigg|\bar{g}^{\rm{tr}}_h(\bs{y})-\bar{g}^{\rm{tr}}(\bs{y})\Bigg|\right)\right)\Bigg|,\nonumber\\
\leq{}&\frac{1}{2}\Bigg|f''\left(\bar{g}^{\rm{tr}}(\bs{y})-\left(C_{\epsilon,d_2}M^{-1+\epsilon}|V_{\rm{HK}}(g^{\rm{tr}}_h(\bs{y},\cdot))|+C_{\rm{s}}h^{\eta_{\rm{s}}}\bar{g}^{\rm{tr}}(\bs{y})\right)\right)\Bigg|.
\end{align}
Next, Lemma~\ref{cor:inverse.g.h} provides that
\begin{align}
\frac{1}{2}{}&\Bigg|f''\left(\bar{g}^{\rm{tr}}(\bs{y})-\left(C_{\epsilon,d_2}M^{-1+\epsilon}|V_{\rm{HK}}(g^{\rm{tr}}_h(\bs{y},\cdot))|+C_{\rm{s}}h^{\eta_{\rm{s}}}\bar{g}^{\rm{tr}}(\bs{y})\right)\right)\Bigg|\nonumber\\
\leq{}&\frac{1}{2}\Bigg|f''\left(\bar{g}^{\rm{tr}}(\bs{y})-\left(C_{\epsilon,d_2}M^{-1+\epsilon}c(TOL)\bar{g}^{\rm{tr}}_h(\bs{y})+C_{\rm{s}}h^{\eta_{\rm{s}}}\bar{g}^{\rm{tr}}(\bs{y})\right)\right)\Bigg|.
\end{align}
Again from Assumption~\ref{asu:inverse.FEM} it follows that
\begin{align}
\frac{1}{2}{}&\Bigg|f''\left(\bar{g}^{\rm{tr}}(\bs{y})-\left(C_{\epsilon,d_2}M^{-1+\epsilon}c(TOL)\bar{g}^{\rm{tr}}_h(\bs{y})+C_{\rm{s}}h^{\eta_{\rm{s}}}\bar{g}^{\rm{tr}}(\bs{y})\right)\right)\Bigg|\nonumber\\
={}&\frac{1}{2}\Bigg|f''\left(\bar{g}^{\rm{tr}}(\bs{y})-\left(C_{\epsilon,d_2}M^{-1+\epsilon}c(TOL)|\bar{g}^{\rm{tr}}_h(\bs{y})-\bar{g}^{\rm{tr}}(\bs{y})+\bar{g}^{\rm{tr}}(\bs{y})|+C_{\rm{s}}h^{\eta_{\rm{s}}}\bar{g}^{\rm{tr}}(\bs{y})\right)\right)\Bigg|,\nonumber\\
\leq{}&\frac{1}{2}\Bigg|f''\left(\bar{g}^{\rm{tr}}(\bs{y})-\left(C_{\epsilon,d_2}M^{-1+\epsilon}c(TOL)(|\bar{g}^{\rm{tr}}_h(\bs{y})-\bar{g}^{\rm{tr}}(\bs{y})|+\bar{g}^{\rm{tr}}(\bs{y}))+C_{\rm{s}}h^{\eta_{\rm{s}}}\bar{g}^{\rm{tr}}(\bs{y})\right)\right)\Bigg|,\nonumber\\
\leq{}&\frac{1}{2}\Bigg|f''\left(\bar{g}^{\rm{tr}}(\bs{y})-\left(C_{\epsilon,d_2}M^{-1+\epsilon}c(TOL)(1+C_{\rm{s}}h^{\eta_{\rm{s}}})\bar{g}^{\rm{tr}}(\bs{y})+C_{\rm{s}}h^{\eta_{\rm{s}}}\bar{g}^{\rm{tr}}(\bs{y})\right)\right)\Bigg|.
\end{align}
Introducing the notation
\begin{equation}
    k(TOL)\coloneqq \frac{C_{\epsilon,d_2}c(TOL)(1+C_{\rm{s}}h^{\eta_{\rm{s}}})}{M^{1-\epsilon}}+C_{\rm{s}}h^{\eta_{\rm{s}}},
\end{equation}
it holds for
\begin{equation}
    h<C_{\rm{s}}^{-\frac{1}{\eta_{\rm{s}}}}
\end{equation}
and
\begin{align}
    M{}&>\left(\frac{C_{\epsilon,d_2}c(TOL)(1+C_{\rm{s}}h^{\eta_{\rm{s}}})}{(1-C_{\rm{s}}h^{\eta_{\rm{s}}})}\right)^{\frac{1}{1-\epsilon}},
\end{align}
that $0\leq k(TOL)<1$. Moreover, for
\begin{equation}
    K(TOL)\coloneqq 1-k(TOL),
\end{equation}
it holds that $0<K(TOL)\leq1$. Thus, we write
\begin{align}
    \frac{1}{2}\Bigg|f''\left(\bar{g}^{\rm{tr}}(\bs{y})-\left(C_{\epsilon,d_2}M^{-1+\epsilon}c(TOL)(1+C_{\rm{s}}h^{\eta_{\rm{s}}})\bar{g}^{\rm{tr}}(\bs{y})+C_{\rm{s}}h^{\eta_{\rm{s}}}\bar{g}^{\rm{tr}}(\bs{y})\right)\right)\Bigg|{}&=\frac{1}{2}\Bigg|f''\left((1-k(TOL))\bar{g}^{\rm{tr}}(\bs{y})\right)\Bigg|,\nonumber\\
{}&=\frac{1}{2}\Bigg|f''\left(K(TOL)\bar{g}^{\rm{tr}}(\bs{y})\right)\Bigg|,
\end{align}
and from Corollary~\ref{cor:bound.nested.integrand}, it follows that there exists $0<\tilde{b}<\infty$ such that
\begin{equation}\label{eq:monotonicity.0}
    \frac{1}{2}\Bigg|f''\left(K(TOL)\bar{g}^{\rm{tr}}(\bs{y})\right)\Bigg|\leq \frac{1}{2}\tilde{b}\min (y_2,1-y_2)^{-\tilde{A}_2}
\end{equation}
for any $\tilde{A}_2>0$. 
Thus, the bound for the second term in~\eqref{eq:f.h.M.0.split} follows as
 \begin{equation}\label{eq:bound.mu}
    \left|\mu(\bs{y})\right| \leq \frac{1}{2}\tilde{b}\min (y_2,1-y_2)^{-\tilde{A}_2}.
 \end{equation}
It follows that for $h<C_{\rm{s}}^{-1/\eta_{\rm{s}}}$ and
\begin{align}
    M{}&>\left(\frac{C_{\epsilon,d_2}c(TOL)(1+C_{\rm{s}}h^{\eta_{\rm{s}}})}{(1-C_{\rm{s}}h^{\eta_{\rm{s}}})}\right)^{\frac{1}{1-\epsilon}},
\end{align}
 there exists $0<b<\infty$ such that
\begin{equation}
    \left|\varphi_{h,M}(y_1,y_2)\right|\leq \left(C_{\epsilon,d_2}^2M^{-2+2\epsilon}+C_{\eta_{{\rm{s}}}}^2h^{2\eta_{\rm{s}}}\right)b\min(y_1, 1-y_1)^{-A_1}\min(y_2, 1-y_2)^{-A_2}
\end{equation}
for any $A_1, A_2>0$, and thus the verification of Condition~\eqref{eq:f.h.M.0}.
For Condition~\eqref{eq:f.h.M.1}, it follows that
\begin{align}\label{eq:cond.y1}
    \left|\frac{\partial}{\partial y_1}\varphi_{h,M}(\bs{y})\right|{}&=\left|\frac{\partial}{\partial y_1}\left(\left(\lambda(\bs{y})\right)^2\mu(\bs{y})\right)\right|,\nonumber\\
    {}&\leq \left|\frac{\partial}{\partial y_1}\left(\lambda(\bs{y})\right)^2\right|\left|\mu(\bs{y})\right|+\left(\lambda(\bs{y})\right)^2\left|\frac{\partial}{\partial y_1}\mu(\bs{y})\right|.
\end{align}
The bound on $(\lambda(\bs{y}))^2$ follows from~\eqref{eq:bound.lambda} and the bound on $\left|\mu(\bs{y})\right|$ follows from~\eqref{eq:bound.mu}. Next, it follows that 
\begin{align}\label{eq:bound.d.dy1}
    \left|\frac{\partial}{\partial y_1}\left(\lambda(\bs{y})\right)^2\right|=2\left|\frac{1}{M}\sum_{m=1}^Mg^{\rm{tr}}_h(\bs{y},\bs{x}^{(m)})-\bar{g}^{\rm{tr}}(\bs{y})\right|\left|\frac{1}{M}\sum_{m=1}^M\frac{\partial}{\partial y_1}g^{\rm{tr}}_h(\bs{y},\bs{x}^{(m)})-\frac{\partial}{\partial y_1}\bar{g}^{\rm{tr}}(\bs{y})\right|.
\end{align}
The Bound on the first term in~\eqref{eq:bound.d.dy1} follows by a similar derivation as the bound in~\eqref{eq:bound.lambda} as
 \begin{equation}
    2\left|\frac{1}{M}\sum_{m=1}^Mg^{\rm{tr}}_h(\bs{y},\bs{x}^{(m)})-\bar{g}^{\rm{tr}}(\bs{y})\right|\leq \left(C_{\epsilon,d_2}M^{-1+\epsilon}+C_{\eta_{{\rm{s}}}}h^{\eta_{\rm{s}}}\right)2\tilde{b}\min (y_1,1-y_1)^{-\tilde{A}_1}\min (y_2,1-y_2)^{-\tilde{A}_2}.
 \end{equation}
 A similar bound applies to the second term in~\eqref{eq:bound.d.dy1}, 
 \begin{equation}
    \left|\frac{1}{M}\sum_{m=1}^M\frac{\partial}{\partial y_1}g^{\rm{tr}}_h(\bs{y},\bs{x}^{(m)})-\frac{\partial}{\partial y_1}\bar{g}^{\rm{tr}}(\bs{y})\right|\leq \left(C_{\epsilon,d_2}M^{-1+\epsilon}+C_{\eta_{{\rm{s},1}}}h^{\eta_{\rm{s},1}}\right)\tilde{b}\min (y_1,1-y_1)^{-\tilde{A}_1-1}\min (y_2,1-y_2)^{-\tilde{A}_2},
 \end{equation}
except for an extra factor $\min (y_1,1-y_1)$ and a possibly different strong rate $\eta_{\rm{s},1}\leq\eta_{\rm{s}}$ and associated constant $C_{\eta_{{\rm{s},1}}}$ because of the partial derivative. It follows from the inequality of arithmetic and geometric means that
\begin{equation}\label{eq:bound.d.y1}
    \left|\frac{\partial}{\partial y_1}\left(\lambda(\bs{y})\right)^2\right|\leq \left(C_{\epsilon,d_2}^2M^{-2+2\epsilon}+C_{\eta_{{\rm{s},1}}}^2h^{2\eta_{\rm{s},1}}\right)4\tilde{b}^2\min (y_1,1-y_1)^{-2\tilde{A}_1-1}\min (y_2,1-y_2)^{-2\tilde{A}_2}.
\end{equation}
Moreover, for the remaining term in~\eqref{eq:cond.y1}, it follows from the Leibniz integral rule that 
\begin{align}
    \left|\frac{\partial}{\partial y_1}\mu(\bs{y})\right|{}&=\left|\frac{\partial}{\partial y_1}\int_0^1f''\left(\bar{g}^{\rm{tr}}(\bs{y})+s \left(\frac{1}{M}\sum_{m=1}^Mg^{\rm{tr}}_h(\bs{y},\bs{x}^{(m)})-\bar{g}^{\rm{tr}}(\bs{y})\right)\right)(1-s)\di{}s\right|\nonumber\\
    {}&= \left|\int_0^1\frac{\partial}{\partial y_1}f''\left(\bar{g}^{\rm{tr}}(\bs{y})+s \left(\frac{1}{M}\sum_{m=1}^Mg^{\rm{tr}}_h(\bs{y},\bs{x}^{(m)})-\bar{g}^{\rm{tr}}(\bs{y})\right)\right)(1-s)\di{}s\right|,\nonumber\\
    {}&= \left|\int_0^1f'''\left(\bar{g}^{\rm{tr}}(\bs{y})+s \left(\frac{1}{M}\sum_{m=1}^Mg^{\rm{tr}}_h(\bs{y},\bs{x}^{(m)})-\bar{g}^{\rm{tr}}(\bs{y})\right)\right)\right.\nonumber\\
    {}&\left.\quad\quad\times\left(\frac{\partial}{\partial y_1}\bar{g}^{\rm{tr}}(\bs{y})+s \left(\frac{1}{M}\sum_{m=1}^M\frac{\partial}{\partial y_1}g^{\rm{tr}}_h(\bs{y},\bs{x}^{(m)})-\frac{\partial}{\partial y_1}\bar{g}^{\rm{tr}}(\bs{y})\right)\right)(1-s)\di{}s\right|.
\end{align}
A similar derivation as for the bound in~\eqref{eq:bound.mu} yields that for 
\begin{equation}
    h<C_{\rm{s},1}^{-\frac{1}{\eta_{\rm{s},1}}}
\end{equation}
and
\begin{align}
    M{}&>\left(\frac{C_{\epsilon,d_2}c(TOL)(1+C_{\rm{s},1}h^{\eta_{\rm{s},1}})}{(1-C_{\rm{s},1}h^{\eta_{\rm{s},1}})}\right)^{\frac{1}{1-\epsilon}},
\end{align}
there exists $0<\tilde{b}<\infty$ such that
\begin{align}
    {}& \left|\int_0^1f'''\left(\bar{g}^{\rm{tr}}(\bs{y})+s \left(\frac{1}{M}\sum_{m=1}^Mg^{\rm{tr}}_h(\bs{y},\bs{x}^{(m)})-\bar{g}^{\rm{tr}}(\bs{y})\right)\right)\right.\nonumber\\
    {}&\left.\quad\quad\times\left(\frac{\partial}{\partial y_1}\bar{g}^{\rm{tr}}(\bs{y})+s \left(\frac{1}{M}\sum_{m=1}^M\frac{\partial}{\partial y_1}g^{\rm{tr}}_h(\bs{y},\bs{x}^{(m)})-\frac{\partial}{\partial y_1}\bar{g}^{\rm{tr}}(\bs{y})\right)\right)(1-s)\di{}s\right|,\nonumber\\
    {}&\leq \frac{1}{2}\left|f'''\left(K(TOL)|\bar{g}^{\rm{tr}}(\bs{y})|\right)\right|\left|K(TOL)\frac{\partial}{\partial y_1}\bar{g}^{\rm{tr}}(\bs{y})\right|,\nonumber\\
    {}&\leq\frac{1}{2}\tilde{b}^2\min(y_1,1-y_1)^{-\tilde{A}_1-1}\min (y_2,1-y_2)^{-\tilde{A}_2}
\end{align}
for any $\tilde{A}_1,\tilde{A}_2>0$. Combining the above results yields that there exists $0<b<\infty$ such that
\begin{equation}
    \left|\frac{\partial}{\partial y_1}\varphi_{h,M}(y_1,y_2)\right|\leq \left(C_{\epsilon,d_2}^2M^{-2+2\epsilon}+C_{\eta_{{\rm{s}},1}}^2h^{2\eta_{{\rm{s}},1}}\right)b\min (y_1,1-y_1)^{-A_1-1}\min (y_2,1-y_2)^{-A_2}
\end{equation}
for any $A_1, A_2>0$. A similar bound applies to the Condition~\eqref{eq:f.h.M.2}, that is,
\begin{equation}
    \left|\frac{\partial}{\partial y_2}\varphi_{h,M}(y_1,y_2)\right|\leq \left(C_{\epsilon,d_2}^2M^{-2+2\epsilon}+C_{\eta_{{\rm{s}},2}}^2h^{2\eta_{{\rm{s}},2}}\right)b\min (y_1,1-y_1)^{-A_1}\min (y_2,1-y_2)^{-A_2-1}.
\end{equation}
Finally, for Condition~\eqref{eq:f.h.M.3}, it follows that
\begin{equation}
    \frac{\partial^2}{\partial y_1\partial y_2}\varphi_{h,M}(\bs{y})=\frac{\partial^2}{\partial y_1\partial y_2}\left(\lambda(\bs{y})\right)^2\mu(\bs{y})+\frac{\partial}{\partial y_1}\left(\lambda(\bs{y})\right)^2\frac{\partial}{\partial y_2}\mu(\bs{y})+\frac{\partial}{\partial y_2}\left(\lambda(\bs{y})\right)^2\frac{\partial}{\partial y_1}\mu(\bs{y})+\left(\lambda(\bs{y})\right)^2\frac{\partial^2}{\partial y_1\partial y_2}\mu(\bs{y}).
\end{equation}
The following bound applies:
\begin{align}\label{eq:bound.d.y1.y2}
    \left|\frac{\partial^2}{\partial y_1\partial y_2}\left(\lambda(\bs{y})\right)^2\right|{}&=\left|\frac{\partial^2}{\partial y_1\partial y_2}\left(\frac{1}{M}\sum_{m=1}^Mg^{\rm{tr}}_h(\bs{y},\bs{x}^{(m)})-\bar{g}^{\rm{tr}}(\bs{y})\right)^2\right|\nonumber\\
    {}&=2\left|\frac{\partial}{\partial y_1}\left(\left(\frac{1}{M}\sum_{m=1}^Mg^{\rm{tr}}_h(\bs{y},\bs{x}^{(m)})-\bar{g}^{\rm{tr}}(\bs{y})\right)\left(\frac{1}{M}\sum_{m=1}^M\frac{\partial}{\partial y_2}g^{\rm{tr}}_h(\bs{y},\bs{x}^{(m)})-\frac{\partial}{\partial y_2}\bar{g}^{\rm{tr}}(\bs{y})\right)\right)\right|,\nonumber\\
    {}&\leq2\left|\left(\frac{1}{M}\sum_{m=1}^M\frac{\partial}{\partial y_1}g^{\rm{tr}}_h(\bs{y},\bs{x}^{(m)})-\frac{\partial}{\partial y_1}\bar{g}^{\rm{tr}}(\bs{y})\right)\right|\left|\left(\frac{1}{M}\sum_{m=1}^M\frac{\partial}{\partial y_2}g^{\rm{tr}}_h(\bs{y},\bs{x}^{(m)})-\frac{\partial}{\partial y_2}\bar{g}^{\rm{tr}}(\bs{y})\right)\right|,\nonumber\\
    {}&\quad+2\left|\left(\frac{1}{M}\sum_{m=1}^Mg^{\rm{tr}}_h(\bs{y},\bs{x}^{(m)})-\bar{g}^{\rm{tr}}(\bs{y})\right)\right|\left|\left(\frac{1}{M}\sum_{m=1}^M\frac{\partial^2}{\partial y_1\partial y_2}g^{\rm{tr}}_h(\bs{y},\bs{x}^{(m)})-\frac{\partial^2}{\partial y_1\partial y_2}\bar{g}^{\rm{tr}}(\bs{y})\right)\right|.
\end{align}
The bound on the first term in \eqref{eq:bound.d.y1.y2} follows by a derivation similar to that in \eqref{eq:bound.d.y1}:
\begin{align}
    {}&2\left|\left(\frac{1}{M}\sum_{m=1}^M\frac{\partial}{\partial y_1}g^{\rm{tr}}_h(\bs{y},\bs{x}^{(m)})-\frac{\partial}{\partial y_1}\bar{g}^{\rm{tr}}(\bs{y})\right)\right|\left|\left(\frac{1}{M}\sum_{m=1}^M\frac{\partial}{\partial y_2}g^{\rm{tr}}_h(\bs{y},\bs{x}^{(m)})-\frac{\partial}{\partial y_2}\bar{g}^{\rm{tr}}(\bs{y})\right)\right|\nonumber\\
    {}&\quad\leq 4\left(C_{\epsilon,d_2}^2M^{-2+2\epsilon}+C_{\eta_{{\rm{s}},d_1}}^2h^{2\eta_{{\rm{s}},d_1}}\right)\tilde{b}^2\min(y_1, 1-y_1)^{-2\tilde{A}_1-1}\min(y_2, 1-y_2)^{-2\tilde{A}_2-1}
\end{align}
for 
\begin{equation}
    h<C_{{\rm{s}},d_1}^{-\frac{1}{\eta_{{\rm{s}},d_1}}}
\end{equation}
and
\begin{align}\label{eq:M.TOL}
    M{}&>\left(\frac{C_{\epsilon,d_2}c(TOL)(1+C_{{\rm{s}},d_1}h^{\eta_{{\rm{s}},d_1}})}{(1-C_{{\rm{s}},d_1}h^{\eta_{{\rm{s}},d_1}})}\right)^{\frac{1}{1-\epsilon}},\nonumber\\
    {}&=\tilde{M}(TOL).
\end{align}
 The same bound also applies to the second term in \eqref{eq:bound.d.y1.y2}.  Moreover, it follows from the Leibniz integral rule that
\begin{align}
    \left|\frac{\partial^2}{\partial y_1\partial y_2}\mu(\bs{y})\right|{}&=\left|\frac{\partial^2}{\partial y_1\partial y_2}\int_0^1f''\left(\bar{g}^{\rm{tr}}(\bs{y})+s \left(\frac{1}{M}\sum_{m=1}^Mg^{\rm{tr}}_h(\bs{y},\bs{x}^{(m)})-\bar{g}^{\rm{tr}}(\bs{y})\right)\right)(1-s)\di{}s\right|\nonumber\\
    {}&\quad\leq \frac{1}{2}f^{''''}\left(K(TOL)\bar{g}^{\rm{tr}}(\bs{y})\right)\left|K(TOL)\frac{\partial}{\partial y_1}\bar{g}^{\rm{tr}}(\bs{y})\right|\left|K(TOL)\frac{\partial}{\partial y_2}\bar{g}^{\rm{tr}}(\bs{y})\right|\nonumber\\
    {}&\quad\quad+\frac{1}{2}f'''\left(K(TOL)\bar{g}^{\rm{tr}}(\bs{y})\right)\left|K(TOL)\frac{\partial^2}{\partial y_1\partial y_2}\bar{g}^{\rm{tr}}(\bs{y})\right|,\nonumber\\
    {}&\quad\leq \tilde{b}^2\min(y_1,1-y_1)^{-2\tilde{A}_1-1}\min (y_2,1-y_2)^{-2\tilde{A}_2-1}.
\end{align}
Combining the previous results, it is demonstrated that for $h<C_{{\rm{s}},d_1}^{-1/\eta_{{\rm{s}},d_1}}$ and $M>\tilde{M}(TOL)$, there exists $0<b<\infty$ such that
\begin{equation}
    \left|\frac{\partial^2}{\partial y_1\partial y_2}\varphi_{h,M}(y_1,y_2)\right|\leq \left(C_{\epsilon,d_2}^2M^{-2+2\epsilon}+C_{\eta_{{\rm{s}},d_1}}^2h^{2\eta_{{\rm{s}},d_1}}\right)b\min (y_1,1-y_1)^{-A_1-1}\min (y_2,1-y_2)^{-A_2-1}
\end{equation}
for any $A_1, A_2>0$. The general case for an arbitrary dimension is a consequence of Fa\`{a} di Bruno's formula.
\end{proof}


\end{document}